\documentclass[a4paper,12pt]{article}
\usepackage{amsmath}
\usepackage{amssymb}

\newtheorem{theorem}{Theorem}[section]
\newtheorem{acknowledgement}[theorem]{Acknowledgement}

\newtheorem{corollary}[theorem]{Corollary}

\newtheorem{definition}[theorem]{Definition}

\newtheorem{lemma}[theorem]{Lemma}

\newtheorem{proposition}[theorem]{Proposition}
\newtheorem{remark}[theorem]{Remark}

\newenvironment{proof}[1][Proof]{\noindent\textbf{#1.} }{\ \rule{0.5em}{0.5em}}
\input{tcilatex}
\begin{document}

\title{Derivations in algebras of operator-valued functions}
\author{A.F. Ber, B. de Pagter and F.A. Sukochev\thanks{%
Research is partially supported by the Australian Research Council}}
\date{}
\maketitle

\begin{abstract}
In this paper we study derivations in subalgebras of $L_{0}^{wo}\left( \nu ;%
\mathcal{L}\left( X\right) \right) $, the algebra of all weak operator
measurable funtions $f:S\rightarrow \mathcal{L}\left( X\right) $, where $%
\mathcal{L}\left( X\right) $ is the Banach algebra of all bounded linear
operators on a Banach space $X$. It is shown, in particular, that all
derivations on $L_{0}^{wo}\left( \nu ;\mathcal{L}\left( X\right) \right) $
are inner whenever $X$ is separable and infinite dimensional. This contrasts
strongly with the fact that $L_{0}^{wo}\left( \nu ;\mathcal{L}\left(
X\right) \right) $ admits non-trivial non-inner derivations whenever $X$ is
finite dimensional and the measure $\nu $ is non-atomic. As an application
of our approach, we study derivations in various algebras of measurable
operators affiliated with von Neumann algebras.
\end{abstract}

\section{Introduction}

It follows from recent results (\cite{BCS1}, \cite{BCS2}) that the algebra $%
L_{0}\left( \nu \right) $ of all (equivalence classes of) $\nu $-measurable
functions on a measure space $\left( S,\Sigma ,\nu \right) $ does not admit
non-trivial derivations if and only if the Boolean algebra of all $\nu $%
-measurable subsets of $S$ is atomic. Thus, if $\nu $ is atomless, then the
algebra $L_{0}\left( \nu \right) $ admits non-trivial derivations $\delta $
which are obviously not inner and also not continuous with respect to the
topology of convergence in measure. If $X$ is a finite dimensional Banach
space, then it is straightforward that any such derivation extends
"coordinatewise" to a non-trivial non-inner derivation on the algebra $%
L_{0}\left( \nu ;\mathcal{L}\left( X\right) \right) $ of all (Bochner) $\nu $%
-measurable functions $f:S\rightarrow \mathcal{L}\left( X\right) $. However,
if $X$ is an infinite dimensional Banach space, then the situation changes
drastically.

Suppose now that $X$ is an infinite dimensional separable Banach space and
let $L_{0}^{wo}\left( \nu ;\mathcal{L}\left( X\right) \right) $ be the
algebra of all weak operator $\nu $-measurable functions $f:S\rightarrow 
\mathcal{L}\left( X\right) $ (see Section \ref{SectVecVal} for precise
definitions). It follows from one of the main results of the present paper
(see Theorem \ref{PThm01}) that every derivation in $L_{0}^{wo}\left( \nu ;%
\mathcal{L}\left( X\right) \right) $ is necessarily inner and hence
continuous with respect to the topology of convergence in measure. Actually,
Theorem \ref{PThm01} shows that similar results hold for a large class of
appropriate subalgebras (so-called admissible subalgebras) of $%
L_{0}^{wo}\left( \nu ;\mathcal{L}\left( X\right) \right) $. It should be
pointed out that, specializing to the case that the underlying set $S$
consists of one point, we recover the classical result of P. Chernoff (\cite%
{Cher}), who showed that every derivation on $\mathcal{L}\left( X\right) $
is inner (answering a question raised by S. Sakai in \cite{Sak2}, Section
4.1).

Interesting applications of our approach and techniques lie in the realm of
the theory of non-commutative integration on semi-finite von Neumann
algebras, initiated by I.E. Segal (\cite{Se}). The classical spaces of
measurable operators associated with an arbitrary pair $\left( \mathcal{M}%
,\tau \right) $ of a semi-finite von Neumann algebra $\mathcal{M}$ and a
faithful normal semi-finite trace $\tau $ are the following:

\begin{enumerate}
\item[(i).] the space of all measurable operators $S\left( \mathcal{M}%
\right) $, introduced by Segal (\cite{Se});

\item[(ii).] the space of all locally measurable operators $LS\left( 
\mathcal{M}\right) $, introduced by S. Sankaran (\cite{San}; see also F.J.
Yeadon's paper \cite{Ye});

\item[(iii).] the space $S\left( \tau \right) =S\left( \mathcal{M},\tau
\right) $ of all $\tau $-measurable operators, first considered by E. Nelson
(\cite{Ne}).
\end{enumerate}

\noindent It should be noted that always $S\left( \tau \right) \subseteq
S\left( \mathcal{M}\right) \subseteq LS\left( \mathcal{M}\right) $ (some of
the relevant definitions may be found at the beginning of Section \ref%
{SectMeas}). For some of the simplest examples of von\ Neumann algebras of
type $I$, we have 
\begin{equation}
S\left( \tau \right) =S\left( \mathcal{M}\right) =LS\left( \mathcal{M}%
\right) .  \label{Peq11}
\end{equation}%
For example, the equalities (\ref{Peq11}) hold in the special case $\left( 
\mathcal{M},\tau \right) =\left( L_{\infty }\left( \nu \right) ,\nu \right) $%
, where the trace $\nu $ is given by integration with respect to a finite
measure $\nu $. Similarly, if $H$ is a Hilbert space and $\mathcal{M}%
=B\left( H\right) $ is the von Neumann algebra of all bounded linear
operators in $H$ (so, $B\left( H\right) =\mathcal{L}\left( X\right) $ as
Banach algebras, with $X=H$), equipped with the standard trace $\tau =\func{%
tr}$, then the equalities (\ref{Peq11}) hold (and these algebras coincide
with the von Neumann algebra $B\left( H\right) $ itself). However, these
equalities do not hold any longer for the von Neumann algebra $\mathcal{M}%
=L_{\infty }\left( \nu \right) \overline{\otimes }B\left( H\right) $,
equipped with the tensor product trace $\tau =\nu \otimes \func{tr}$,
whenever $\nu $ is non-atomic and $H$ is infinite dimensional (see Remark %
\ref{PRem02}). As will be seen in Section \ref{SectMeas}, in this latter
case, our general results for derivations on subalgebras of $%
L_{0}^{wo}\left( \nu ;\mathcal{L}\left( X\right) \right) $ have implications
for each of the algebras listed above.

Firstly, it will be shown that the algebra $L_{0}^{wo}\left( \nu ;B\left(
H\right) \right) $ may be identified with the algebra $LS\left( \mathcal{M}%
\right) $ of all locally measurable operators affiliated with $\mathcal{M}%
=L_{\infty }\left( \nu \right) \overline{\otimes }B\left( H\right) $ (see
Proposition \ref{PProp08}) and therefore, every derivation on $LS\left( 
\mathcal{M}\right) $ is inner. This result has already been obtained in \cite%
{AAK}, through methods based on the theory of Kaplansky-Hilbert modules. Our
approach is different and we believe that the methods based on the
realization of the algebra $LS\left( \mathcal{M}\right) $ as a space of
operator-valued functions provide additional insight into the matter. For
example, the following observation, which seems to have escaped the
attention, exhibits itself: topology of convergence locally in measure in $%
LS\left( \mathcal{M}\right) $, as introduced by Yeadon in \cite{Ye},
coincides with the natural topology of convergence in measure on sets of
finite measure in $L_{0}^{wo}\left( \nu ;B\left( H\right) \right) $ (see
Remark \ref{PRem03}).

The effectiveness of our approach is illustrated further by showing that the
results of Section \ref{SectAdmis} may be applied to show that derivations
on the algebras $S\left( \mathcal{M}\right) $ and $S\left( \tau \right) $
are also inner (see Corollary \ref{PCor03}). These applications are not
straightforward and require some additional techniques (see Section \ref%
{SectMeas}, in particular, Proposition \ref{AProp02} and its proof), which
may be of independent potential interest. Furthermore, the notion of
admissible subalgebra of $L_{0}^{wo}\left( \nu ;\mathcal{L}\left( X\right)
\right) $ (see Section \ref{SectAdmis}) is flexible enough to include also
the ideal $S_{0}\left( \tau \right) $ in $S\left( \tau \right) $ of all $%
\tau $-compact operators (see Proposition \ref{PProp07}). Consequently, also
derivations on $S_{0}\left( \tau \right) $ are discussed (see Corollary \ref%
{PCor02}).

In Section \ref{SectNotTerm} some basic notation and terminology is
introduced and in Section \ref{SectVecVal} some necessary results concerning
vector-valued and operator valued functions are established. Section \ref%
{SectDer} exhibits some elementary algebraic properties of derivations,
which will be used in the subsequent sections. In Section \ref{SectAdmis},
which contains one of the main results of the paper (Theorem \ref{PThm01}),
we study in detail derivations on admissible subalgebras of the algebra $%
L_{0}^{wo}\left( \nu ;\mathcal{L}\left( X\right) \right) $. Finally, in the
last section, the results of Section \ref{SectAdmis} will be applied to
derivations on algebras of measurable operators affiliated with a von
Neumann algebra $\mathcal{M}=L_{\infty }\left( \nu \right) \overline{\otimes 
}B\left( H\right) $.

\begin{acknowledgement}
The authors thank Sh.A. Ayupov for the preprint \cite{AAK} and V.I. Chilin
for his continuing interest. The third author would like to thank the Delft
Institute of Applied Mathematics for its hospitality and the Thomas
Stieltjes Institute for Mathematics for its support during his stay in Delft
as a Visiting Stieltjes Professor in June and July 2008.
\end{acknowledgement}

\section{Notation and terminology\label{SectNotTerm}}

If $X$ is a (complex) Banach space, then $X^{\ast }$ denotes its Banach dual
space. The value of a functional $x^{\ast }\in X^{\ast }$ at $x\in X$ is
denoted by $\left\langle x,x^{\ast }\right\rangle $. The Banach algebra of
all bounded linear operators in $X$ is denoted by $\mathcal{L}\left(
X\right) $ and its multiplicative unit is the identity operator $\mathbf{1}=%
\mathbf{1}_{X}$. The subalgebra of $\mathcal{L}\left( X\right) $ consisting
of all finite rank operators is denoted by $\mathcal{F}\left( X\right) $.
Furthermore, $\mathcal{K}\left( X\right) $ is the Banach algebra of all
compact operators. Evidently, $\mathcal{F}\left( X\right) \subseteq \mathcal{%
K}\left( X\right) \subseteq \mathcal{L}\left( X\right) $. If $p\in \mathcal{L%
}\left( X\right) $ is a projection (that is, $p^{2}=p$), then its range $%
\limfunc{Ran}\left( p\right) $ will also be denoted by $X_{p}$. Hence, $%
X=X_{p}\oplus X_{p^{\bot }}$, where $p^{\bot }=\mathbf{1}-p$. The algebra $%
\mathcal{L}\left( X_{p}\right) $ may be identified with a subalgebra of $%
\mathcal{L}\left( X\right) $ as follows. Given $T\in \mathcal{L}\left(
X_{p}\right) $, define $\tilde{T}\in \mathcal{L}\left( X\right) $ by setting 
$\tilde{T}x=Tx$, $x\in X_{p}$, and $\tilde{T}x=0$, $x\in X_{p^{\bot }}$.
Note that 
\begin{equation*}
\left\Vert T\right\Vert _{\mathcal{L}\left( X_{p}\right) }\leq \left\Vert 
\tilde{T}\right\Vert _{\mathcal{L}\left( X\right) }\leq \left\Vert
p\right\Vert _{\mathcal{L}\left( X\right) }\left\Vert T\right\Vert _{%
\mathcal{L}\left( X_{p}\right) }.
\end{equation*}%
The map $T\longmapsto \tilde{T}$ is an algebraic isomorphism from $\mathcal{L%
}\left( X_{p}\right) $ onto the subalgebra 
\begin{equation*}
p\mathcal{L}\left( X\right) p=\left\{ pSp:S\in \mathcal{L}\left( X\right)
\right\} =\left\{ S\in \mathcal{L}\left( X\right) :S=pSp\right\}
\end{equation*}%
of $\mathcal{L}\left( X\right) $. It should be observed that the
identification of the algebra $\mathcal{L}\left( X_{p}\right) $ with a
subalgebra of $\mathcal{L}\left( X\right) $ depends on the choice of the
specific projection $p$ onto $X_{p}$.

\section{Spaces of vector-valued functions\label{SectVecVal}}

Next, we recall some facts concerning vector-valued and operator-valued
functions, and introduce some relevant notation. For the general theory of
vector-valued measurable functions, the reader is referred to \cite{DU}.
Given a non-empty set $S$ and a (complex) vector space $V$, if $%
f:S\rightarrow \mathbb{C}$ and $v\in V$, then the function $f\otimes
v:S\rightarrow V$ is defined by setting $\left( f\otimes v\right) \left(
s\right) =f\left( s\right) v$, $s\in S$. If, in addition, $\Sigma $ is a $%
\sigma $-algebra of subsets of $S$, then $\func{sim}\left( \Sigma ;V\right) $
denotes the vector space of all $\Sigma $-simple $V$-valued functions, that
is, $f\in \func{sim}\left( \Sigma ;V\right) $ if and only if $f$ is of the
form $f=\sum_{j=1}^{n}\chi _{A_{j}}\otimes v_{j}$, with $A_{j}\in \Sigma $
and $v_{j}\in V$ ($j=1,\ldots ,n$; $n\in \mathbb{N}$). Note that, without
loss of generality, it may be assumed that $A_{i}\cap A_{j}=\emptyset $
whenever $i\neq j$.

Let $Y$ be a (complex) Banach space and $\left( S,\Sigma ,\nu \right) $ be a
complete $\sigma $-finite measure space (so, $\Sigma $ equals the $\sigma $%
-algebra of all $\nu $-measurable sets). Assume that $\Gamma $ is a linear
subspace of $Y^{\ast }$ which contains a countable norming subset for $Y$
(that is, there exists a sequence $\left\{ y_{n}^{\ast }\right\}
_{n=1}^{\infty }$ in $\Gamma $ such that $\left\Vert y\right\Vert
_{Y}=\sup_{n}\left\vert \left\langle y,y_{n}^{\ast }\right\rangle
\right\vert $ for all $y\in Y$; note that this implies that $\left\Vert
y_{n}^{\ast }\right\Vert _{Y^{\ast }}\leq 1$ for all $n$). A function $%
f:S\rightarrow Y$ is called $\nu $\textit{-measurable with respect to }$%
\Gamma $ if the scalar functions $s\longmapsto \left\langle f\left( s\right)
,y^{\ast }\right\rangle $, $s\in S$, are $\nu $-measurable for all $y^{\ast
}\in \Gamma $ (see e.g. \cite{DU}, Section II.1). The linear space of all
such $Y$-valued measurable functions is denoted by $\mathcal{L}_{0}^{\Gamma
}\left( \nu ;Y\right) $ and $\mathcal{L}_{\infty }^{\Gamma }\left( \nu
;Y\right) $ is the subspace of $\mathcal{L}_{0}^{\Gamma }\left( \nu
;Y\right) $ consisting of all $\nu $-essentially bounded functions. The
following simple observation will be used.

\begin{lemma}
\label{NLem04}If $f\in \mathcal{L}_{0}^{\Gamma }\left( \nu ;Y\right) $, then
the function $s\longmapsto \left\Vert f\left( s\right) \right\Vert _{Y}$, $%
s\in S$, is $\nu $-measurable.
\end{lemma}

\begin{proof}
If $\left\{ y_{n}^{\ast }\right\} _{n=1}^{\infty }\subseteq \Gamma $ is a
norming subset, then 
\begin{equation*}
\left\Vert f\left( s\right) \right\Vert _{Y}=\sup_{n}\left\vert \left\langle
f\left( s\right) ,y_{n}^{\ast }\right\rangle \right\vert ,\ \ \ s\in S.
\end{equation*}%
Since the functions $s\longmapsto \left\vert \left\langle f\left( s\right)
,y_{n}^{\ast }\right\rangle \right\vert $, $s\in S$, are $\nu $-measurable
for all $n$, the statement of the lemma is now clear.\medskip
\end{proof}

A function $f\in \mathcal{L}_{0}^{\Gamma }\left( \nu ;Y\right) $ is called a 
$\nu $\textit{-null function} if $f\left( s\right) =0$ $\nu $-a.e. on $S$.
Using the assumption on $\Gamma $, the following lemma is easily proved.

\begin{lemma}
\label{NLem03}If $f\in \mathcal{L}_{0}^{\Gamma }\left( \nu ;Y\right) $, then
the following two statements are equivalent:

\begin{enumerate}
\item[(i).] $f$ is a $\nu $-null function;

\item[(ii).] for every $y^{\ast }\in \Gamma $ the scalar function $%
s\longmapsto \left\langle f\left( s\right) ,y^{\ast }\right\rangle $, $s\in
S $, is a $\nu $-null function.
\end{enumerate}
\end{lemma}

The space of all $\nu $-null functions is denoted by $\mathcal{N}\left( \nu
;Y\right) $. Note that $\mathcal{N}\left( \nu ;Y\right) $ is a linear
subspace of both $\mathcal{L}_{0}^{\Gamma }\left( \nu ;Y\right) $ and $%
\mathcal{L}_{\infty }^{\Gamma }\left( \nu ;Y\right) $. The spaces $%
L_{0}^{\Gamma }\left( \nu ;Y\right) $ and $L_{\infty }^{\Gamma }\left( \nu
;Y\right) $ are defined by setting 
\begin{equation*}
L_{0}^{\Gamma }\left( \nu ;Y\right) =\mathcal{L}_{0}^{\Gamma }\left( \nu
;Y\right) \diagup \mathcal{N}\left( \nu ;Y\right) ,\ \ \ L_{\infty }^{\Gamma
}\left( \nu ;Y\right) =\mathcal{L}_{\infty }^{\Gamma }\left( \nu ;Y\right)
\diagup \mathcal{N}\left( \nu ;Y\right) ,
\end{equation*}%
respectively. Equipped with the norm $\left\Vert \cdot \right\Vert _{\infty
} $, given by 
\begin{equation*}
\left\Vert f\right\Vert _{\infty }=\func{esssup}_{s\in S}\left\Vert f\left(
s\right) \right\Vert _{Y},\ \ \ f\in L_{\infty }^{\Gamma }\left( \nu
;Y\right) ,
\end{equation*}%
it is easily verified that $L_{\infty }^{\Gamma }\left( \nu ;Y\right) $ is a
Banach space. Note that the algebraic tensor products $L_{0}\left( \nu
\right) \otimes Y$ and $L_{\infty }\left( \nu \right) \otimes Y$ may be
identified with linear subspaces of $L_{0}^{\Gamma }\left( \nu ;Y\right) $
and $L_{\infty }^{\Gamma }\left( \nu ;Y\right) $, respectively.

\begin{remark}
\label{NRem02}In the present paper we shall be dealing with the following
two special cases of this general setting. Suppose that $X$ is a separable
(complex) Banach space. Let $\left\{ x_{n}\right\} _{n=1}^{\infty }$ be
dense in the closed unit ball of $X$ and suppose that $x_{n}^{\ast }\in
X^{\ast }$ is such that $\left\Vert x_{n}^{\ast }\right\Vert _{X^{\ast }}=1$
and $\left\langle x_{n},x_{n}^{\ast }\right\rangle =\left\Vert
x_{n}\right\Vert _{X}$, for each $n\in \mathbb{N}$.

\begin{enumerate}
\item[(a).] The sequence $\left\{ x_{n}^{\ast }\right\} _{n=1}^{\infty }$ is
norming for $X$. Hence, $\Gamma =X^{\ast }$ satisfies the requirements for $%
Y=X$. It follows from the Pettis measurability criterion (see e.g. \cite{DU}%
, Theorem II.1.2) that the $\nu $-measurable functions with respect to $%
\Gamma $ are precisely the \emph{Bochner }$\nu $\emph{-measurable functions}%
. In this case, the spaces $L_{0}^{X^{\ast }}\left( \nu ;X\right) $ and $%
L_{\infty }^{X^{\ast }}\left( \nu ;X\right) $ will be denoted simply by $%
L_{0}\left( \nu ;X\right) $ and $L_{\infty }\left( \nu ;X\right) $,
respectively.

\item[(b).] Let $\mathcal{L}\left( X\right) $ be the Banach algebra of all
bounded linear operators in $X$. For $x\in X$ and $x^{\ast }\in X^{\ast }$,
define the linear functional $x\otimes x^{\ast }\in \mathcal{L}\left(
X\right) ^{\ast }$ by setting 
\begin{equation*}
\left\langle T,x\otimes x^{\ast }\right\rangle =\left\langle Tx,x^{\ast
}\right\rangle ,\ \ \ T\in \mathcal{L}\left( X\right) ,
\end{equation*}%
and let $\Gamma =X\otimes X^{\ast }$, the linear span of $\left\{ x\otimes
x^{\ast }:x\in X,x^{\ast }\in X^{\ast }\right\} $ in $\mathcal{L}\left(
X\right) ^{\ast }$. It is clear that $\left\{ x_{n}\otimes x_{m}^{\ast
}:n,m\in \mathbb{N}\right\} $ is a countable norming set for $Y=\mathcal{L}%
\left( X\right) $. The $\mathcal{L}\left( X\right) $-valued $\nu $%
-measurable functions with respect to $X\otimes X^{\ast }$ will be called 
\emph{weak operator }$\nu $\emph{-measurable functions} and the space $%
\mathcal{L}_{0}^{X\otimes X^{\ast }}\left( \nu ;\mathcal{L}\left( X\right)
\right) $ is denoted by $\mathcal{L}_{0}^{wo}\left( \nu ;\mathcal{L}\left(
X\right) \right) $. The spaces $L_{0}^{X\otimes X^{\ast }}\left( \nu ;%
\mathcal{L}\left( X\right) \right) $ and $L_{\infty }^{X\otimes X^{\ast
}}\left( \nu ;\mathcal{L}\left( X\right) \right) $ are denoted by $%
L_{0}^{wo}\left( \nu ;\mathcal{L}\left( X\right) \right) $ and $L_{\infty
}^{wo}\left( \nu ;\mathcal{L}\left( X\right) \right) $, respectively.
Similar notation will be used if $\mathcal{L}\left( X\right) $ is replaced
by any closed subspace $\mathfrak{B}$ of $\mathcal{L}\left( X\right) $. It
should be observed that, by the Pettis measurability criterion, a function $%
f:S\rightarrow \mathcal{L}\left( X\right) $ is weak operator $\nu $%
-measurable if and only if for each $x\in X$ the $X$-valued function $%
s\longmapsto f\left( s\right) x$, $s\in S$, is Bochner $\nu $-measurable.
\end{enumerate}
\end{remark}

Returning to the general setting, the space $L_{0}^{\Gamma }\left( \nu
;Y\right) $ is equipped with the vector space topology $\mathcal{T}_{0}$ of 
\textit{convergence in measure on sets of finite measure}. A neighbourhood
base at $0$ for $\mathcal{T}_{0}$ is given by the sets 
\begin{equation*}
V\left( A;\varepsilon ,\delta \right) =\left\{ f\in L_{0}^{\Gamma }\left(
\nu ;Y\right) :\nu \left\{ s\in A:\left\Vert f\left( s\right) \right\Vert
_{Y}>\varepsilon \right\} <\delta \right\} ,
\end{equation*}%
where $\varepsilon ,\delta >0$ and $A$ ranges over all $\nu $-measurable
subsets of $S$ satisfying $\nu \left( A\right) <\infty $. The topology $%
\mathcal{T}_{0}$ is metrizable, induced by the translation invariant metric $%
d_{0}$ given by%
\begin{equation*}
d_{0}\left( f,g\right) =\sum_{n=1}^{\infty }\frac{1}{2^{n}\mu \left(
A_{n}\right) }\int_{A_{n}}\frac{\left\Vert f\left( s\right) -g\left(
s\right) \right\Vert _{Y}}{1+\left\Vert f\left( s\right) -g\left( s\right)
\right\Vert _{Y}}d\nu \left( s\right) ,\ \ f,g\in L_{0}^{wo}\left( \nu
;Y\right) ,
\end{equation*}%
where $\left\{ A_{n}\right\} _{n=1}^{\infty }$ is a partition of $S$ into $%
\nu $-measurable sets satisfying $0<\nu \left( A_{n}\right) <\infty $ for
all $n$ (note: the metric $d_{0}$ depends on the choice of the partition $%
\left\{ A_{n}\right\} _{n=1}^{\infty }$). A sequence $\left\{ f_{n}\right\}
_{n=1}^{\infty }$ in $L_{0}^{\Gamma }\left( \nu ;Y\right) $ converges to $%
f\in L_{0}^{\Gamma }\left( \nu ;Y\right) $ with respect to $\mathcal{T}_{0}$
if and only if every subsequence $\left\{ f_{n_{k}}\right\} _{k=1}^{\infty }$
has a subsequence $\left\{ f_{n_{k_{j}}}\right\} _{j=1}^{\infty }$ such that 
$f_{n_{k_{j}}}\left( s\right) \rightarrow f\left( s\right) $ with respect to
the norm in $Y$ as $j\rightarrow \infty $, $\nu $-a.e. on $S$. This is also
equivalent to 
\begin{equation*}
\lim_{n\rightarrow \infty }\nu \left\{ s\in A:\left\Vert f_{n}\left(
s\right) -f\left( s\right) \right\Vert _{Y}\geq \varepsilon \right\} =0
\end{equation*}%
for every $\varepsilon >0$ and every $\nu $-measurable set $A\subseteq S$
with $\nu \left( A\right) <\infty $. The reader is advised to keep in mind
that, in general, $L_{0}^{\Gamma }\left( \nu ;Y\right) $ is not a Banach
algebra. The following observation will be important.

\begin{lemma}
\label{NLem01}The space $L_{0}^{\Gamma }\left( \nu ;Y\right) $ is complete
with respect to the topology $\mathcal{T}_{0}$.
\end{lemma}

\begin{proof}
For the sake of convenience, we present the proof for the case that $\nu
\left( S\right) <\infty $. The extension to the general $\sigma $-finite
case is easy. If $\nu \left( S\right) <\infty $, then it may be assumed that
the metric $d_{0}$ is given by 
\begin{equation*}
d_{0}\left( f,g\right) =\int_{S}\frac{\left\Vert f\left( s\right) -g\left(
s\right) \right\Vert _{Y}}{1+\left\Vert f\left( s\right) -g\left( s\right)
\right\Vert _{Y}}d\nu \left( s\right) ,\ \ f,g\in L_{0}^{wo}\left( \nu
;Y\right) .
\end{equation*}%
Suppose that $\left\{ f_{n}\right\} _{n=1}^{\infty }$ is a Cauchy sequence
in $L_{0}^{\Gamma }\left( \nu ;Y\right) $. By passing to a subsequence, it
may be assumed, without loss of generality, that $d_{0}\left(
f_{n+1},f_{n}\right) \leq 2^{-n}$ for all $n\in \mathbb{N}$. This implies
that 
\begin{equation*}
\int_{S}\sum_{n=1}^{\infty }\frac{\left\Vert f_{n+1}\left( s\right)
-f_{n}\left( s\right) \right\Vert _{Y}}{1+\left\Vert f_{n+1}\left( s\right)
-f_{n}\left( s\right) \right\Vert _{Y}}d\nu \left( s\right) <\infty
\end{equation*}%
and so, there exists a $\nu $-measurable set $E\subseteq S$ such that $\nu
\left( S\diagdown E\right) =0$ and%
\begin{equation*}
\sum_{n=1}^{\infty }\frac{\left\Vert f_{n+1}\left( s\right) -f_{n}\left(
s\right) \right\Vert _{Y}}{1+\left\Vert f_{n+1}\left( s\right) -f_{n}\left(
s\right) \right\Vert _{Y}}<\infty ,\ \ s\in E.
\end{equation*}%
It is easy to see that this implies that 
\begin{equation*}
\sum_{n=1}^{\infty }\left\Vert f_{n+1}\left( s\right) -f_{n}\left( s\right)
\right\Vert _{Y}<\infty ,\ \ \ s\in E.
\end{equation*}%
Since $Y$ is a Banach space, the series $\sum_{n=1}^{\infty }\left(
f_{n+1}\left( s\right) -f_{n}\left( s\right) \right) $ is norm convergent in 
$Y$ for all $s\in E$. Defining%
\begin{equation*}
f\left( s\right) =f_{1}\left( s\right) +\sum_{n=1}^{\infty }\left(
f_{n+1}\left( s\right) -f_{n}\left( s\right) \right) ,\ \ \ s\in E
\end{equation*}%
and $f\left( s\right) =0$, $s\in S\diagdown E$, it is clear that $%
f_{n}\left( s\right) \rightarrow f\left( s\right) $ $\nu $-a.e. on $S$ and $%
f\in L_{0}^{wo}\left( \nu ;Y\right) $ (as $f$ is the $\nu $-a.e. pointwise
limit of a sequence of functions in $L_{0}^{wo}\left( \nu ;Y\right) $). It
will be shown next that $d_{0}\left( f_{n},f\right) \rightarrow 0$ as $%
n\rightarrow \infty $. Given $\varepsilon >0$, there exists $N\in \mathbb{N}$
such that 
\begin{equation*}
\int_{S}\frac{\left\Vert f_{n}\left( s\right) -f_{m}\left( s\right)
\right\Vert _{Y}}{1+\left\Vert f_{n}\left( s\right) -f_{m}\left( s\right)
\right\Vert _{Y}}d\nu \left( s\right) \leq \varepsilon
\end{equation*}%
for all $n,m\geq N$. Letting $m\rightarrow \infty $, it follows from the
dominated convergence theorem that 
\begin{equation*}
\int_{S}\frac{\left\Vert f_{n}\left( s\right) -f\left( s\right) \right\Vert
_{Y}}{1+\left\Vert f_{n}\left( s\right) -f\left( s\right) \right\Vert _{Y}}%
d\nu \left( s\right) \leq \varepsilon ,
\end{equation*}%
that is, $d_{0}\left( f_{n},f\right) \leq \varepsilon $ for all $n\geq N$.
This suffices to prove the lemma.\medskip
\end{proof}

If $a\in L_{0}\left( \nu \right) $ and $f\in L_{0}^{\Gamma }\left( \nu
;Y\right) $, then the function $af\in L_{0}^{\Gamma }\left( \nu ;Y\right) $
is defined by setting $\left( af\right) \left( s\right) =a\left( s\right)
f\left( s\right) $, $\nu $-a.e. on $S$. With respect to this action, $%
L_{0}^{\Gamma }\left( \nu ;Y\right) $ is an $L_{0}\left( \nu \right) $%
-module. Similarly, $L_{\infty }^{\Gamma }\left( \nu ;Y\right) $ has the
structure of an $L_{\infty }\left( \nu \right) $-module.

\begin{definition}
\label{NDef01}A linear operator $T:L_{0}^{\Gamma }\left( \nu ;Y\right)
\rightarrow L_{0}^{\Gamma }\left( \nu ;Y\right) $ is called $L_{0}\left( \nu
\right) $\emph{-linear} if $T\left( af\right) =aT\left( f\right) $ for all $%
f\in L_{0}^{\Gamma }\left( \nu ;Y\right) $ and all $a\in L_{0}\left( \nu
\right) $. The unital algebra of all $L_{0}\left( \nu \right) $-linear
operators in $L_{0}^{\Gamma }\left( \nu ;Y\right) $ which are continuous
with respect to $\mathcal{T}_{0}$ is denoted by $\mathcal{L}_{L_{0}}\left(
L_{0}^{\Gamma }\left( \nu ;Y\right) \right) $. Similarly, the unital algebra
of all $L_{\infty }\left( \nu \right) $-linear operators in $L_{\infty
}^{\Gamma }\left( \nu ;Y\right) $ which are continuous (with respect to the
norm in $L_{\infty }^{\Gamma }\left( \nu ;Y\right) $), is denoted by $%
\mathcal{L}_{L_{\infty }}\left( L_{\infty }^{\Gamma }\left( \nu ;Y\right)
\right) $.
\end{definition}

\begin{remark}
\label{NRem01}It may be of some interest to observe that $\mathcal{L}%
_{L_{\infty }}\left( L_{\infty }^{\Gamma }\left( \nu ;Y\right) \right) $
coincides with the algebra of all $L_{\infty }\left( \nu \right) $-linear
operators in $L_{\infty }^{\Gamma }\left( \nu ;Y\right) $ which are
continuous with respect to the topology $\mathcal{T}_{0}$. Indeed, if $T$ is
an $L_{\infty }\left( \nu \right) $-linear operator which is $\mathcal{T}%
_{0} $-continuous, then it is an immediate consequence of the closed graph
theorem that $T$ is continuous with respect to $\left\Vert \cdot \right\Vert
_{\infty }$.

For the proof of the converse implication, suppose that $T$ is an $L_{\infty
}\left( \nu \right) $-linear operator in $L_{\infty }^{\Gamma }\left( \nu
;Y\right) $ which is bounded with respect to $\left\Vert \cdot \right\Vert
_{\infty }$. We claim that%
\begin{equation}
\left\Vert \left( Tf\right) \left( s\right) \right\Vert _{Y}\leq \left\Vert
T\right\Vert _{\mathcal{L}\left( L_{\infty }^{\Gamma }\left( \nu ;Y\right)
\right) }\left\Vert f\left( s\right) \right\Vert _{Y},\ \ \ \nu \text{-a.e.
on }S,  \label{Neq02}
\end{equation}%
for all $f\in L_{\infty }^{\Gamma }\left( \nu ;Y\right) $. Indeed, suppose
that $f\in L_{\infty }^{\Gamma }\left( \nu ;Y\right) $ is such that (\ref%
{Neq02}) does not hold. Then there exists $\alpha >\left\Vert T\right\Vert _{%
\mathcal{L}\left( L_{\infty }^{\Gamma }\left( \nu ;Y\right) \right) }$ such
that the set 
\begin{equation*}
A=\left\{ s\in S:\left\Vert \left( Tf\right) \left( s\right) \right\Vert
_{Y}>\alpha \left\Vert f\left( s\right) \right\Vert _{Y}\right\}
\end{equation*}%
satisfies $\nu \left( A\right) >0$. Since $T\left( \chi _{A}f\right) =\chi
_{A}Tf$, it is clear that $\chi _{A}f\neq 0$ and 
\begin{equation*}
\left\Vert \left( T\left( \chi _{A}f\right) \right) \left( s\right)
\right\Vert _{Y}\geq \alpha \left\Vert \left( \chi _{A}f\right) \left(
s\right) \right\Vert _{Y},\ \ \ s\in S.
\end{equation*}%
Consequently, 
\begin{equation*}
\alpha \left\Vert \left( \chi _{A}f\right) \left( s\right) \right\Vert
_{Y}\leq \left\Vert T\left( \chi _{A}f\right) \right\Vert _{\infty }\leq
\left\Vert T\right\Vert _{\mathcal{L}\left( L_{\infty }^{\Gamma }\left( \nu
;Y\right) \right) }\left\Vert \chi _{A}f\right\Vert _{\infty },\ \ \ \nu 
\text{-a.e. on }S,
\end{equation*}%
and so, 
\begin{equation*}
\alpha \left\Vert \chi _{A}f\right\Vert _{\infty }\leq \left\Vert
T\right\Vert _{\mathcal{L}\left( L_{\infty }^{\Gamma }\left( \nu ;Y\right)
\right) }\left\Vert \chi _{A}f\right\Vert _{\infty }.
\end{equation*}%
This implies that $\alpha \leq \left\Vert T\right\Vert _{\mathcal{L}\left(
L_{\infty }^{\Gamma }\left( \nu ;Y\right) \right) }$, which is a
contradiction. Therefore, (\ref{Neq02}) holds, from which it is clear that $%
T $ is continuous with respect to $\mathcal{T}_{0}$.

It should be observed that it follows from (\ref{Neq02}) that any $T\in 
\mathcal{L}_{L_{\infty }}\left( L_{\infty }^{\Gamma }\left( \nu ;Y\right)
\right) $ has the property that $Tf_{n}\rightarrow Tf$ $\nu $-a.e. on $S$,
whenever $f\in L_{\infty }^{\Gamma }\left( \nu ;Y\right) $ and $\left\{
f_{n}\right\} _{n=1}^{\infty }$ is a sequence in $L_{\infty }^{\Gamma
}\left( \nu ;Y\right) $ satisfying $f_{n}\rightarrow f$ $\nu $-a.e. on $S$.
\end{remark}

For the proof of the next proposition, the following simple and well known
observation is needed (cf. \cite{Za2}, Lemma 96.6). Recall that a set $%
D\subseteq L_{0}\left( \nu \right) ^{+}$ is called order bounded if there
exists $a\in L_{0}\left( \nu \right) ^{+}$ such that $f\leq a$ for all $f\in
D$ (here $L_{0}\left( \nu \right) ^{+}$ denotes the positive cone of $%
L_{0}\left( \nu \right) $).

\begin{lemma}
\label{NLem02}Suppose that $D\subseteq L_{0}\left( \nu \right) ^{+}$ is
upwards directed. If $D$ is not order bounded in $L_{0}\left( \nu \right)
^{+}$, then there exists a sequence $\left\{ f_{n}\right\} _{n=1}^{\infty }$
in $D$ and a $\nu $-measurable set $E\subseteq S$ such that $\nu \left(
E\right) >0$ and $f_{n}\geq n\chi _{E}$ for all $n\in \mathbb{N}$.
\end{lemma}

\begin{proof}
Since $L_{0}\left( \nu \right) $ is order separable (see e.g. \cite{LZ},
Example 23.3 (iv)), there exists an increasing sequence $\left\{
g_{n}\right\} _{n=1}^{\infty }$ in $D$ which has the same upper bounds as
the set $D$. Hence, $\left\{ g_{n}\right\} _{n=1}^{\infty }$ is not order
bounded in $L_{0}\left( \nu \right) $. Define the function $h:S\rightarrow %
\left[ 0,\infty \right] $ by setting $h\left( s\right) =\sup_{n\in \mathbb{N}%
}g_{n}\left( s\right) $, $s\in S$. Since $\left\{ g_{n}\right\}
_{n=1}^{\infty }$ is not order bounded in $L_{0}\left( \nu \right) $, there
exists a $\nu $-measurable set $H$ such that%
\begin{equation*}
H\subseteq \left\{ s\in S:h\left( s\right) =\infty \right\}
\end{equation*}%
and $0<\nu \left( H\right) <\infty $. For $k,n\in \mathbb{N}$ define 
\begin{equation*}
G_{n,k}=\left\{ s\in H:g_{n}\left( s\right) \leq k\right\} .
\end{equation*}%
Since $G_{n,k}\downarrow _{n}\emptyset $ for each $k\in \mathbb{N}$, there
exists $n_{k}\in \mathbb{N}$ such that $\nu \left( G_{n_{k},k}\right) \leq
2^{-\left( k+1\right) }\nu \left( H\right) $. Defining 
\begin{equation*}
G=\bigcup_{k=1}^{\infty }G_{n_{k},k},
\end{equation*}%
it follows that $\nu \left( G\right) \leq \frac{1}{2}\nu \left( H\right) $
and so, $\nu \left( H\diagdown G\right) \geq \frac{1}{2}\nu \left( H\right)
>0$. Furthermore, 
\begin{equation*}
g_{n_{k}}\geq k\chi _{H\diagdown G_{n_{k},k}}\geq k\chi _{H\diagdown G},\ \
\ k\in \mathbb{N}.
\end{equation*}%
Hence, the set $E=H\diagdown G$ and the functions $f_{k}=g_{n_{k}}$, $%
k=1,2,\ldots $, have the desired properties.\medskip
\end{proof}

\begin{proposition}
\label{NProp01}An $L_{0}\left( \nu \right) $-linear operator $%
T:L_{0}^{\Gamma }\left( \nu ;Y\right) \rightarrow L_{0}^{\Gamma }\left( \nu
;Y\right) $ is continuous with respect to $\mathcal{T}_{0}$ if and only if
there exists a function $0\leq a\in L_{0}\left( \nu \right) $ such that 
\begin{equation}
\left\Vert \left( Tf\right) \left( s\right) \right\Vert _{Y}\leq a\left(
s\right) \left\Vert f\left( s\right) \right\Vert _{Y},\ \ \ \nu \text{-a.e.
on }S,  \label{Neq03}
\end{equation}%
for all $f\in L_{0}^{\Gamma }\left( \nu ;Y\right) $.
\end{proposition}

\begin{proof}
If there exists a function $0\leq a\in L_{0}\left( \nu \right) $ satisfying (%
\ref{Neq03}), then it is clear that $T$ is continuous with respect to the
topology $\mathcal{T}_{0}$ (see the discussion preceding Lemma \ref{NLem01}).

Assuming that $T$ is $\mathcal{T}_{0}$-continuous, define the subset $%
D\subseteq L_{0}\left( \nu \right) ^{+}$ by setting 
\begin{equation*}
D=\left\{ \left\Vert \left( Tf\right) \left( \cdot \right) \right\Vert
_{Y}:f\in L_{\infty }^{\Gamma }\left( \nu ;Y\right) ,\left\Vert f\right\Vert
_{\infty }\leq 1\right\} .
\end{equation*}%
It should be observed that $D$ is upwards directed. Indeed, suppose that $%
f,g\in L_{\infty }^{\Gamma }\left( \nu ;Y\right) $ with $\left\Vert
f\right\Vert _{\infty }$, $\left\Vert g\right\Vert _{\infty }\leq 1$ and let 
\begin{equation*}
A=\left\{ s\in S:\left\Vert \left( Tf\right) \left( s\right) \right\Vert
_{Y}\geq \left\Vert \left( Tg\right) \left( s\right) \right\Vert
_{Y}\right\} .
\end{equation*}%
Defining $h\in L_{\infty }^{\Gamma }\left( \nu ;Y\right) $ by setting $%
h=\chi _{A}f+\chi _{S\diagdown A}g$, it is clear that $\left\Vert
h\right\Vert _{\infty }\leq 1$. Since $T$ is $L_{0}$-linear, it follows that 
$Th=\chi _{A}Tf+\chi _{S\diagdown A}Tg$ and so, 
\begin{equation*}
\left\Vert \left( Th\right) \left( s\right) \right\Vert _{Y}=\chi _{A}\left(
s\right) \left\Vert \left( Tf\right) \left( s\right) \right\Vert _{Y}+\chi
_{S\diagdown A}\left( s\right) \left\Vert \left( Tg\right) \left( s\right)
\right\Vert _{Y},\ \ \ s\in S.
\end{equation*}%
This implies that $\left\Vert Th\right\Vert _{Y}\geq \left\Vert
Tf\right\Vert _{Y}\vee \left\Vert Tg\right\Vert _{Y}$ in $L_{0}\left( \nu
\right) $. Hence, $D$ is upwards directed.

Suppose that $D$ is not order bounded in $L_{0}\left( \nu \right) ^{+}$. It
follows from Lemma \ref{NLem02} that there exists a sequence $\left\{
f_{n}\right\} _{n=1}^{\infty }$ in $L_{\infty }^{\Gamma }\left( \nu
;Y\right) $, with $\left\Vert f_{n}\right\Vert _{\infty }\leq 1$ for all $n$%
, and a $\nu $-measurable set $E\subseteq S$, with $0<\nu \left( E\right)
<\infty $, such that 
\begin{equation}
\left\Vert Tf_{n}\right\Vert _{Y}\geq n\chi _{E},\ \ \ n\in \mathbb{N}.
\label{Neq07}
\end{equation}%
Defining $g_{n}=\left( 1/n\right) f_{n}$, $n\in \mathbb{N}$, it is clear
that $\left\Vert g_{n}\right\Vert _{\infty }\rightarrow 0$ and hence, $%
g_{n}\rightarrow 0$ with respect to $\mathcal{T}_{0}$. Since $T$ is assumed
to be $\mathcal{T}_{0}$-continuous, this implies that $Tg_{n}\rightarrow 0$
with respect to $\mathcal{T}_{0}$. On the other hand, it follows from (\ref%
{Neq07}) that $\left\Vert Tg_{n}\right\Vert _{Y}\geq \chi _{E}$ for all $%
n\in \mathbb{N}$, which clearly is a contradiction.

Consequently, there exists a function $0\leq a\in L_{0}\left( \nu \right) $
such that 
\begin{equation*}
\left\Vert \left( Tf\right) \left( s\right) \right\Vert _{Y}\leq a\left(
s\right) ,\ \ \ \ \nu \text{-a.e. on }S,
\end{equation*}%
for all $f\in L_{\infty }^{\Gamma }\left( \nu ;Y\right) $ satisfying $%
\left\Vert f\right\Vert _{\infty }\leq 1$. If $f\in L_{0}^{\Gamma }\left(
\nu ;Y\right) $ is arbitrary, then $f$ may be written as $f=b\tilde{f}$,
where $\tilde{f}\in L_{\infty }^{\Gamma }\left( \nu ;Y\right) $ satisfies $%
\left\Vert \tilde{f}\left( s\right) \right\Vert _{Y}\leq 1$, $s\in S$, and $%
b\in L_{0}\left( \nu \right) $ satisfies $0\leq b\left( s\right) \leq
\left\Vert f\left( s\right) \right\Vert _{Y}$, $s\in S$ (indeed, define $%
\tilde{f}\left( s\right) =f\left( s\right) /\left\Vert f\left( s\right)
\right\Vert _{Y}$ and $b\left( s\right) =\left\Vert f\left( s\right)
\right\Vert _{Y}$ whenever $f\left( s\right) \neq 0$ and let $\tilde{f}%
\left( s\right) =0$ and $b\left( s\right) =0$ whenever $f\left( s\right) =0$%
). Since $T$ is $L_{0}$-linear, it follows that 
\begin{equation*}
\left\Vert \left( Tf\right) \left( s\right) \right\Vert _{Y}=b\left(
s\right) \left\Vert \left( T\tilde{f}\right) \left( s\right) \right\Vert
_{Y}\leq a\left( s\right) \left\Vert f\left( s\right) \right\Vert _{Y},\ \ \
s\in S.
\end{equation*}%
The proof is complete.\medskip
\end{proof}

In the sequel, we will also use the following simple observation.

\begin{proposition}
\label{PProp01}Suppose that $\mathfrak{A}$ is a linear subspace of $%
L_{0}^{\Gamma }\left( \nu ;Y\right) $, which is also a module over $%
L_{\infty }\left( \nu \right) $ (that is, $af\in \mathfrak{A}$ for all $a\in
L_{\infty }\left( \nu \right) $ and $f\in \mathfrak{A}$). Furthermore,
assume that for every $f\in L_{0}^{\Gamma }\left( \nu ;Y\right) $ there
exists a $\nu $-measurable partition $\left\{ A_{n}\right\} _{n=1}^{\infty }$
of $S$ such that $\chi _{A_{n}}f\in \mathfrak{A}$ for all $n\in \mathbb{N}$.
If $T_{0}:\mathfrak{A}\rightarrow L_{0}^{\Gamma }\left( \nu ;Y\right) $ is
an $L_{\infty }\left( \nu \right) $-linear operator, then $T_{0}$ has a
unique extension to an $L_{0}$-linear operator $T:L_{0}^{\Gamma }\left( \nu
;Y\right) \rightarrow L_{0}^{\Gamma }\left( \nu ;Y\right) $. The operator $T$
is uniquely determined by the property that $\chi _{A}Tf=T_{0}\left( \chi
_{A}f\right) $ for all $f\in L_{0}^{\Gamma }\left( \nu ;Y\right) $ and all $%
\nu $-measurable sets $A\subseteq S$ satisfying $\chi _{A}f\in \mathfrak{A}$.
\end{proposition}

\begin{proof}
It will be shown first that for every $f\in L_{0}^{\Gamma }\left( \nu
;Y\right) $ there exists a unique $Tf\in L_{0}^{\Gamma }\left( \nu ;Y\right) 
$ such that $\chi _{A}Tf=T_{0}\left( \chi _{A}f\right) $ for all $\nu $%
-measurable sets $A\subseteq S$ satisfying $\chi _{A}f\in \mathfrak{A}$. To
this end, let $f\in L_{0}^{\Gamma }\left( \nu ;Y\right) $ be given. By the
assumption on $\mathfrak{A}$, there exists a $\nu $-measurable partition $%
\left\{ A_{n}\right\} _{n=1}^{\infty }$ of $S$ such that $\chi _{A_{n}}f\in 
\mathfrak{A}$ for all $n$. Since $T_{0}$ is $L_{\infty }$-linear, it follows
that 
\begin{equation*}
T_{0}\left( \chi _{A_{n}}f\right) =\chi _{A_{n}}T_{0}\left( \chi
_{A_{n}}f\right)
\end{equation*}%
and so, it is clear that $T_{0}\left( \chi _{A_{n}}f\right) $ is supported
on $A_{n}$. Therefore, 
\begin{equation*}
Tf=\sum_{n=1}^{\infty }T_{0}\left( \chi _{A_{n}}f\right)
\end{equation*}%
is well defined as a pointwise sum and $Tf\in L_{0}^{\Gamma }\left( \nu
;Y\right) $. If $A\subseteq S$ is any $\nu $-measurable set such that $\chi
_{A}f\in \mathfrak{A}$, then 
\begin{equation*}
\chi _{A}Tf=\sum_{n=1}^{\infty }\chi _{A}T_{0}\left( \chi _{A_{n}}f\right)
=\sum_{n=1}^{\infty }\chi _{A_{n}}T_{0}\left( \chi _{A}f\right) =T_{0}\left(
\chi _{A}f\right) .
\end{equation*}%
To prove the uniqueness of $Tf$, suppose that $g\in L_{0}^{\Gamma }\left(
\nu ;Y\right) $ is such that $\chi _{A}g=T_{0}\left( \chi _{A}f\right) $ for
all $\nu $-measurable sets $A\subseteq S$ satisfying $\chi _{A}f\in 
\mathfrak{A}$. This implies, in particular, that $\chi _{A_{n}}g=T_{0}\left(
\chi _{A_{n}}f\right) $ for all $n$ and so, 
\begin{equation*}
g=\sum_{n=1}^{\infty }\chi _{A_{n}}g=\sum_{n=1}^{\infty }T_{0}\left( \chi
_{A_{n}}f\right) =Tf.
\end{equation*}%
It follows, in particular, that the definition of $Tf$ does not depend on
the particular choice of the partition $\left\{ A_{n}\right\} _{n=1}^{\infty
}$ with the property that $\chi _{A_{n}}f\in \mathfrak{A}$ for all $n$.

To show that $T$ is linear, let $f,g\in L_{0}^{\Gamma }\left( \nu ;Y\right) $
be given. A moment's reflection shows that there exists a $\nu $-measurable
partition $\left\{ A_{n}\right\} _{n=1}^{\infty }$ of $S$ with the property
that $\chi _{A_{n}}f,\chi _{A_{n}}g\in \mathfrak{A}$ for all $n$. If $\alpha
,\beta \in \mathbb{C}$, then 
\begin{eqnarray*}
\chi _{A_{n}}T\left( \alpha f+\beta g\right) &=&T_{0}\left( \alpha \chi
_{A_{n}}f+\beta \chi _{A_{n}}g\right) =\alpha T_{0}\left( \chi
_{A_{n}}f\right) +\beta T_{0}\left( \chi _{A_{n}}g\right) \\
&=&\alpha \chi _{A_{n}}Tf+\beta \chi _{A_{n}}Tg=\chi _{A_{n}}\left( \alpha
Tf+\beta Tg\right)
\end{eqnarray*}%
for all $n$ and hence, $T\left( \alpha f+\beta g\right) =\alpha Tf+\beta Tg$.

Finally, it will be shown that $T$ is $L_{0}$-linear. Suppose that $g\in
L_{0}\left( \nu \right) $ and that $f\in L_{0}^{\Gamma }\left( \nu ;Y\right) 
$ and let $\left\{ A_{n}\right\} _{n=1}^{\infty }$ be a $\nu $-measurable
partition of $S$ such that $\chi _{A_{n}}g\in L_{\infty }\left( \nu \right) $
and $\chi _{A_{n}}f\in \mathfrak{A}$ for all $n$. Since $\chi _{A_{n}}gf\in 
\mathfrak{A}$, it follows that 
\begin{equation*}
\chi _{A_{n}}T\left( gf\right) =T_{0}\left( \chi _{A_{n}}gf\right) =\chi
_{A_{n}}gT_{0}\left( \chi _{A_{n}}f\right) =\chi _{A_{n}}gTf
\end{equation*}%
for all $n$. Hence, $T\left( gf\right) =gTf$, which completes the proof of
the proposition.\medskip
\end{proof}

It should be observed that, in particular, the subspace $\mathfrak{A}%
=L_{\infty }^{\Gamma }\left( \nu ;Y\right) $ satisfies the hypothesis of the
above proposition. Indeed, if $f\in L_{0}^{\Gamma }\left( \nu ;Y\right) $,
then the sets 
\begin{equation*}
A_{n}=\left\{ s\in S:n-1\leq \left\Vert f\left( s\right) \right\Vert
_{Y}<n\right\} ,\ \ \ n\in \mathbb{N}
\end{equation*}%
are $\nu $-measurable (see Lemma \ref{NLem04}) and satisfy $f\chi
_{A_{n}}\in L_{\infty }^{\Gamma }\left( \nu ;Y\right) $ for all $n$.\medskip

It is assumed now that $X$ is a \textit{separable Banach} space and the
notation and terminology introduced in Remark \ref{NRem02} will be employed.

\begin{lemma}
\label{Lem01}If $F:S\rightarrow X$ is Bochner $\nu $-measurable and $%
f:S\rightarrow \mathcal{L}\left( X\right) $ is weak operator $\nu $%
-measurable, then the function $s\longmapsto f\left( s\right) F\left(
s\right) $, $s\in S$, is Bochner $\nu $-measurable.
\end{lemma}

\begin{proof}
First suppose that $F\in \func{sim}\left( \Sigma ;X\right) $, that is, there
exist $A_{1},\ldots ,A_{n}\in \Sigma $ and $x_{1},\ldots ,x_{n}\in X$ such
that $F=\sum_{j=1}^{n}\chi _{A_{j}}\otimes x_{j}$. This implies that 
\begin{equation*}
f\left( s\right) F\left( s\right) =\sum_{j=1}^{n}\chi _{A_{j}}\left(
s\right) f\left( s\right) x_{j},\ \ \ s\in S.
\end{equation*}%
Since the functions $s\longmapsto f\left( s\right) x_{j}$, $s\in S$, are
Bochner $\nu $-measurable (see Remark \ref{NRem02} (b)), it is clear that
the function $s\longmapsto f\left( s\right) F\left( s\right) $, $s\in S$, is
Bochner $\nu $-measurable. The general result now follows from the fact that
the pointwise limit of a sequence of Bochner $\nu $-measurable functions is
again Bochner $\nu $-measurable (see \cite{DU}, Section II.1). \medskip
\end{proof}

\begin{corollary}
With respect to pointwise multiplication, the space \linebreak $\mathcal{L}%
_{0}^{wo}\left( \nu ;\mathcal{L}\left( X\right) \right) $ is a unital
algebra. Consequently, $L_{0}^{wo}\left( \nu ;\mathcal{L}\left( X\right)
\right) $ is also a unital algebra.
\end{corollary}

\begin{proof}
Let $f,g\in \mathcal{L}_{0}^{wo}\left( \nu ;\mathcal{L}\left( X\right)
\right) $ be given. If $x\in X$, then the function $s\longmapsto g\left(
s\right) x$, $s\in S$, is Bochner $\nu $-measurable (see Remark \ref{NRem02}
(b)). Hence, by Lemma \ref{Lem01}, the function $s\longmapsto f\left(
s\right) g\left( s\right) x$, $s\in S$, is Bochner $\nu $-measurable.
Consequently (see Remark \ref{NRem02} (b)), the function $s\longmapsto
f\left( s\right) g\left( s\right) $ is weak operator $\nu $-measurable, that
is, $fg\in \mathcal{L}_{0}^{wo}\left( \nu ;\mathcal{L}\left( X\right)
\right) $. Since the $\nu $-null functions form a two-sided ideal in $%
\mathcal{L}_{0}^{wo}\left( \nu ;\mathcal{L}\left( X\right) \right) $, it is
now evident that $L_{0}^{wo}\left( \nu ;\mathcal{L}\left( X\right) \right) $
is an algebra.\medskip
\end{proof}

Next, the algebras $L_{\infty }^{wo}\left( \nu ;\mathcal{L}\left( X\right)
\right) $ and $L_{0}^{wo}\left( \nu ;\mathcal{L}\left( X\right) \right) $
will be identified with algebras of linear operators on the spaces $%
L_{\infty }\left( \nu ;X\right) $ and $L_{0}\left( \nu ;X\right) $,
respectively. Suppose that $f\in L_{\infty }^{wo}\left( \nu ;\mathcal{L}%
\left( X\right) \right) $. For $F\in L_{\infty }\left( \nu ;X\right) $,
define $m_{f}F:S\rightarrow X$ by setting 
\begin{equation*}
\left( m_{f}F\right) \left( s\right) =f\left( s\right) F\left( s\right) ,\ \
\ s\in S.
\end{equation*}%
It follows from Lemma \ref{Lem01} that $m_{f}F$ is Bochner measurable.
Moreover, 
\begin{equation}
\left\Vert \left( m_{f}F\right) \left( s\right) \right\Vert _{X}\leq
\left\Vert f\left( s\right) \right\Vert _{\mathcal{L}\left( X\right)
}\left\Vert F\left( s\right) \right\Vert _{X}\leq \left\Vert f\right\Vert
_{\infty }\left\Vert F\right\Vert _{\infty },\ \ \ \nu \text{-a.e. on }S.
\label{eq01}
\end{equation}%
Consequently, $m_{f}:L_{\infty }\left( \nu ;X\right) \rightarrow L_{\infty
}\left( \nu ;X\right) $ is a bounded linear map satisfying $\left\Vert
m_{f}\right\Vert _{\mathcal{L}\left( L_{\infty }\left( \nu ,X\right) \right)
}\leq \left\Vert f\right\Vert _{\infty }$. The next lemma shows that the
latter inequality is actually an equality.

\begin{lemma}
\label{Lem02}If $f\in L_{\infty }^{wo}\left( \nu ;\mathcal{L}\left( X\right)
\right) $, then $\left\Vert m_{f}\right\Vert _{\mathcal{L}\left( L_{\infty
}\left( \nu ,X\right) \right) }=\left\Vert f\right\Vert _{\infty }$.
\end{lemma}

\begin{proof}
If $x\in X$, then $m_{f}\left( \mathbf{1}\otimes x\right) \left( s\right)
=f\left( s\right) x$ $\nu $-a.e. on $S$. Hence, if $\left\Vert x\right\Vert
_{X}\leq 1$, then 
\begin{eqnarray*}
\left\Vert f\left( s\right) x\right\Vert _{X} &\leq &\left\Vert m_{f}\left( 
\mathbf{1}\otimes x\right) \right\Vert _{\infty }\leq \left\Vert
m_{f}\right\Vert _{\mathcal{L}\left( L_{\infty }\left( \nu ,X\right) \right)
}\left\Vert \left( \mathbf{1}\otimes x\right) \right\Vert _{\infty } \\
&=&\left\Vert m_{f}\right\Vert _{\mathcal{L}\left( L_{\infty }\left( \nu
,X\right) \right) }\left\Vert x\right\Vert _{X}\leq \left\Vert
m_{f}\right\Vert _{\mathcal{L}\left( L_{\infty }\left( \nu ,X\right) \right)
},\ \ \ \nu \text{-a.e. on }S,
\end{eqnarray*}%
where the exceptional set depends on $x$. Since $X$ is separable, there
exists a countable set $\left\{ x_{n}\right\} _{n=1}^{\infty }$ which is
dense in the unit ball of $X$. It follows from the above that $\left\Vert
f\left( s\right) x_{n}\right\Vert _{X}\leq \left\Vert m_{f}\right\Vert _{%
\mathcal{L}\left( L_{\infty }\left( \nu ,X\right) \right) }$ for all $n$, $%
\nu $-a.e. on $S$. Consequently, 
\begin{equation*}
\left\Vert f\left( s\right) \right\Vert _{\mathcal{L}\left( X\right)
}=\sup_{n}\left\Vert f\left( s\right) x_{n}\right\Vert _{X}\leq \left\Vert
m_{f}\right\Vert _{\mathcal{L}\left( L_{\infty }\left( \nu ,X\right) \right)
},\ \ \ \nu \text{-a.e. on }S,
\end{equation*}%
and so, $\left\Vert f\right\Vert _{\infty }\leq \left\Vert m_{f}\right\Vert
_{\mathcal{L}\left( L_{\infty }\left( \nu ,X\right) \right) }$. The proof is
complete.\medskip
\end{proof}

If $f\in L_{\infty }^{wo}\left( \nu ;\mathcal{L}\left( X\right) \right) $,
then it is clear that $m_{f}\in \mathcal{L}_{L_{\infty }}\left( L_{\infty
}\left( \nu ;X\right) \right) $ (see Definition \ref{NDef01}) and, by Lemma %
\ref{Lem02}, it is clear that the map $f\longmapsto m_{f}$ is a linear
isometry from $L_{\infty }^{wo}\left( \nu ;\mathcal{L}\left( X\right)
\right) $ into $\mathcal{L}_{L_{\infty }}\left( L_{\infty }\left( \nu
;X\right) \right) $. The next objective is to show that this isometry is
surjective. The proof of this proposition is modeled after the proof of
Lemma 6 in \cite{dPR}, where $L_{\infty }^{wo}\left( \nu ;\mathcal{L}\left(
X\right) \right) $ is considered as an algebra of operators on the space $%
L_{p}\left( \nu ;X\right) $, $1\leq p<\infty $.

\begin{proposition}
\label{Propo02}If $T\in \mathcal{L}_{L_{\infty }}\left( L_{\infty }\left(
\nu ;X\right) \right) $, then there exists a unique $f\in L_{\infty
}^{wo}\left( \nu ;\mathcal{L}\left( X\right) \right) $ such that $T=m_{f}$.
\end{proposition}

\begin{proof}
To avoid any possible confusion, in the present proof, we use the pedantic
notation $\dot{F}$ for the equivalence class (with respect to the measure $%
\nu $) of a $\nu $-measurable function $F$ on $S$.

Let $T\in \mathcal{L}_{L_{\infty }}\left( L_{\infty }\left( \nu ;X\right)
\right) $ be given. Since 
\begin{equation*}
\left\Vert T\left( \mathbf{1}\otimes x\right) \right\Vert _{\infty }\leq
\left\Vert T\right\Vert _{\mathcal{L}\left( L_{\infty }\left( \nu ;X\right)
\right) }\left\Vert x\right\Vert _{X},\ \ \ x\in X,
\end{equation*}%
we may choose $F_{x}\in \mathcal{L}_{\infty }\left( \nu ;X\right) $ for each 
$x\in X$, such that 
\begin{equation*}
\dot{F}_{x}=T\left( \mathbf{1}\otimes x\right) \ \ \text{and}\ \ \left\Vert
F_{x}\left( s\right) \right\Vert _{X}\leq \left\Vert T\right\Vert _{\mathcal{%
L}\left( L_{\infty }\left( \nu ;X\right) \right) }\left\Vert x\right\Vert
_{X},\ \ \ s\in S.
\end{equation*}%
Since the Banach space $X$ is assumed to be separable, there exists a
countable linearly independent subset $\left\{ x_{n}\right\} _{n=1}^{\infty
} $ such that $\overline{\func{span}}\left\{ x_{n}\right\} _{n=1}^{\infty
}=X $. Let $W$ be the $\mathbb{Q}+i\mathbb{Q}$-linear span of the set $%
\left\{ x_{n}\right\} _{n=1}^{\infty }$ and note that $W$ is countable and
dense in $X$. Since the vectors $\left\{ x_{n}\right\} _{n=1}^{\infty }$ are
linearly independent, for each $s\in S$ there exists a unique $\mathbb{Q}+i%
\mathbb{Q}$-linear map $f_{0}\left( s\right) :W\rightarrow X$ such that $%
f_{0}\left( s\right) x_{n}=F_{x_{n}}\left( s\right) $ for all $n\in \mathbb{N%
}$. If $w\in W$, then $w=\sum_{j=1}^{n}\alpha _{j}x_{j}$ with $\alpha
_{j}\in \mathbb{Q}+i\mathbb{Q}$, $1\leq j\leq n$. Hence, 
\begin{equation*}
\left( f_{0}\left( \cdot \right) w\right) ^{\cdot }=\sum_{j=1}^{n}\alpha _{j}%
\dot{F}_{x_{n}}=\sum_{j=1}^{n}\alpha _{j}T\left( \mathbf{1}\otimes
x_{n}\right) =T\left( \mathbf{1}\otimes w\right) =\dot{F}_{w}.
\end{equation*}%
Since $W$ is countable, this implies that there exists a $\nu $-measurable
set $A\subseteq S$ such that $\nu \left( S\diagdown A\right) =0$ and 
\begin{equation}
f_{0}\left( s\right) w=F_{w}\left( s\right) ,\ \ \ s\in A,w\in W.
\label{eq03}
\end{equation}%
In particular, 
\begin{equation*}
\left\Vert f_{0}\left( s\right) w\right\Vert _{X}\leq \left\Vert
T\right\Vert _{\mathcal{L}\left( L_{\infty }\left( \nu ;X\right) \right)
}\left\Vert w\right\Vert _{X},\ \ \ s\in A,w\in W,
\end{equation*}%
which implies that for each $s\in A$, the $\mathbb{Q}+i\mathbb{Q}$-linear
map $f_{0}\left( s\right) :W\rightarrow X$ extends uniquely to a bounded ($%
\mathbb{C}$-) linear map $f\left( s\right) :X\rightarrow X$ satisfying $%
\left\Vert f\left( s\right) \right\Vert _{\mathcal{L}\left( X\right) }\leq
\left\Vert T\right\Vert _{\mathcal{L}\left( L_{\infty }\left( \nu ;X\right)
\right) }$. Setting $f\left( s\right) =0$ whenever $s\in S\diagdown A$, this
defines a function $f:S\rightarrow \mathcal{L}\left( X\right) $. A moment's
reflection shows that the function $f$ is weak operator $\nu $-measurable
and so, $f\in \mathcal{L}_{\infty }^{wo}\left( \nu ;\mathcal{L}\left(
X\right) \right) $.

We claim that 
\begin{equation}
\left( f\left( \cdot \right) x\right) ^{\cdot }=T\left( \mathbf{1}\otimes
x\right) ,\ \ \ x\in X.  \label{eq04}
\end{equation}%
Indeed, given $x\in X$, there exists a sequence $\left\{ w_{k}\right\}
_{k=1}^{\infty }$ in $W$ such that $\left\Vert x-w_{k}\right\Vert
_{X}\rightarrow 0$ as $k\rightarrow \infty $. Since $f\left( s\right) \in 
\mathcal{L}\left( X\right) $, this implies that $f\left( s\right)
w_{k}\rightarrow f\left( s\right) x$ for all $s\in S$. On the other hand, $%
w_{k}\rightarrow x$ in $X$ implies that $\mathbf{1}\otimes w_{k}\rightarrow 
\mathbf{1}\otimes x$ in $L_{\infty }\left( \nu ;X\right) $ and so, $T\left( 
\mathbf{1}\otimes w_{k}\right) \rightarrow T\left( \mathbf{1}\otimes
x\right) $, that is, $\dot{F}_{w_{k}}\rightarrow \dot{F}_{x}$ in $L_{\infty
}\left( \nu ;X\right) $. Consequently, $F_{w_{k}}\left( s\right) \rightarrow
F_{x}\left( s\right) $ $\nu $-a.e. on $S$. It follows from (\ref{eq03}) that 
$f\left( s\right) w_{k}=f_{0}\left( s\right) w_{k}=F_{w_{k}}\left( s\right) $
for all $s\in A$ and all $k$. Therefore, $F_{x}\left( s\right) =f\left(
s\right) x$ $\nu $-a.e. on $S$, that is, $\left( f\left( \cdot \right)
x\right) ^{\cdot }=\dot{F}_{x}=T\left( \mathbf{1}\otimes x\right) $, by
which the claim (\ref{eq04}) is proved.

Let $m_{f}\in \mathcal{L}_{L_{\infty }}\left( L_{\infty }\left( \nu
;X\right) \right) $ be the operator corresponding to \linebreak $\dot{f}\in
L_{\infty }^{wo}\left( \nu ;\mathcal{L}\left( X\right) \right) $. It will be
shown next that $T=m_{f}$. It follows from (\ref{eq04}) that $m_{f}\left( 
\mathbf{1}\otimes x\right) =T\left( \mathbf{1}\otimes x\right) $ for all $%
x\in X$. By the $L_{\infty }\left( \nu \right) $-linearity of both
operators, this implies that 
\begin{equation}
m_{f}\left( g\otimes x\right) =T\left( g\otimes x\right) ,\ \ \ g\in
L_{\infty }\left( \nu \right) ,\ \ \ x\in X.  \label{Peq02}
\end{equation}%
Consequently, $m_{f}G=TG$ for all $G\in L_{\infty }\left( \nu \right)
\otimes X$. Given $F\in L_{\infty }\left( \nu ;X\right) $, it follows from
the definition of Bochner measurability that there exists a sequence $%
\left\{ F_{n}\right\} _{n=1}^{\infty }$ in $L_{\infty }\left( \nu \right)
\otimes X$ such that $F_{n}\left( s\right) \rightarrow F\left( s\right) $ $%
\nu $-a.e. on $S$. as $n\rightarrow \infty $. This clearly implies that $%
m_{f}F_{n}\rightarrow m_{f}F$ $\nu $-a.e. Furthermore, as observed at the
end of Remark \ref{NRem01}, also $TF_{n}\rightarrow TF$ $\nu $-a.e. It
follows from (\ref{Peq02}) that, $m_{f}F_{n}=TF_{n}$ for all $n$ and hence, $%
m_{f}F=TF$. This shows that $T=m_{f}$ and the proof is complete.\bigskip
\end{proof}

In combination with Lemma \ref{Lem02}, the above proposition yields
immediately the following result.

\begin{corollary}
The map $f\longmapsto m_{f}$, $f\in L_{\infty }^{wo}\left( \nu ;\mathcal{L}%
\left( X\right) \right) $, is an isometric unital algebra isomorphism from $%
L_{\infty }^{wo}\left( \nu ;\mathcal{L}\left( X\right) \right) $ onto $%
\mathcal{L}_{L_{\infty }}\left( L_{\infty }\left( \nu ;X\right) \right) $.
\end{corollary}

Next, operators on the space $L_{0}\left( \nu ;X\right) $ will be
considered. If \linebreak $f\in L_{0}^{wo}\left( \nu ;\mathcal{L}\left(
X\right) \right) $, then it follows from Lemma \ref{Lem01} that, for every $%
F\in L_{0}\left( \nu ;X\right) $, the function $s\longmapsto f\left(
s\right) F\left( s\right) $, $s\in S$, defines an element of $L_{0}\left(
\nu ;X\right) $. Defining the map $m_{f}:L_{0}\left( \nu ;X\right)
\rightarrow L_{0}\left( \nu ;X\right) $ by setting 
\begin{equation*}
\left( m_{f}F\right) \left( s\right) =f\left( s\right) F\left( s\right) ,\ \
\ s\in S,
\end{equation*}%
it is clear that $m_{f}$ is an $L_{0}$-linear operator. Moreover, 
\begin{equation*}
\left\Vert \left( m_{f}F\right) \left( s\right) \right\Vert _{X}\leq
\left\Vert f\left( s\right) \right\Vert _{\mathcal{L}\left( X\right)
}\left\Vert F\left( s\right) \right\Vert _{X},\ \ \ s\in S,
\end{equation*}%
and so, it follows from Proposition \ref{NProp01} that $m_{f}\in \mathcal{L}%
_{L_{0}}\left( L_{0}\left( \nu ;X\right) \right) $. If $f\in
L_{0}^{wo}\left( \nu ;\mathcal{L}\left( X\right) \right) $ and $m_{f}=0$,
then $f=0$. Indeed, $m_{f}=0$ implies that $m_{f}\left( \mathbf{1}\otimes
x\right) =0$ for all $x\in X$, that is, $f\left( s\right) x=0$ $\nu $-a.e.
on $S$ for all $x\in X$. Since $X$ is separable, it follows that $f$ is a $%
\nu $-null function (cf. Lemma \ref{NLem03}). Consequently, the map $%
f\longmapsto m_{f}$ is a unital algebra isomorphism from $L_{0}^{wo}\left(
\nu ;\mathcal{L}\left( X\right) \right) $ into $\mathcal{L}_{L_{0}}\left(
L_{0}\left( \nu ;X\right) \right) $. It will be shown next that this map is
actually surjective.

\begin{proposition}
\label{PProp05}If $T\in \mathcal{L}_{L_{0}}\left( L_{0}\left( \nu ;X\right)
\right) $, then there exists a unique $f\in L_{0}^{wo}\left( \nu ;\mathcal{L}%
\left( X\right) \right) $ such that $T=m_{f}$.
\end{proposition}

\begin{proof}
Let $T\in \mathcal{L}_{L_{0}}\left( L_{0}\left( \nu ;X\right) \right) $ be
given. It follows from Proposition \ref{NProp01} that there exists $0\leq
a\in L_{0}\left( \nu \right) $ such that $\left\Vert \left( TF\right) \left(
s\right) \right\Vert _{X}\leq a\left( s\right) \left\Vert F\left( s\right)
\right\Vert _{X}$ $\nu $-a.e. on $S$ for all $F\in L_{0}\left( \nu ;X\right) 
$. Defining 
\begin{equation*}
A_{n}=\left\{ s\in S:n-1\leq a\left( s\right) <n\right\} ,\ \ \ n\in \mathbb{%
N},
\end{equation*}%
it is clear that $\left\{ A_{n}\right\} _{n=1}^{\infty }$ is a partition of $%
S$ and $a\left( s\right) \leq n$ $\nu $-a.e. on $S$. Since $T\left( \chi
_{A_{n}}F\right) =\chi _{A_{n}}T\left( F\right) $ for all $F\in L_{0}\left(
\nu ;X\right) $, it follows that $T$ leaves the subspace $L_{0}\left(
A_{n},\nu ;X\right) $ invariant. Moreover, 
\begin{equation*}
\left\Vert \left( TF\right) \left( s\right) \right\Vert _{X}\leq n\left\Vert
F\left( s\right) \right\Vert _{X},\ \ \ s\in A_{n},F\in L_{0}\left(
A_{n},\nu ;X\right) .
\end{equation*}%
This implies that $T$ leaves $L_{\infty }\left( A_{n},\nu ;X\right) $
invariant and its restriction to $L_{\infty }\left( A_{n},\nu ;X\right) $ is
bounded (and evidently, $L_{\infty }\left( A_{n},\nu \right) $-linear).
Hence, it follows from Proposition \ref{Propo02} that there exists $f_{n}\in
L_{\infty }^{wo}\left( A_{n},\nu ;\mathcal{L}\left( X\right) \right) $ such
that $\left\Vert f_{n}\right\Vert _{\infty }\leq n$ and 
\begin{equation}
\left( TF\right) \left( s\right) =f_{n}\left( s\right) F\left( s\right) ,\ \
\ s\in A_{n},F\in L_{\infty }\left( A_{n},\nu ;X\right) .  \label{eq05}
\end{equation}%
It should be observed that this implies that 
\begin{equation}
\left( TF\right) \left( s\right) =f_{n}\left( s\right) F\left( s\right) ,\ \
\ s\in A_{n},F\in L_{0}\left( A_{n},\nu ;X\right) .  \label{eq06}
\end{equation}%
Indeed, let $F\in L_{0}\left( A_{n},\nu ;X\right) $ be given and define 
\begin{equation*}
B_{k}=\left\{ s\in A_{n}:\left\Vert F\left( s\right) \right\Vert _{X}\leq
k\right\} ,\ \ \ k\in \mathbb{N}.
\end{equation*}%
Since $\chi _{B_{k}}F\in L_{\infty }\left( A_{n},\nu ;X\right) $, it follows
from (\ref{eq05}) that 
\begin{equation*}
\chi _{B_{k}}\left( s\right) \left( TF\right) \left( s\right) =T\left( \chi
_{B_{k}}F\right) \left( s\right) =\chi _{B_{k}}\left( s\right) f_{n}\left(
s\right) F\left( s\right) ,\ \ \ s\in A_{n}.
\end{equation*}%
Since $B_{k}\uparrow A_{n}$ as $k\rightarrow \infty $, this implies that (%
\ref{eq06}) holds.

Define $f\in L_{0}^{wo}\left( \nu ;\mathcal{L}\left( X\right) \right) $ by
setting $f\left( s\right) =f_{n}\left( s\right) $ whenever $s\in A_{n}$. If $%
F\in L_{0}\left( \nu ;X\right) $, then $\chi _{A_{n}}F\in L_{0}\left(
A_{n},\nu ;X\right) $ for all $n$ and so, it follows from (\ref{eq06}) that 
\begin{eqnarray*}
\chi _{A_{n}}\left( s\right) \left( TF\right) \left( s\right) &=&T\left(
\chi _{A_{n}}F\right) \left( s\right) =\chi _{A_{n}}\left( s\right)
f_{n}\left( s\right) F\left( s\right) \\
&=&\chi _{A_{n}}\left( s\right) f\left( s\right) F\left( s\right) ,\ \ \
s\in S,
\end{eqnarray*}%
for all $n\in \mathbb{N}$. Since $\bigcup\nolimits_{n=1}^{\infty }A_{n}=S$,
this implies that $\left( TF\right) \left( s\right) =f\left( s\right)
F\left( s\right) $, $\nu $-a.e. on $S$. The proof is complete.\bigskip
\end{proof}

\begin{corollary}
\label{Cor04}The map $f\longmapsto m_{f}$, $f\in L_{0}^{wo}\left( \nu ;%
\mathcal{L}\left( X\right) \right) $ is an algebra isomorphism from $%
L_{0}^{wo}\left( \nu ;\mathcal{L}\left( X\right) \right) $ onto $\mathcal{L}%
_{L_{0}}\left( L_{0}\left( \nu ;X\right) \right) $.
\end{corollary}

\section{Derivations\label{SectDer}}

In this section, some simple and well known facts concerning derivations are
recalled. For convenience of the reader, we include their short proofs. For
an account of derivations in Banach algebras, the reader is referred to \cite%
{Da}. Additional details concerning derivations in $C^{\ast }$ and $W^{\ast
} $-algebras may be found in \cite{Sak1}, \cite{Sak2} and \cite{BR}.

Let $A$ be a complex algebra. For $x,y\in A$, the \textit{commutator} $\left[
x,y\right] $ is defined by setting 
\begin{equation*}
\left[ x,y\right] =xy-yx.
\end{equation*}%
The \textit{center} of $A$ is denoted by $\mathcal{Z}\left( A\right) $, that
is, 
\begin{equation*}
\mathcal{Z}\left( A\right) =\left\{ x\in A:\left[ x,y\right] =0\ \ \forall
y\in A\right\} .
\end{equation*}%
Evidently, $\mathcal{Z}\left( A\right) $ is an abelian subalgebra of $A$.

A linear map $\delta :A\rightarrow A$ is called a \textit{derivation} if 
\begin{equation*}
\delta \left( xy\right) =\delta \left( x\right) y+x\delta \left( y\right) ,\
\ \ x,y\in A.
\end{equation*}%
If $w\in A$, then the map $\delta _{w}:A\rightarrow A$, given by $\delta
_{w}\left( x\right) =\left[ w,x\right] $, $x\in A$, is a derivation. A
derivation of this form is called \textit{inner}. It is clear that any inner
derivation on $A$ vanishes on the center. Furthermore, if $w_{1},w_{2}\in A$%
, then $\delta _{w_{1}}=\delta _{w_{2}}$ if and only if $w_{1}-w_{2}\in 
\mathcal{Z}\left( A\right) $.

\begin{lemma}
\label{PLem03}Suppose that $A$ is unital (with unit element $\mathbf{1}$).
If $\delta :A\rightarrow A$ is a derivation, then the following two
statements are equivalent:

\begin{enumerate}
\item[(i).] $\delta \left( z\right) =0$ for all $z\in \mathcal{Z}\left(
A\right) $;

\item[(ii).] $\delta $ is $\mathcal{Z}$-linear, that is, $\delta \left(
zx\right) =z\delta \left( x\right) $ for all $x\in A$ and $z\in \mathcal{Z}%
\left( A\right) $.
\end{enumerate}
\end{lemma}

\begin{proof}
If $\delta \left( z\right) =0$ for all $z\in \mathcal{Z}\left( A\right) $,
then $\delta \left( zx\right) =z\delta \left( x\right) $ for all $x\in A$
and $z\in \mathcal{Z}\left( A\right) $. Hence, (i) implies (ii). If $\delta $
is $\mathcal{Z}$-linear and $z\in \mathcal{Z}\left( A\right) $, then $\delta
\left( z\right) x=\delta \left( zx\right) -z\delta \left( x\right) =0$ for
all $x\in A$. Taking $x=\mathbf{1}$, this implies that $\delta \left(
z\right) =0$. Hence, (ii) implies (i).\medskip
\end{proof}

\begin{lemma}
\label{PLem04}If $\delta :A\rightarrow A$ is a derivation, then $\delta
\left( z\right) \in \mathcal{Z}\left( A\right) $ for all $z\in \mathcal{Z}%
\left( A\right) $.
\end{lemma}

\begin{proof}
If $z\in \mathcal{Z}\left( A\right) $, then 
\begin{equation*}
x\delta \left( z\right) =\delta \left( xz\right) -\delta \left( x\right)
z=\delta \left( zx\right) -z\delta \left( x\right) =\delta \left( z\right) x
\end{equation*}%
for all $x\in A$. Hence, $\delta \left( z\right) \in \mathcal{Z}\left(
A\right) $. \medskip
\end{proof}

The set of all \textit{idempotents} in $A$ is denoted by $\mathcal{I}\left(
A\right) $, that is, 
\begin{equation*}
\mathcal{I}\left( A\right) =\left\{ p\in A:p^{2}=p\right\} .
\end{equation*}%
In particular, the set of all \textit{central idempotents} in $A$ is denoted
by $\mathcal{I}\left( \mathcal{Z}\left( A\right) \right) $.

\begin{lemma}
\label{PLem05}Let $\delta :A\rightarrow A$ be a derivation.

\begin{enumerate}
\item[(i).] If $p\in \mathcal{I}\left( A\right) $, then $p\delta \left(
p\right) p=0$.

\item[(ii).] If $p\in \mathcal{I}\left( \mathcal{Z}\left( A\right) \right) $%
, then $\delta \left( p\right) =0$.

\item[(iii).] If $p\in \mathcal{I}\left( \mathcal{Z}\left( A\right) \right) $%
, then $\delta \left( px\right) =p\delta \left( x\right) $ for all $x\in A$.

\item[(iv).] If $p\in \mathcal{I}\left( A\right) $ and $x\in A$ are such
that $px=xp$, then $p\delta \left( x\right) p=p\delta \left( xp\right) p$.
\end{enumerate}
\end{lemma}

\begin{proof}
(i). If $p\in \mathcal{I}\left( A\right) $, then 
\begin{equation}
\delta \left( p\right) =\delta \left( p^{2}\right) =\delta \left( p\right)
p+p\delta \left( p\right) .  \label{Peq01}
\end{equation}%
Multiplying both sides of this equation by $p$ yields that $p\delta \left(
p\right) p=2p\delta \left( p\right) p$ and so, $p\delta \left( p\right) p=0$.

(ii). If $p\in \mathcal{I}\left( \mathcal{Z}\left( A\right) \right) $, then
it follows from (i) that $p\delta \left( p\right) =\delta \left( p\right)
p=0 $ and so, (\ref{Peq01}) implies that $\delta \left( p\right) =0$.

(iii). Since $\delta \left( px\right) =\delta \left( p\right) x+p\delta
\left( x\right) $, statement (iii) follows immediately from (ii).

(iv). It follows from $\delta \left( xp\right) =\delta \left( x\right)
p+x\delta \left( p\right) $ that 
\begin{equation*}
p\delta \left( xp\right) p=p\delta \left( x\right) p+px\delta \left(
p\right) p.
\end{equation*}%
Since $px\delta \left( p\right) p=xp\delta \left( p\right) p$, it follows
from (i) that $px\delta \left( p\right) p=0$ and so, $p\delta \left(
xp\right) p=p\delta \left( x\right) p$.\medskip
\end{proof}

The following lemma may be considered as an extension of (iii) of Lemma \ref%
{PLem05}.

\begin{lemma}
\label{PLem06}Suppose that $A$ is an algebra which is also a module over a
commutative unital algebra $B$, such that $\mathbf{1}_{B}x=x$, $x\in A$, and 
$b\left( xy\right) =\left( bx\right) y=x\left( by\right) $ for all $b\in B$
and $x,y\in A$. Assume furthermore that $x\in A$ and $xy=0$ for all $y\in A$
imply that $x=0$. If $\delta :A\rightarrow A$ is a derivation, then $\delta
\left( px\right) =p\delta \left( x\right) $ for all $x\in A$ and $p\in 
\mathcal{I}\left( B\right) $.
\end{lemma}

\begin{proof}
Let $p\in \mathcal{I}\left( B\right) $ and $x\in A$ be given. If $y\in A$,
then $px\left( \mathbf{1}_{B}-p\right) y=0$ and so, 
\begin{equation*}
\left( \mathbf{1}_{B}-p\right) \delta \left( px\right) y+px\delta \left(
\left( \mathbf{1}_{B}-p\right) y\right) =0.
\end{equation*}%
Multiplying this identity on the left by $\mathbf{1}_{B}-p$ shows that $%
\left( \mathbf{1}_{B}-p\right) \delta \left( px\right) y=0$. Since this
holds for all $y\in A$, it follows that $\left( \mathbf{1}_{B}-p\right)
\delta \left( px\right) =0$, that is $\delta \left( px\right) =p\delta
\left( px\right) $. Interchanging $p$ and $\mathbf{1}_{B}-p$, this also
implies that $p\delta \left( \left( \mathbf{1}_{B}-p\right) x\right) =0$.
Furthermore, 
\begin{equation*}
\delta \left( x\right) =\delta \left( \mathbf{1}_{B}x\right) =\delta \left(
px\right) +\delta \left( \left( \mathbf{1}_{B}-p\right) x\right) .
\end{equation*}%
Multiplying this identity by $p$, it follows that 
\begin{equation*}
p\delta \left( x\right) =p\delta \left( px\right) +p\delta \left( \left( 
\mathbf{1}_{B}-p\right) x\right) =\delta \left( px\right) .
\end{equation*}%
The proof of the lemma is complete.\medskip
\end{proof}

Also, the following general observation will be used.

\begin{lemma}
\label{DLem01}Let $A$ be an algebra and $\delta :A\rightarrow A$ be a
derivation. For $p\in \mathcal{I}\left( A\right) $ define the subalgebra $%
A_{p}$ of $A$ by setting 
\begin{equation*}
A_{p}=pAp=\left\{ pxp:x\in A\right\} .
\end{equation*}%
Defining $\delta _{p}:A_{p}\rightarrow A_{p}$ by setting $\delta _{p}\left(
x\right) =p\delta \left( x\right) p$, $x\in A_{p}$, the map $\delta _{p}$ is
a derivation in $A_{p}$.
\end{lemma}

\begin{proof}
It should be noted that an element $x\in A$ belongs to $A_{p}$ if and only
if $x=px=xp$. Therefore, if $x,y\in A_{p}$, then 
\begin{eqnarray*}
\delta _{p}\left( xy\right) &=&p\delta \left( xy\right) p=p\delta \left(
x\right) yp+px\delta \left( y\right) p \\
&=&p\delta \left( x\right) py+xp\delta \left( y\right) p=\delta _{p}\left(
x\right) y+x\delta _{p}\left( y\right) .
\end{eqnarray*}%
This proof is complete.\medskip
\end{proof}

The following observation will be used in the next section.

\begin{lemma}
\label{DLem02}Let $\left( S,\Sigma ,\nu \right) $ be a finite measure space.
If $0\neq \delta :L_{\infty }\left( \nu \right) \rightarrow L_{0}\left( \nu
\right) $ is a linear map satisfying $\delta \left( ab\right) =\delta \left(
a\right) b+a\delta \left( b\right) $, $a,b\in L_{\infty }\left( \nu \right) $%
, then there exists a $\nu $-measurable set $B\subseteq S$ with $\nu \left(
B\right) >0$ such that for every $\varepsilon >0$ there exists $%
a_{\varepsilon }\in L_{\infty }\left( \nu \right) $ satisfying $\left\vert
a_{\varepsilon }\right\vert \leq \varepsilon \mathbf{1}$ and $\left\vert
\delta \left( a_{\varepsilon }\right) \right\vert \geq \chi _{B}$.
\end{lemma}

\begin{proof}
If $A\subseteq S$ is $\nu $-measurable, then $\delta \left( \chi _{A}\right)
=\delta \left( \chi _{A}^{2}\right) =2\chi _{A}\delta \left( \chi
_{A}\right) $ and so, $\delta \left( \chi _{A}\right) =0$. Since $\delta
\neq 0$, there exists $a\in L_{\infty }\left( \nu \right) $ such that $%
\delta \left( a\right) \neq 0$. Multiplying $a$ by an appropriate scalar, it
may be assumed that $\left\vert \delta \left( a\right) \right\vert \geq \chi
_{B}$ for some $\nu $-measurable set $B\subseteq S$ with $\nu \left(
B\right) >0$. Given $\varepsilon >0$, there exist disjoint $\nu $-measurable
sets $\left\{ A_{n}\right\} _{n=1}^{N}$ and scalars $\left\{ \alpha
_{n}\right\} _{n=1}^{N}$ such that the function $b=\sum_{n=1}^{N}\alpha
_{n}\chi _{A_{n}}$ satisfies $\left\vert a-b\right\vert \leq \varepsilon 
\mathbf{1}$. From the observations above, it follows that $\delta \left(
b\right) =0$. Consequently, the function $a_{\varepsilon }=a-b$ satisfies
the desired conditions.\medskip
\end{proof}

\section{Admissible subalgebras of $L_{0}^{wo}\left( \protect\nu ;\mathcal{L}%
\left( X\right) \right) \label{SectAdmis}$}

Following P.R. Chernoff (\cite{Cher}), a subalgebra $\mathfrak{U}$ of $%
\mathcal{L}\left( X\right) $ is termed \textit{standard} whenever $\mathcal{F%
}\left( X\right) \subseteq \mathfrak{U}$. If $\mathfrak{U}\subseteq \mathcal{%
L}\left( X\right) $ is a standard subalgebra, then $\mathcal{I}\left( 
\mathfrak{U}\right) $ denotes the set of all projections in $\mathfrak{U}$,
that is, 
\begin{equation*}
\mathcal{I}\left( \mathfrak{U}\right) =\left\{ p\in \mathfrak{U}%
:p^{2}=p\right\} .
\end{equation*}%
Throughout the present section, it is assumed that $X$ is a \textit{%
separable Banach space} and that $\left( S,\Sigma ,\nu \right) $ is a $%
\sigma $-finite measure space.

\begin{definition}
\label{DDef01}Let $\mathfrak{U}$ be a closed standard algebra of $\mathcal{L}%
\left( X\right) $. A subalgebra $\mathfrak{A}$ of $L_{0}^{wo}\left( \nu ;%
\mathcal{L}\left( X\right) \right) $ is called $\mathfrak{U}$-\emph{%
admissible} if:

\begin{enumerate}
\item[(i).] $\mathfrak{A}$ is an $L_{\infty }\left( \nu \right) $-submodule
of $L_{0}^{wo}\left( \nu ;\mathfrak{U}\right) $, that is, $\mathfrak{A}%
\subseteq L_{0}^{wo}\left( \nu ;\mathfrak{U}\right) $ and $gf\in \mathfrak{A}
$ whenever $g\in L_{\infty }\left( \nu \right) $ and $f\in \mathfrak{A}$;

\item[(ii).] If $f\in L_{\infty }^{wo}\left( \nu ;\mathcal{L}\left( X\right)
\right) $ and if for every $\varepsilon >0$ there exists a projection $p\in 
\mathcal{I}\left( \mathcal{F}\left( X\right) \right) $ such that $\left\Vert
f\left( \mathbf{1}\otimes p^{\bot }\right) \right\Vert _{\infty }\leq
\varepsilon $, then $\chi _{A}f\in \mathfrak{A}$ for all $\nu $-measurable
sets $A\subseteq S$ with $\nu \left( A\right) <\infty $;

\item[(iii).] for every $f\in L_{0}^{wo}\left( \nu ;\mathfrak{U}\right) $
there exists a $\nu $-measurable partition $\left\{ A_{n}\right\}
_{n=1}^{\infty }$ of $S$ such that $\chi _{A_{n}}f\in \mathfrak{A}$ for all $%
n\in \mathbb{N}$.
\end{enumerate}
\end{definition}

It should be noted that condition (ii) implies that $L_{\infty }^{wo}\left(
A,\nu ;\mathcal{L}\left( X_{p}\right) \right) \subseteq \mathfrak{A}$ for
all $\nu $-measurable sets $A\subseteq S$ with $\nu \left( A\right) <\infty $
and all $p\in \mathcal{I}\left( \mathcal{F}\left( X\right) \right) $ (where $%
\mathcal{L}\left( X_{p}\right) $ is identified with a closed subalgebra of $%
\mathcal{L}\left( X\right) $, as described in Section \ref{SectNotTerm}).
Furthermore, if $\mathfrak{U}\subseteq \mathcal{L}\left( X\right) $ is a
closed standard algebra, then it is evident that $L_{0}^{wo}\left( \nu ;%
\mathfrak{U}\right) $ and $L_{\infty }^{wo}\left( \nu ;\mathfrak{U}\right) $
are $\mathfrak{U}$-admissible subalgebras of $L_{0}^{wo}\left( \nu ;\mathcal{%
L}\left( X\right) \right) $. The main objective in the present section is to
show that, if $X$ is \textit{infinite dimensional}, then any derivation $%
\delta :\mathfrak{A}\rightarrow \mathfrak{A}$ is given by $\delta \left(
f\right) =\left[ w,f\right] $, $f\in \mathfrak{A}$, for some $w\in
L_{0}^{wo}\left( \nu ;\mathcal{L}\left( X\right) \right) $.

\begin{remark}
\label{DRem01}Let $\left( S,\Sigma ,\nu \right) $ be a $\sigma $-finite
measure space and $H$ be a separable Hilbert space. The von Neumann algebra $%
\mathcal{M}=L_{\infty }\left( \nu \right) \overline{\otimes }B\left(
H\right) $ is equipped with the tensor product trace $\tau $ (with respect
to the integral in $L_{\infty }\left( \nu \right) $ and the standard trace
in $B\left( H\right) $. As will be shown in Section \ref{SectMeas}, the
algebra $LS\left( \mathcal{M}\right) $ of all locally measurable operators
may be identified with the algebra $L_{0}^{wo}\left( \nu ;B\left( H\right)
\right) $. As we shall see in Proposition \ref{PProp07}, the above
definition of admissible algebras has been chosen such that:

\begin{enumerate}
\item[(a).] The algebra $S\left( \mathcal{M}\right) $ of all measurable
operators affiliated with $\mathcal{M}$ and the algebra $S\left( \tau
\right) $ of all $\tau $-measurable operators are both $\mathfrak{U}$%
-admissible with respect to $\mathfrak{U}=B\left( H\right) .$

\item[(b).] The algebra $S_{0}\left( \tau \right) $ of all $\tau $-compact
operators is $\mathfrak{U}$-admissible with respect to $\mathfrak{U}=K\left(
H\right) $, the algebra of all compact operators in $H$.
\end{enumerate}
\end{remark}

Some simple observations first.

\begin{lemma}
\label{PLem07}Suppose that $\mathfrak{A}\subseteq L_{0}^{wo}\left( \nu ;%
\mathcal{L}\left( X\right) \right) $ satisfies (i) and (ii) of Definition %
\ref{DDef01}. If $f\in \mathfrak{A}$ is such that $fg=0$ for all $g\in 
\mathfrak{A}$, then $f=0$.
\end{lemma}

\begin{proof}
Let $\left\{ x_{n}\right\} _{n=1}^{\infty }$ be a dense subset of $X$. For
each $n\in \mathbb{N}$, let $x_{n}^{\ast }\in X^{\ast }$ satisfy $%
\left\langle x_{n},x_{n}^{\ast }\right\rangle =1$. Define the rank one
projection $p_{n}$ in $X$ by setting $p_{n}\left( x\right) =\left\langle
x,x_{n}^{\ast }\right\rangle x_{n}$, $x\in X$. Let the $\nu $-measurable set 
$A\subseteq S$ with $\nu \left( A\right) <\infty $ be fixed. It follows from
assumption (ii) that $\chi _{A}\otimes p_{n}$ $\in \mathfrak{A}$ and so, $%
f\left( \chi _{A}\otimes p_{n}\right) =0$ for all $n$. This implies, in
particular, that 
\begin{equation*}
f\left( s\right) x_{n}=f\left( s\right) \chi _{A}\left( s\right) p_{n}x_{n}=0%
\text{\ \ \ }\nu \text{-a.e. on }A
\end{equation*}%
for all $n\in \mathbb{N}$. Since $\left\{ x_{n}\right\} _{n=1}^{\infty }$
dense in $X$, this implies that $f\left( s\right) =0$ $\nu $-a.e. on $A$.
Using that $\nu $ is a $\sigma $-finite measure, it follows that $f=0$.
\medskip
\end{proof}

A combination of Lemmas \ref{PLem06} and \ref{PLem07} yields the following
result.

\begin{corollary}
\label{DCor01}If $\mathfrak{A}\subseteq L_{0}^{wo}\left( \nu ;\mathcal{L}%
\left( X\right) \right) $ satisfies (i) and (ii) and if $\delta :\mathfrak{A}%
\rightarrow \mathfrak{A}$ is a derivation, then $\delta \left( \chi
_{A}f\right) =\chi _{A}\delta \left( f\right) $ for all $f\in \mathfrak{A}$
and all $\nu $-measurable subsets $A\subseteq S$.
\end{corollary}

Suppose now that $\mathfrak{A}\subseteq L_{0}^{wo}\left( \nu ;\mathcal{L}%
\left( X\right) \right) $ satisfies (i) and (ii) of Definition \ref{DDef01}
and let $\delta :\mathfrak{A}\rightarrow \mathfrak{A}$ be a derivation. It
will be shown that $\delta $ is $L_{\infty }$-linear, that is, $\delta
\left( af\right) =a\delta \left( f\right) $ for all $a\in L_{\infty }\left(
\nu \right) $ and $f\in \mathfrak{A}$. \textit{It is assumed, up to Lemma %
\ref{DLem06}, in addition, that} $\nu \left( S\right) <\infty $.

If $p\in \mathcal{I}\left( \mathcal{F}\left( X\right) \right) $, then it
follows from (ii) and the discussion in Section \ref{SectNotTerm} that 
\begin{equation*}
L_{\infty }^{wo}\left( \nu ;\mathcal{L}\left( X_{p}\right) \right) \subseteq 
\mathfrak{A}_{p}=\left( \mathbf{1}\otimes p\right) \mathfrak{A}\left( 
\mathbf{1}\otimes p\right) \subseteq L_{0}^{wo}\left( \nu ;\mathcal{L}\left(
X_{p}\right) \right) .
\end{equation*}%
It is not difficult to show that the commutant of $L_{\infty }^{wo}\left(
\nu ;\mathcal{L}\left( X_{p}\right) \right) $ in \linebreak $%
L_{0}^{wo}\left( \nu ;\mathcal{L}\left( X_{p}\right) \right) $ is equal to $%
L_{0}\left( \nu \right) \otimes \mathbb{C}p$ and so, the centre of $%
\mathfrak{A}_{p}$ is given by $Z\left( \mathfrak{A}_{p}\right) =Z_{p}\otimes 
\mathbb{C}p$, where $Z_{p}\subseteq L_{0}\left( \nu \right) $ is given by 
\begin{equation*}
Z_{p}=\left\{ a\in L_{0}\left( \nu \right) :a\otimes p\in \mathfrak{A}%
\right\} .
\end{equation*}%
Note that $Z_{p}$ is a subalgebra of $L_{0}\left( \nu \right) $ satisfying $%
L_{\infty }\left( \nu \right) \subseteq Z_{p}$. It follows from Lemma \ref%
{DLem01} that the map $\delta _{p}:f\longmapsto \left( \mathbf{1}\otimes
p\right) \delta \left( f\right) \left( \mathbf{1}\otimes p\right) $, $f\in 
\mathfrak{A}_{p}$, defines a derivation on $\mathfrak{A}_{p}$. The
derivation $\delta _{p}$ leaves the centre $Z_{p}\otimes \mathbb{C}p$
invariant (see Lemma \ref{PLem04}). Consequently, there exists a unique
linear map $\delta _{p}^{z}:L_{\infty }\left( \nu \right) \rightarrow
L_{0}\left( \nu \right) $ such that 
\begin{equation}
\delta _{p}\left( a\otimes p\right) =\delta _{p}^{z}\left( a\right) \otimes
p,\ \ \ a\in L_{\infty }\left( \nu \right) ,  \label{Deq01}
\end{equation}%
and satisfying 
\begin{equation*}
\delta _{p}^{z}\left( ab\right) =\delta _{p}^{z}\left( a\right) b+a\delta
_{p}^{z}\left( a\right) ,\ \ \ a,b\in L_{\infty }\left( \nu \right) .
\end{equation*}

\begin{lemma}
\label{DLem03}If $0\neq p,q\in \mathcal{I}\left( \mathcal{F}\left( X\right)
\right) $ are such that $q=pq=qp$, then $\delta _{p}^{z}\left( a\right)
=\delta _{q}^{z}\left( a\right) $ for all $a\in L_{\infty }\left( \nu
\right) $.
\end{lemma}

\begin{proof}
If $a\in L_{\infty }\left( \nu \right) $, then 
\begin{eqnarray*}
\delta _{q}^{z}\left( a\right) \otimes q &=&\delta _{q}\left( a\otimes
q\right) =\left( \mathbf{1}\otimes q\right) \delta \left( a\otimes q\right)
\left( \mathbf{1}\otimes q\right) \\
&=&\left( \mathbf{1}\otimes q\right) \delta \left( \left( \mathbf{1}\otimes
q\right) \left( a\otimes p\right) \right) \left( \mathbf{1}\otimes q\right)
\\
&=&\left( \mathbf{1}\otimes q\right) \left[ \delta \left( \mathbf{1}\otimes
q\right) \left( a\otimes p\right) +\left( \mathbf{1}\otimes q\right) \delta
\left( a\otimes p\right) \right] \left( \mathbf{1}\otimes q\right) .
\end{eqnarray*}%
Using Lemma \ref{PLem05} (i), it follows that 
\begin{equation*}
\left( \mathbf{1}\otimes q\right) \delta \left( \mathbf{1}\otimes q\right)
\left( a\otimes p\right) \left( \mathbf{1}\otimes q\right) =\left( \mathbf{1}%
\otimes q\right) \delta \left( \mathbf{1}\otimes q\right) \left( \mathbf{1}%
\otimes q\right) \left( a\otimes p\right) =0
\end{equation*}%
and so, 
\begin{eqnarray*}
\delta _{q}^{z}\left( a\right) \otimes q &=&\left( \mathbf{1}\otimes
q\right) \delta \left( a\otimes p\right) \left( \mathbf{1}\otimes q\right) \\
&=&\left( \mathbf{1}\otimes q\right) \left( \mathbf{1}\otimes p\right)
\delta \left( a\otimes p\right) \left( \mathbf{1}\otimes p\right) \left( 
\mathbf{1}\otimes q\right) \\
&=&\left( \mathbf{1}\otimes q\right) \delta _{p}\left( a\otimes p\right)
\left( \mathbf{1}\otimes q\right) \\
&=&\left( \mathbf{1}\otimes q\right) \left( \delta _{p}^{z}\left( a\right)
\otimes p\right) \left( \mathbf{1}\otimes q\right) =\delta _{p}^{z}\left(
a\right) \otimes q.
\end{eqnarray*}%
Since $q\neq 0$, this implies that $\delta _{q}^{z}\left( a\right) =\delta
_{p}^{z}\left( a\right) $. The proof of the lemma is complete.\medskip
\end{proof}

Assume now that $X$ is infinite dimensional. Since the Banach space $X$ is
separable, there exist sequences $\left\{ x_{n}\right\} _{n=1}^{\infty }$ in 
$X$ and $\left\{ x_{n}^{\ast }\right\} _{n=1}^{\infty }$ in $X^{\ast }$ such
that (see e.g. \cite{LT1}, Proposition 1.f.3):

\begin{enumerate}
\item $\left\langle x_{n},x_{m}^{\ast }\right\rangle =\delta _{nm}$ for all $%
n,m\in \mathbb{N}$;

\item $\func{span}\left\{ x_{n}\right\} _{n=1}^{\infty }$ is dense in $X$;

\item $\left\{ x_{n}^{\ast }\right\} _{n=1}^{\infty }$ separates the points
of $X$.
\end{enumerate}

For $n\in \mathbb{N}$, define the projection $p_{n}\in \mathcal{I}\left( 
\mathcal{F}\left( X\right) \right) $ by setting 
\begin{equation*}
p_{n}x=\sum_{k=1}^{n}\left\langle x,x_{k}^{\ast }\right\rangle x_{k},\ \ \
x\in X.
\end{equation*}%
It should be observed that $p_{m}p_{n}=p_{n}p_{m}=p_{m}$ whenever $m\leq n$.
Indeed, if $x\in X$, then 
\begin{eqnarray*}
p_{m}p_{n}x &=&\sum_{k=1}^{m}\left\langle p_{n}x,x_{k}^{\ast }\right\rangle
x_{k}=\sum_{k=1}^{m}\left\langle \sum_{j=1}^{n}\left\langle x,x_{j}^{\ast
}\right\rangle x_{j},x_{k}^{\ast }\right\rangle x_{k} \\
&=&\sum_{k=1}^{m}\sum_{j=1}^{n}\left\langle x,x_{j}^{\ast }\right\rangle
\left\langle x_{j},x_{k}^{\ast }\right\rangle
x_{k}=\sum_{k=1}^{m}\sum_{j=1}^{n}\delta _{j,k}\left\langle x,x_{j}^{\ast
}\right\rangle x_{k} \\
&=&\sum_{k=1}^{m}\left\langle x,x_{k}^{\ast }\right\rangle x_{k}=p_{m}x,
\end{eqnarray*}%
and similarly, $p_{n}p_{m}x=p_{m}x$. In particular, that $\left\{
p_{n}\right\} _{n=1}^{\infty }$ is a commuting family of projections. Also
the following observation will be used.

\begin{lemma}
\label{DLem05}If $f\in L_{0}\left( \nu ;\mathcal{L}\left( X\right) \right) $
is such that $\left( \mathbf{1}\otimes p_{n}\right) f\left( \mathbf{1}%
\otimes p_{n}\right) =0$ for all $n\in \mathbb{N}$, then $f=0$.
\end{lemma}

\begin{proof}
First assume that $f\in L_{0}\left( \nu ;\mathcal{L}\left( X\right) \right) $
is such that $f\left( \mathbf{1}\otimes p_{n}\right) =0$ for all $n\in 
\mathbb{N}$, that is, 
\begin{equation*}
\sum_{k=1}^{n}\left\langle x,x_{k}^{\ast }\right\rangle f\left( s\right)
x_{k}=0,\ \ \ \nu \text{-a.e. on }S,\ x\in X,n\in \mathbb{N}.
\end{equation*}%
Taking $x=x_{n}$, this implies that $f\left( s\right) x_{n}=0$ for all $n\in 
\mathbb{N}$ and $\nu $-a.e. on $S$. Since $\func{span}\left\{ x_{n}\right\}
_{n=1}^{\infty }$ is dense in $X$, it follows that $f\left( s\right) =0$ $%
\nu $-a.e. on $S$.

Suppose next that $f\in L_{0}\left( \nu ;\mathcal{L}\left( X\right) \right) $
is such that $\left( \mathbf{1}\otimes p_{n}\right) f=0$ for all $n\in 
\mathbb{N}$, that is, 
\begin{equation*}
\sum_{k=1}^{n}\left\langle f\left( s\right) x,x_{k}^{\ast }\right\rangle
x_{k}=0,\ \ \ \nu \text{-a.e. on }S,\ x\in X,n\in \mathbb{N}.
\end{equation*}%
Since $\left\{ x_{k}\right\} _{k=1}^{\infty }$ are linearly independent,
this implies that $\left\langle f\left( s\right) x,x_{k}^{\ast
}\right\rangle =0$ $\nu $-a.e. on $S$ for all $k\in \mathbb{N}$. Using that $%
\left\{ x_{k}^{\ast }\right\} _{k=1}^{\infty }$ separates the points of $X$,
it follows that $f\left( s\right) x=0$ $\nu $-a.e. on $S$ for all $x\in X$
and so, $f=0$ (since $X$ is separable).

Assume now that $\left( \mathbf{1}\otimes p_{n}\right) f\left( \mathbf{1}%
\otimes p_{n}\right) =0$ for all $n\in \mathbb{N}$. Since 
\begin{equation*}
\left( \mathbf{1}\otimes p_{m}\right) \left( \mathbf{1}\otimes p_{n}\right) =%
\mathbf{1}\otimes p_{m}
\end{equation*}%
whenever $m\leq n$, it follows that $\left( \mathbf{1}\otimes p_{m}\right)
f\left( \mathbf{1}\otimes p_{n}\right) =0$ for all $m,n\in \mathbb{N}$.
Fixing $m\in \mathbb{N}$, it follows from the first part of the proof that $%
\left( \mathbf{1}\otimes p_{m}\right) f=0$ for all $m$. Applying the second
part of the present proof, this implies that $f=0$.\medskip
\end{proof}

Since $p_{m}p_{n}=p_{m}$ whenever $m\leq n$, it follows from Lemma \ref%
{DLem03} that $\delta _{p_{n}}^{z}=\delta _{p_{m}}^{z}$ for all $m,n\in 
\mathbb{N}$. Consequently, there exists a linear map $\delta ^{z}:L_{\infty
}\left( \nu \right) \rightarrow L_{0}\left( \nu \right) $ satisfying 
\begin{equation*}
\delta ^{z}\left( ab\right) =\delta ^{z}\left( a\right) b+a\delta ^{z}\left(
b\right) ,\ \ \ a,b\in L_{\infty }\left( \nu \right) ,
\end{equation*}%
and 
\begin{equation*}
\left( \mathbf{1}\otimes p_{n}\right) \delta \left( a\otimes p_{n}\right)
\left( \mathbf{1}\otimes p_{n}\right) =\delta ^{z}\left( a\right) \otimes
p_{n},\ \ \ a\in L_{\infty }\left( \nu \right) ,n\in \mathbb{N}.
\end{equation*}

\begin{lemma}
\label{DLem04}In the situation as above, it follows that $\delta ^{z}=0$.
\end{lemma}

\begin{proof}
Suppose that $\delta ^{z}\neq 0$. Let $q_{n}\in \mathcal{I}\left( \mathcal{F}%
\left( X\right) \right) $ be defined by $q_{n}=p_{n+1}-p_{n}$ for all $n\in 
\mathbb{N}$. Note that $q_{n}\neq 0$ for all $n$. It follows from Lemma \ref%
{DLem02} that there exists a $\nu $-measurable subset $B\subseteq S$ with $%
\nu \left( B\right) >0$ such that for every $n\in \mathbb{N}$ there exists $%
a_{n}\in L_{\infty }\left( \nu \right) $ satisfying $\left\vert
a_{n}\right\vert \leq n^{-2}\left\Vert q_{n}\right\Vert _{\mathcal{L}\left(
X\right) }^{-1}$ and 
\begin{equation}
\left\vert \delta ^{z}\left( a_{n}\right) \right\vert \geq n\left\Vert
q_{n}\right\Vert _{\mathcal{L}\left( X\right) }\chi _{B}.  \label{Deq02}
\end{equation}

Since $\left\vert a_{n}\left( s\right) \right\vert \leq n^{-2}\left\Vert
q_{n}\right\Vert _{\mathcal{L}\left( X\right) }^{-1}$, $s\in S$, it follows
that the series 
\begin{equation*}
f\left( s\right) =\sum_{n=1}^{\infty }a_{n}\left( s\right) q_{n}
\end{equation*}%
defines an element $f\left( s\right) \in \mathcal{K}\left( X\right) $. It is
also clear that $f\in L_{0}\left( \nu ;\mathcal{K}\left( X\right) \right)
\subseteq L_{0}\left( \nu ;\mathcal{L}\left( X\right) \right) $. Given $%
\varepsilon >0$, there exists $N\in \mathbb{N}$ such that 
\begin{equation*}
\sum_{n=N}^{\infty }\left\vert a_{n}\left( s\right) \right\vert \left\Vert
q_{n}\right\Vert _{\mathcal{L}\left( X\right) }\leq \varepsilon ,s\in S,
\end{equation*}%
and hence, the projection $p_{N}\in \mathcal{I}\left( \mathcal{F}\left(
X\right) \right) $ satisfies $\left\Vert f\left( \mathbf{1}\otimes
p_{N}^{\bot }\right) \right\Vert _{\infty }\leq \varepsilon $. Therefore, it
follows from (ii) that $f\in \mathfrak{A}$ (recall that it is assumed at
present that $\nu \left( S\right) <\infty $).

Since $\left( \mathbf{1}\otimes q_{n}\right) f=\left( a_{n}\otimes
q_{n}\right) =f\left( \mathbf{1}\otimes q_{n}\right) $, it follows that 
\begin{eqnarray*}
\left( \mathbf{1}\otimes q_{n}\right) \delta \left( f\right) \left( \mathbf{1%
}\otimes q_{n}\right) &=&\left( \mathbf{1}\otimes q_{n}\right) \delta \left(
f\left( \mathbf{1}\otimes q_{n}\right) \right) \left( \mathbf{1}\otimes
q_{n}\right) \\
&=&\left( \mathbf{1}\otimes q_{n}\right) \delta \left( a_{n}\otimes
q_{n}\right) \left( \mathbf{1}\otimes q_{n}\right) \\
&=&\delta _{q_{n}}\left( a_{n}\otimes q_{n}\right) =\delta
_{q_{n}}^{z}\left( a_{n}\right) \otimes q_{n}
\end{eqnarray*}%
for all $n\in \mathbb{N}$. Observing that $q_{n}p_{n+1}=p_{n+1}q_{n}=q_{n}$,
it follows from Lemma \ref{DLem03} that $\delta _{q_{n}}^{z}=\delta
_{p_{n+1}}^{z}$, that is, $\delta _{q_{n}}^{z}=\delta ^{z}$ on $L_{\infty
}\left( \nu \right) $. Consequently, 
\begin{equation*}
\left( \mathbf{1}\otimes q_{n}\right) \delta \left( f\right) \left( \mathbf{1%
}\otimes q_{n}\right) =\delta ^{z}\left( a_{n}\right) \otimes q_{n},\ \ \
n\in \mathbb{N}.
\end{equation*}%
If $s\in B$, then 
\begin{equation*}
\left\vert \delta ^{z}\left( a_{n}\right) \left( s\right) \right\vert
\left\Vert q_{n}\right\Vert _{\mathcal{L}\left( X\right) }=\left\Vert
q_{n}\delta \left( f\right) \left( s\right) q_{n}\right\Vert _{\mathcal{L}%
\left( X\right) }\leq \left\Vert q_{n}\right\Vert _{\mathcal{L}\left(
X\right) }^{2}\left\Vert \delta \left( f\right) \left( s\right) \right\Vert
_{\mathcal{L}\left( X\right) }
\end{equation*}%
and so, (\ref{Deq02}) implies that 
\begin{equation*}
n\left\Vert q_{n}\right\Vert _{\mathcal{L}\left( X\right) }^{2}\leq
\left\Vert q_{n}\right\Vert _{\mathcal{L}\left( X\right) }^{2}\left\Vert
\delta \left( f\right) \left( s\right) \right\Vert _{\mathcal{L}\left(
X\right) },\ \ \ n\in \mathbb{N}.
\end{equation*}%
Since $q_{n}\neq 0$, this contradicts the fact that $\left\Vert \delta
\left( f\right) \left( s\right) \right\Vert _{\mathcal{L}\left( X\right)
}<\infty $ $\nu $-a.e. on $S$. Therefore, it may be concluded that $\delta
^{z}=0$.\medskip
\end{proof}

Note that the result of Lemma \ref{DLem04} may be formulated as%
\begin{equation}
\left( \mathbf{1}\otimes p_{n}\right) \delta \left( a\otimes p_{n}\right)
\left( \mathbf{1}\otimes p_{n}\right) =0,\ \ \ a\in L_{\infty }\left( \nu
\right) ,n\in \mathbb{N}.  \label{Deq03}
\end{equation}

\begin{lemma}
\label{DLem06}Same situation as above. The derivation $\delta $ is $%
L_{\infty }\left( \nu \right) $-linear on $\mathfrak{A}$, that is, $\delta
\left( af\right) =a\delta \left( f\right) $ for all $f\in \mathfrak{A}$ and $%
a\in L_{\infty }\left( \nu \right) $.
\end{lemma}

\begin{proof}
It should be observed first that $\delta \left( a\otimes p_{n}\right)
=a\delta \left( \mathbf{1}\otimes p_{n}\right) $ for all $a\in L_{\infty
}\left( \nu \right) $ and $n\in \mathbb{N}$. Indeed, if $m\geq n$, then 
\begin{eqnarray*}
\delta \left( a\otimes p_{n}\right) &=&\delta \left( \left( \mathbf{1}%
\otimes p_{n}\right) \left( a\otimes p_{m}\right) \right) \\
&=&\delta \left( \mathbf{1}\otimes p_{n}\right) \left( a\otimes p_{m}\right)
+\left( \mathbf{1}\otimes p_{n}\right) \delta \left( a\otimes p_{m}\right)
\end{eqnarray*}%
and so, 
\begin{equation*}
\delta \left( a\otimes p_{n}\right) \left( \mathbf{1}\otimes p_{m}\right)
=a\delta \left( \mathbf{1}\otimes p_{n}\right) \left( \mathbf{1}\otimes
p_{m}\right) +\left( \mathbf{1}\otimes p_{n}\right) \delta \left( a\otimes
p_{m}\right) \left( \mathbf{1}\otimes p_{m}\right) .
\end{equation*}%
It follows from (\ref{Deq03}) that 
\begin{equation*}
\left( \mathbf{1}\otimes p_{n}\right) \delta \left( a\otimes p_{m}\right)
\left( \mathbf{1}\otimes p_{m}\right) =\left( \mathbf{1}\otimes p_{n}\right)
\left( \mathbf{1}\otimes p_{m}\right) \delta \left( a\otimes p_{m}\right)
\left( \mathbf{1}\otimes p_{m}\right) =0
\end{equation*}%
and hence, 
\begin{equation*}
\delta \left( a\otimes p_{n}\right) \left( \mathbf{1}\otimes p_{m}\right)
=a\delta \left( \mathbf{1}\otimes p_{n}\right) \left( \mathbf{1}\otimes
p_{m}\right)
\end{equation*}%
for all $m\geq n$. Now, Lemma \ref{DLem05} implies that $\delta \left(
a\otimes p_{n}\right) =a\delta \left( \mathbf{1}\otimes p_{n}\right) $.

Suppose next that $f\in \mathfrak{A}$ and $a\in L_{\infty }\left( \nu
\right) $. If $n\in \mathbb{N}$, then 
\begin{eqnarray}
\delta \left( a\left( \mathbf{1}\otimes p_{n}\right) f\right) &=&\delta
\left( a\left( \mathbf{1}\otimes p_{n}\right) \right) f+a\left( \mathbf{1}%
\otimes p_{n}\right) \delta \left( f\right)  \label{Deq04} \\
&=&a\delta \left( \mathbf{1}\otimes p_{n}\right) f+a\left( \mathbf{1}\otimes
p_{n}\right) \delta \left( f\right) .  \notag
\end{eqnarray}%
On the other hand, 
\begin{eqnarray}
\delta \left( a\left( \mathbf{1}\otimes p_{n}\right) f\right) &=&\delta
\left( \left( \mathbf{1}\otimes p_{n}\right) af\right)  \label{Deq05} \\
&=&\delta \left( \mathbf{1}\otimes p_{n}\right) af+\left( \mathbf{1}\otimes
p_{n}\right) \delta \left( af\right)  \notag \\
&=&a\delta \left( \mathbf{1}\otimes p_{n}\right) f+\left( \mathbf{1}\otimes
p_{n}\right) \delta \left( af\right) .  \notag
\end{eqnarray}%
Comparing (\ref{Deq04}) and (\ref{Deq05}), it follows that 
\begin{equation*}
\left( \mathbf{1}\otimes p_{n}\right) \left( a\delta \left( f\right) \right)
=\left( \mathbf{1}\otimes p_{n}\right) \delta \left( af\right) ,\ \ \ n\in 
\mathbb{N}.
\end{equation*}%
Consequently, Lemma \ref{DLem05} implies that $a\delta \left( f\right)
=\delta \left( af\right) $. The proof is complete.\medskip
\end{proof}

\begin{proposition}
\label{PProp03}Suppose that $X$ is an infinite dimensional separable Banach
space and that $\mathfrak{A}$ is a subalgebra of $L_{0}^{wo}\left( \nu ;%
\mathcal{L}\left( X\right) \right) $ satisfying conditions (i) and (ii) of
Definition \ref{DDef01}. If $\delta :\mathfrak{A}\rightarrow \mathfrak{A}$
is a derivation, then $\delta $ is $L_{\infty }\left( \nu \right) $-linear,
that is, $\delta \left( af\right) =a\delta \left( f\right) $ for all $f\in 
\mathfrak{A}$ and $a\in L_{\infty }\left( \nu \right) $.
\end{proposition}

\begin{proof}
The case that $\nu \left( S\right) <\infty $ has been obtained in Lemma \ref%
{DLem06}. Using Corollary \ref{DCor01}, the extension to the case that $\nu $
is $\sigma $-finite follows via a standard argument.\medskip
\end{proof}

\begin{remark}
\label{PRem01}In the situation of Proposition \ref{PProp03}, define the
subalgebra $M_{\mathfrak{A}}$ of $L_{0}\left( \nu \right) $ by 
\begin{equation*}
M_{\mathfrak{A}}=\left\{ a\in L_{0}\left( \nu \right) :af\in \mathfrak{A}%
\text{ \ for all }f\in \mathfrak{A}\right\} .
\end{equation*}%
By (i) of Definition \ref{DDef01}, $L_{\infty }\left( \nu \right) \subseteq
M_{\mathfrak{A}}$. We claim that $\delta \left( af\right) =a\delta \left(
f\right) $ for all $f\in \mathfrak{A}$ and $a\in M_{\mathfrak{A}}$. Indeed,
given $a\in M_{\mathfrak{A}}$, there exists a $\nu $-measurable partition $%
\left\{ A_{n}\right\} _{n=1}^{\infty }$ of $S$ such that $a\chi _{A_{n}}\in
L_{\infty }\left( \nu \right) $ for all $n$. Using Corollary \ref{DCor01}
and Proposition \ref{PProp03}, it follows that 
\begin{equation*}
\chi _{A_{n}}\delta \left( af\right) =\delta \left( \chi _{A_{n}}af\right)
=\chi _{A_{n}}a\delta \left( f\right)
\end{equation*}%
for all $n\in \mathbb{N}$. This implies that $\delta \left( af\right)
=a\delta \left( f\right) $.

In particular, if $\mathfrak{A=}L_{0}^{wo}\left( \nu ;\mathfrak{U}\right) $
for some closed standard subalgebra $\mathfrak{U}$ of $\mathcal{L}\left(
X\right) $, then $\delta $ is $L_{0}\left( \nu \right) $-linear.
\end{remark}

If $\mathfrak{U}$ is a closed standard subalgebra of $\mathcal{L}\left(
X\right) $, then $\mathcal{T}_{0}$ denotes the topology in $L_{0}^{wo}\left(
\nu ;\mathfrak{U}\right) $ of convergence in measure on sets of finite
measure (see Section \ref{SectVecVal}).

\begin{proposition}
\label{PProp04}Let $\mathfrak{U}$ be a closed standard subalgebra of $%
\mathcal{L}\left( X\right) $. If $\delta :L_{0}^{wo}\left( \nu ;\mathfrak{U}%
\right) \rightarrow L_{0}^{wo}\left( \nu ;\mathfrak{U}\right) $ is an $%
L_{0}\left( \nu \right) $-linear derivation, then $\delta $ is continuous
with respect to $\mathcal{T}_{0}$.
\end{proposition}

\begin{proof}
It follows from Lemma \ref{NLem01} that $L_{0}^{wo}\left( \nu ;\mathfrak{U}%
\right) $ is a complete metrizable topological vector space with respect to $%
\mathcal{T}_{0}$. Therefore, it is sufficient to show that $\delta $ is
closed (see e.g. \cite{Ru}, Theorem 2.15). Accordingly, suppose that $%
\left\{ f_{n}\right\} _{n=1}^{\infty }$ is a sequence in $L_{0}^{wo}\left(
\nu ;\mathfrak{U}\right) $ such that $f_{n}\overset{\mathcal{T}_{0}}{%
\rightarrow }0$ and $\delta \left( f_{n}\right) \overset{\mathcal{T}_{0}}{%
\rightarrow }g$ as $n\rightarrow \infty $ for some $g\in L_{0}^{wo}\left(
\nu ;\mathfrak{U}\right) $. Let $p,q\in \mathcal{I}\left( \mathcal{F}\left(
X\right) \right) $ be two rank one projections and note that $p,q\in 
\mathfrak{U}$. There exist $x_{0},y_{0}\in X$ and $x_{0}^{\ast },y_{0}^{\ast
}\in X^{\ast }$ satisfying $\left\langle x_{0},x_{0}^{\ast }\right\rangle
=\left\langle y_{0},y_{0}^{\ast }\right\rangle =1$, such that $%
px=\left\langle x,x_{0}^{\ast }\right\rangle x_{0}$ and $qx=\left\langle
x,y_{0}^{\ast }\right\rangle y_{0}$, $x\in X$. Defining the rank one
operator $v\in \mathcal{F}\left( X\right) $ by setting $vx=\left\langle
x,y_{0}^{\ast }\right\rangle x_{0}$, $x\in X$, it follows that 
\begin{equation}
pTq=\left\langle Ty_{0},x_{0}^{\ast }\right\rangle v,\ \ \ T\in \mathcal{L}%
\left( X\right) .  \label{Peq03}
\end{equation}%
It follows from (\ref{Peq03}) that 
\begin{equation*}
\left( \left( \mathbf{1}\otimes p\right) f_{n}\left( \mathbf{1}\otimes
q\right) \right) \left( s\right) =pf_{n}\left( s\right) q=a_{n}\left(
s\right) v,\ \ \ s\in S,
\end{equation*}%
where $a_{n}\in L_{0}\left( \nu \right) $. Since $f_{n}\overset{\mathcal{T}%
_{0}}{\rightarrow }0$ implies that $\left( \mathbf{1}\otimes p\right)
f_{n}\left( \mathbf{1}\otimes q\right) \overset{\mathcal{T}_{0}}{\rightarrow 
}0$, it is clear that $a_{n}\rightarrow 0$ as $n\rightarrow \infty $ in
measure on sets of finite measure. Hence, $a_{n}\delta \left( v\right) 
\overset{\mathcal{T}_{0}}{\rightarrow }0$ in $L_{0}^{wo}\left( \nu ;%
\mathfrak{U}\right) $. Furthermore, $\delta $ is $L_{0}$-linear and so, 
\begin{equation*}
\delta \left( \left( \mathbf{1}\otimes p\right) f_{n}\left( \mathbf{1}%
\otimes q\right) \right) =a_{n}\delta \left( v\right)
\end{equation*}%
for all $n$. Consequently, $\delta \left( \left( \mathbf{1}\otimes p\right)
f_{n}\left( \mathbf{1}\otimes q\right) \right) \overset{\mathcal{T}_{0}}{%
\rightarrow }0$ as $n\rightarrow \infty $. Since, 
\begin{multline*}
\delta \left( \left( \mathbf{1}\otimes p\right) f_{n}\left( \mathbf{1}%
\otimes q\right) \right) \\
=\delta \left( \mathbf{1}\otimes p\right) f_{n}\left( \mathbf{1}\otimes
q\right) +\left( \mathbf{1}\otimes p\right) \delta \left( f_{n}\right)
\left( \mathbf{1}\otimes q\right) +\left( \mathbf{1}\otimes p\right)
f_{n}\delta \left( \mathbf{1}\otimes q\right) ,
\end{multline*}%
$\delta \left( \mathbf{1}\otimes p\right) f_{n}\left( \mathbf{1}\otimes
q\right) \overset{\mathcal{T}_{0}}{\rightarrow }0$ and $\left( \mathbf{1}%
\otimes p\right) f_{n}\delta \left( \mathbf{1}\otimes q\right) \overset{%
\mathcal{T}_{0}}{\rightarrow }0$ as $n\rightarrow \infty $, it follows that 
\begin{equation*}
\left( \mathbf{1}\otimes p\right) \delta \left( f_{n}\right) \left( \mathbf{1%
}\otimes q\right) \overset{\mathcal{T}_{0}}{\rightarrow }0,\ \ \
n\rightarrow \infty .
\end{equation*}%
On the other hand, $\delta \left( f_{n}\right) \overset{\mathcal{T}_{0}}{%
\rightarrow }g$ implies that 
\begin{equation*}
\left( \mathbf{1}\otimes p\right) \delta \left( f_{n}\right) \left( \mathbf{1%
}\otimes q\right) \overset{\mathcal{T}_{0}}{\rightarrow }\left( \mathbf{1}%
\otimes p\right) g\left( \mathbf{1}\otimes q\right) ,\ \ \ n\rightarrow
\infty
\end{equation*}%
and so, $\left( \mathbf{1}\otimes p\right) g\left( \mathbf{1}\otimes
q\right) =0$.

Since every finite rank projection in $X$ may be written as a sum of
finitely many rank one projections, it follows that $\left( \mathbf{1}%
\otimes p\right) g\left( \mathbf{1}\otimes p\right) =0$ for all finite rank
projections $p$ in $X$. Consequently, via Lemma \ref{DLem05} it may be
concluded that $g=0$. This shows that $\delta $ is closed and hence, $\delta 
$ is $\mathcal{T}_{0}$-continuous. \medskip
\end{proof}

\begin{corollary}
\label{PProp041}Assuming that $X$ is infinite dimensional, let $\mathfrak{U}$
be a closed standard subalgebra of $\mathcal{L}\left( X\right) $. If $\delta
:L_{0}^{wo}\left( \nu ;\mathfrak{U}\right) \rightarrow L_{0}^{wo}\left( \nu ;%
\mathfrak{U}\right) $ is a derivation, then $\delta $ is continuous with
respect to $\mathcal{T}_{0}$.
\end{corollary}

\begin{proof}
It follows from Remark \ref{PRem01} that $\delta $ is $L_{0}$-linear and so,
the result is an immediate consequence of \ref{PProp04}.
\end{proof}

\begin{proposition}
\label{PProp06}Let $\mathfrak{U}$ be a closed standard subalgebra of $%
\mathcal{L}\left( X\right) $. If $\delta :L_{0}^{wo}\left( \nu ;\mathfrak{U}%
\right) \rightarrow L_{0}^{wo}\left( \nu ;\mathfrak{U}\right) $ is an $%
L_{0}\left( \nu \right) $-linear derivation, then there exists $w\in
L_{0}\left( \nu ;\mathcal{L}\left( X\right) \right) $ such that $\delta
\left( f\right) =wf-fw$ for all $f\in \mathfrak{A}$.
\end{proposition}

\begin{proof}
It follows from Proposition \ref{PProp04} that $\delta $ is continuous with
respect to $\mathcal{T}_{0}$. Consequently, by Proposition \ref{NProp01},
there exists $0\leq a\in L_{0}\left( \nu \right) $ such that 
\begin{equation}
\left\Vert \delta \left( f\right) \left( s\right) \right\Vert _{\mathcal{L}%
\left( X\right) }\leq a\left( s\right) \left\Vert f\left( s\right)
\right\Vert _{\mathcal{L}\left( X\right) },\ \ \ \nu \text{-a.e. on }S,
\label{eqB05}
\end{equation}%
for all $f\in L_{0}^{wo}\left( \nu ;\mathfrak{U}\right) $. Fix $x_{0}\in X$
with $\left\Vert x_{0}\right\Vert _{X}=1$ and let $x_{0}^{\ast }\in X^{\ast
} $ be such that $\left\Vert x_{0}^{\ast }\right\Vert _{X^{\ast }}=1$ and $%
\left\langle x_{0},x_{0}^{\ast }\right\rangle =1$. Define $F_{0}\in
L_{0}\left( \nu ;X\right) $ by $F_{0}=\mathbf{1}\otimes x_{0}$.

For $F\in L_{0}\left( \nu ;X\right) $ and $s\in S$, define the rank one
operator $g_{F}\left( s\right) \in \mathcal{F}\left( X\right) \subseteq 
\mathfrak{U}$ by setting 
\begin{equation*}
g_{F}\left( s\right) x=\left\langle x,x_{0}^{\ast }\right\rangle F\left(
s\right) ,\ \ \ x\in X.
\end{equation*}%
It is clear that 
\begin{equation}
\left\Vert g_{F}\left( s\right) \right\Vert _{\mathcal{L}\left( X\right)
}\leq \left\Vert F\left( s\right) \right\Vert _{X}.  \label{eqB04}
\end{equation}%
Furthermore, if $x\in X$ and $x^{\ast }\in X^{\ast }$, then the function 
\begin{equation*}
s\longmapsto \left\langle g_{F}\left( s\right) x,x^{\ast }\right\rangle
=\left\langle x,x_{0}^{\ast }\right\rangle \left\langle F\left( s\right)
,x^{\ast }\right\rangle ,\ \ \ s\in S,
\end{equation*}%
is $\nu $-measurable. Therefore, $g_{F}\in L_{0}^{wo}\left( \nu ;\mathfrak{U}%
\right) $. Identifying the elements of $L_{0}^{wo}\left( \nu ;\mathfrak{U}%
\right) $ with $L_{0}\left( \nu \right) $-linear operators on $L_{0}\left(
\nu ;X\right) $ (see Corollary \ref{Cor04}), it is clear that $g_{F}F_{0}=F$%
. Furthermore, the map $F\longmapsto g_{F}$, $F\in L_{0}\left( \nu ;X\right) 
$, is linear from $L_{0}\left( \nu ;X\right) $ into $L_{0}^{wo}\left( \nu ;%
\mathfrak{U}\right) $. It should be observed that 
\begin{equation}
g_{fF}=fg_{F},\ \ \ f\in L_{0}^{wo}\left( \nu ;\mathfrak{U}\right) ,F\in
L_{0}\left( \nu ;X\right) .  \label{eqB01}
\end{equation}%
Indeed, 
\begin{equation*}
g_{fF}\left( s\right) x=\left\langle x,x_{0}^{\ast }\right\rangle f\left(
s\right) F\left( s\right) =f\left( s\right) \left( \left\langle
x,x_{0}^{\ast }\right\rangle F\left( s\right) \right) =f\left( s\right)
g_{F}\left( s\right) x
\end{equation*}%
for $x\in X$, $s\in S$, and so, (\ref{eqB01}) holds. Similarly, it follows
that 
\begin{equation}
g_{bF}=bg_{F},\ \ \ b\in L_{0}\left( \nu \right) ,F\in L_{0}\left( \nu
;X\right) ,  \label{Peq04}
\end{equation}%
that is, the map $F\longmapsto g_{F}$, $F\in L_{0}\left( \nu ;X\right) $, is 
$L_{0}\left( \nu \right) $-linear.

Define the linear operator $T:L_{0}\left( \nu ;X\right) \rightarrow
L_{0}\left( \nu ;X\right) $ by setting 
\begin{equation}
TF=\delta \left( g_{F}\right) F_{0},\ \ \ F\in L_{0}\left( \nu ;X\right) .
\label{Peq05}
\end{equation}%
Since $\delta $ is $L_{0}$-linear, it follows, in particularly, from (\ref%
{Peq04}) that $T$ is also $L_{0}$-linear. Furthermore, it follows from (\ref%
{eqB05}) and (\ref{eqB04}) that 
\begin{eqnarray*}
\left\Vert \left( TF\right) \left( s\right) \right\Vert _{X} &=&\left\Vert
\delta \left( g_{F}\right) \left( s\right) x_{0}\right\Vert _{X}\leq
\left\Vert \delta \left( g_{F}\right) \left( s\right) \right\Vert _{\mathcal{%
L}\left( X\right) } \\
&\leq &a\left( s\right) \left\Vert g_{F}\left( s\right) \right\Vert _{%
\mathcal{L}\left( X\right) }\leq a\left( s\right) \left\Vert F\left(
s\right) \right\Vert _{X},\ \ \ s\in S,
\end{eqnarray*}%
and so, $T$ is continuous with respect to the topology $\mathcal{T}_{0}$
(see Proposition \ref{NProp01}). Therefore, it follows from Proposition \ref%
{PProp05} that there exists $w\in L_{0}^{wo}\left( \nu ;\mathcal{L}\left(
X\right) \right) $ such that 
\begin{equation}
TF=wF,\ \ \ F\in L_{0}\left( \nu ;X\right) .  \label{Peq06}
\end{equation}%
Finally, it will be shown that $\delta \left( f\right) =\left[ w,f\right] $
for all $f\in L_{0}^{wo}\left( \nu ;\mathfrak{U}\right) $. If $f\in
L_{0}^{wo}\left( \nu ;\mathfrak{U}\right) $ and $F\in L_{0}\left( \nu
;X\right) $, then it follows from (\ref{Peq05}), (\ref{Peq06}) and (\ref%
{eqB01}) that 
\begin{eqnarray*}
\left[ w,f\right] F &=&wfF-fwF=\delta \left( g_{fF}\right) F_{0}-f\delta
\left( g_{F}\right) F_{0} \\
&=&\delta \left( fg_{F}\right) F_{0}-f\delta \left( g_{F}\right) F_{0} \\
&=&\delta \left( f\right) g_{F}F_{0}+f\delta \left( g_{F}\right)
F_{0}-f\delta \left( g_{F}\right) F_{0} \\
&=&\delta \left( f\right) g_{F}F_{0}=\delta \left( f\right) F.
\end{eqnarray*}%
Since this holds for all $F\in L_{0}\left( \nu ;X\right) $, we may conclude
that $\delta \left( f\right) =\left[ w,f\right] $. The proof is
complete.\medskip
\end{proof}

\begin{corollary}
\label{PProp061}Suppose that $X$ is infinite dimensional and let $\mathfrak{U%
} $ be a closed standard subalgebra of $\mathcal{L}\left( X\right) $. If $%
\delta :L_{0}^{wo}\left( \nu ;\mathfrak{U}\right) \rightarrow
L_{0}^{wo}\left( \nu ;\mathfrak{U}\right) $ is a derivation, then there
exists $w\in L_{0}\left( \nu ;\mathcal{L}\left( X\right) \right) $ such that 
$\delta \left( f\right) =wf-fw$ for all $f\in \mathfrak{A}$.
\end{corollary}

\begin{corollary}
Let $X$ be a separable Banach space and let $\mathcal{A}$ be the algebra $%
L_{0}^{wo}\left( \nu ;\mathcal{L}\left( X\right) \right) $ or $%
L_{0}^{wo}\left( \nu ;\mathcal{K}\left( X\right) \right) $. If $\delta :%
\mathcal{A}\rightarrow \mathcal{A}$ is an $L_{0}\left( \nu \right) $-linear
derivation, then there exists $w\in L_{0}^{wo}\left( \nu ;\mathcal{L}\left(
X\right) \right) $ such that $\delta \left( f\right) =wf-fw$ for all $f\in 
\mathcal{A}$.
\end{corollary}

\begin{corollary}
Let $X$ be an infinite dimensional separable Banach space and let $\mathcal{A%
}$ be the algebra $L_{0}^{wo}\left( \nu ;\mathcal{L}\left( X\right) \right) $
or $L_{0}^{wo}\left( \nu ;\mathcal{K}\left( X\right) \right) $. If $\delta :%
\mathcal{A}\rightarrow \mathcal{A}$ is a derivation, then there exists $w\in
L_{0}^{wo}\left( \nu ;\mathcal{L}\left( X\right) \right) $ such that $\delta
\left( f\right) =wf-fw$ for all $f\in \mathcal{A}$.
\end{corollary}

\begin{theorem}
\label{PThm01}Let $X$ be a separable Banach space and $\left( S,\Sigma ,\nu
\right) $ be a $\sigma $-finite measure space. Suppose that $\mathfrak{U}%
\subseteq \mathcal{L}\left( X\right) $ is a closed standard algebra and that 
$\mathfrak{A}\subseteq L_{0}^{wo}\left( \nu ;\mathcal{L}\left( X\right)
\right) $ is a $\mathfrak{U}$-admissible subalgebra. If $\delta :\mathfrak{A}%
\rightarrow \mathfrak{A}$ is an $L_{\infty }\left( \nu \right) $-linear
derivation, then there exists $w\in L_{0}^{wo}\left( \nu ;\mathcal{L}\left(
X\right) \right) $ such that $\delta \left( f\right) =wf-fw$ for all $f\in 
\mathfrak{A}$.
\end{theorem}

\begin{proof}
It follows from condition (iii) of Definition \ref{DDef01} and Proposition %
\ref{PProp01} that $\delta $ extends uniquely to an $L_{0}\left( \nu \right) 
$-linear map $\hat{\delta}:L_{0}^{wo}\left( \nu ;\mathfrak{U}\right)
\rightarrow L_{0}^{wo}\left( \nu ;\mathfrak{U}\right) $. Recall from
Proposition \ref{PProp01} that $\chi _{A}\hat{\delta}\left( f\right) =\delta
\left( \chi _{A}f\right) $ for all $f\in L_{0}^{wo}\left( \nu ;\mathfrak{U}%
\right) $ and all $\nu $-measurable sets $A$ satisfying $\chi _{A}f\in 
\mathfrak{A}$. We claim that $\hat{\delta}$ is a derivation. Indeed, suppose
that $f,g\in L_{0}^{wo}\left( \nu ;\mathfrak{U}\right) $. It follows from
condition (iii) on the algebra $\mathfrak{A}$ that there exists a measurable
partition $\left\{ A_{n}\right\} _{n=1}^{\infty }$ such that $\chi
_{A_{n}}f,\chi _{A_{n}}g\in \mathfrak{A}$ for all $n$. Using that $\hat{%
\delta}$ is $L_{0}\left( \nu \right) $-linear, we find that%
\begin{eqnarray*}
\chi _{A_{n}}\hat{\delta}\left( fg\right) &=&\hat{\delta}\left( \chi
_{A_{n}}fg\right) =\delta \left( \chi _{A_{n}}f\chi _{A_{n}}g\right) \\
&=&\delta \left( \chi _{A_{n}}f\right) \chi _{A_{n}}g+\chi _{A_{n}}f\delta
\left( \chi _{A_{n}}g\right) \\
&=&\chi _{A_{n}}\hat{\delta}\left( f\right) g+\chi _{A_{n}}f\hat{\delta}%
\left( g\right)
\end{eqnarray*}%
for all $n\in \mathbb{N}$. Hence, $\hat{\delta}\left( fg\right) =\hat{\delta}%
\left( f\right) g+f\hat{\delta}\left( g\right) $, which shows that $\hat{%
\delta}$ is a derivation in $L_{0}^{wo}\left( \nu ;\mathfrak{U}\right) $.
Consequently, the result of the theorem follows immediately from Proposition %
\ref{PProp06}. \medskip
\end{proof}

\begin{corollary}
\label{PThm011}Let $X$ be an infinite dimensional separable Banach space and 
$\left( S,\Sigma ,\nu \right) $ be a $\sigma $-finite measure space. Suppose
that $\mathfrak{U}\subseteq \mathcal{L}\left( X\right) $ is a closed
standard algebra and that $\mathfrak{A}\subseteq L_{0}^{wo}\left( \nu ;%
\mathcal{L}\left( X\right) \right) $ is a $\mathfrak{U}$-admissible
subalgebra. If $\delta :\mathfrak{A}\rightarrow \mathfrak{A}$ is a
derivation, then there exists $w\in L_{0}^{wo}\left( \nu ;\mathcal{L}\left(
X\right) \right) $ such that $\delta \left( f\right) =wf-fw$ for all $f\in 
\mathfrak{A}$.
\end{corollary}

\section{Applications to algebras of measurable operators\label{SectMeas}}

In the present section, we shall apply the results of the previous section
to algebras of measurable operators with respect to certain von Neumann
algebras. For details on von Neumann algebra theory, the reader is referred
to e.g. \cite{Dix}, \cite{KR}, \cite{Sak1} or \cite{Ta}. General facts
concerning measurable operators may be found in \cite{Se}, \cite{Te} and 
\cite{Ye} (see also \cite{Ta2}, Chapter IX). For the convenience of the
reader, some of the basic definitions are recalled.

Let $\mathcal{M}$ be a von Neumann algebra on a Hilbert space $H$. The set
of all self-adjoint projections in $\mathcal{M}$ is denoted by $P\left( 
\mathcal{M}\right) $. The commutant of $\mathcal{M}$ is denoted by $\mathcal{%
M}^{\prime }$. A linear operator $x:\mathfrak{D}\left( x\right) \rightarrow
H $, where the domain $\mathfrak{D}\left( x\right) $ of $x$ is a linear
subspace of $H$, is said to be affiliated with $\mathcal{M}$ if $yx\subseteq
xy$ for all $y\in \mathcal{M}^{\prime }$ (which is denoted by $x\eta 
\mathcal{M}$). A linear operator $x:\mathfrak{D}\left( x\right) \rightarrow
H $ is termed measurable with respect to $\mathcal{M}$ if $x$ is closed,
densely defined, affiliated with $\mathcal{M}$ and there exists a sequence $%
\left\{ p_{n}\right\} _{n=1}^{\infty }$ in $P\left( \mathcal{M}\right) $
such that $p_{n}\uparrow \mathbf{1}$, $p_{n}\left( H\right) \subseteq 
\mathfrak{D}\left( x\right) $ and $p_{n}^{\bot }$ is a finite projection
(with respect to $\mathcal{M}$) for all $n$. It should be noted that the
condition $p_{n}\left( H\right) \subseteq \mathfrak{D}\left( x\right) $
implies that $xp_{n}\in \mathcal{M}$. The collection of all measurable
operators with respect to $\mathcal{M}$ is denoted by $S\left( \mathcal{M}%
\right) $, which is a unital $\ast $-algebra with respect to strong sums and
products (denoted simply by $x+y$ and $xy$ for all $x,y\in S\left( \mathcal{M%
}\right) $). Furthermore, a closed densely defined linear operator $x:%
\mathfrak{D}\left( x\right) \rightarrow H$ is called locally measurable
(with respect to $\mathcal{M}$) if $x\eta \mathcal{M}$ and there exists a
sequence $\left\{ p_{n}\right\} _{n=1}^{\infty }$ in $P\left( \mathcal{Z}%
\left( \mathcal{M}\right) \right) $ such that $p_{n}\uparrow \mathbf{1}$ and 
$xp_{n}\in S\left( \mathcal{M}\right) $ for all $n$. The collection of all
locally measurable operators is denoted by $LS\left( \mathcal{M}\right) $,
which is a unital $\ast $-algebra with respect to strong sums and products.
Evidently, $S\left( \mathcal{M}\right) $ is a $\ast $-subalgebra of $%
LS\left( \mathcal{M}\right) $. For details we refer e.g. to \cite{Ye}.

If the von Neumann $\mathcal{M}$ is semi-finite, equipped with a semi-finite
faithful normal trace $\tau $, then an operator $x\in S\left( \mathcal{M}%
\right) $ is called $\tau $-measurable if there exists a sequence $\left\{
p_{n}\right\} _{n=1}^{\infty }$ in $P\left( \mathcal{M}\right) $ such that $%
p_{n}\uparrow \mathbf{1}$, $p_{n}\left( H\right) \subseteq \mathfrak{D}%
\left( x\right) $ and $\tau \left( p_{n}^{\bot }\right) <\infty $ for all $n$%
. The collection $S\left( \tau \right) $ of all $\tau $-measurable operators
is a unital $\ast $-subalgebra of $S\left( \mathcal{M}\right) $. For details
see e.g. \cite{Te}. An operator $x\in S\left( \tau \right) $ is said to $%
\tau $-compact if for every $\varepsilon >0$ there exists a projection $p\in
P\left( \mathcal{M}\right) $ such that $\left\Vert xp\right\Vert _{B\left(
H\right) }<\varepsilon $ and $\tau \left( p^{\bot }\right) <\infty $. The
collection $S_{0}\left( \tau \right) $ of all $\tau $-compact operators is a
two-sided ideal in $S\left( \tau \right) $.

Now, we specialize to the case where $H$ is a separable Hilbert space and $%
\mathcal{M}=L_{\infty }\left( \nu \right) \overline{\otimes }B\left(
H\right) $, the von Neumann algebra tensor product of $L_{\infty }\left( \nu
\right) $ and $B\left( H\right) $ acting on the Hilbert space tensor product 
$L_{2}\left( \nu \right) \overline{\otimes }H$. Here, $\left( S,\Sigma ,\nu
\right) $ is a $\sigma $-finite measure space and $L_{\infty }\left( \nu
\right) $ is considered as a von Neumann algebra on the Hilbert space $%
L_{2}\left( \nu \right) $, acting via multiplication. The von Neumann
algebra $B\left( H\right) $ is equipped with its standard trace $\func{tr}$
and $L_{\infty }\left( \nu \right) $ is equipped with the trace $\nu $,
given by $\nu \left( f\right) =\int_{S}fd\nu $, $0\leq f\in L_{\infty
}\left( \nu \right) $. Let $\tau =\nu \otimes \func{tr}$ be the tensor
product trace on $\mathcal{M}$, that is, $\tau $ is the unique semi-finite
normal faithful trace $\tau $ satisfying $\tau \left( f\otimes x\right) =\nu
\left( f\right) \func{tr}\left( x\right) $ for all $0\leq f\in L_{\infty
}\left( \nu \right) $ and $0\leq x\in B\left( H\right) $. Identifying the
tensor product $L_{2}\left( \nu \right) \overline{\otimes }H$ with the
Bochner space $L_{2}\left( \nu ;H\right) $, the von Neumann algebra $%
L_{\infty }\left( \nu \right) \overline{\otimes }B\left( H\right) $ may be
identified with the space $L_{\infty }^{wo}\left( \nu ;B\left( H\right)
\right) $, where the function $f\in L_{\infty }^{wo}\left( \nu ;B\left(
H\right) \right) $ corresponds to the operator $x_{f}$ on $L_{2}\left( \nu
;H\right) $ given by $\left( x_{f}F\right) \left( s\right) =f\left( s\right)
F\left( s\right) $, $\nu $-a.e. on $S$ for all $F\in L_{2}\left( \nu
;H\right) $ (see e.g. \cite{Ta}). It should be observed that if $f\in
L_{\infty }^{wo}\left( \nu ;B\left( H\right) \right) $, then $f\geq 0$ (that
is, $x_{f}\geq 0$) if and only if $f\left( s\right) \geq 0$ in $B\left(
H\right) $ $\nu $-a.e. on $S$. Using this observation, it is not difficult
to show that a sequence $\left\{ f_{n}\right\} _{n=1}^{\infty }$ in $%
L_{\infty }^{wo}\left( \nu ;B\left( H\right) \right) $ satisfies $%
f_{n}\downarrow _{n}0$ (in the von Neumann algebra $L_{\infty }^{wo}\left(
\nu ;B\left( H\right) \right) $) if and only if $f_{n}\left( s\right)
\downarrow _{n}0$ (in $B\left( H\right) $) $\nu $-a.e. on $S$. The (tensor
product) trace on $L_{\infty }^{wo}\left( \nu ;B\left( H\right) \right) $ is
given by 
\begin{equation}
\tau \left( f\right) =\int_{S}\func{tr}\left( f\left( s\right) \right) d\nu
\left( s\right) ,\ \ \ 0\leq f\in L_{\infty }^{wo}\left( \nu ;B\left(
H\right) \right) .  \label{Peq10}
\end{equation}%
Indeed, using the preceding observations, it is not difficult to show that (%
\ref{Peq10}) defines a semi-finite normal faithful trace $\tau $ on $%
L_{\infty }^{wo}\left( \nu ;B\left( H\right) \right) $ satisfying $\tau
\left( f\otimes x\right) =\nu \left( f\right) \func{tr}\left( x\right) $, $%
0\leq f\in L_{\infty }\left( \nu \right) $, $0\leq x\in B\left( H\right) $.

The following observation will be used.

\begin{lemma}
\label{PLem08}A projection $p\in \mathcal{M}\cong L_{\infty }^{wo}\left( \nu
;B\left( H\right) \right) $ is finite if and only if $p\left( s\right) $ is
a finite projection in $B\left( H\right) $ for $\nu $-almost all $s\in S$.
\end{lemma}

\begin{proof}
First assume that $p\left( s\right) $ is a finite projection in $B\left(
H\right) $ for $\nu $-almost all $s\in S$ and suppose that $q\in P\left(
L_{\infty }^{wo}\left( \nu ;B\left( H\right) \right) \right) $ is such that $%
q\leq p$ and $q\sim p$. Since $q\sim p$, there exists a partial isometry $%
v\in L_{\infty }^{wo}\left( \nu ;B\left( H\right) \right) $ such that $%
v^{\ast }v=p$ and $vv^{\ast }=q$. This implies that $v\left( s\right) ^{\ast
}v\left( s\right) =p\left( s\right) $ and $v\left( s\right) v\left( s\right)
^{\ast }=q\left( s\right) $ $\nu $-a.e. on $S$. Since $p\left( s\right) $ is
finite $\nu $-a.e., it follows that $q\left( s\right) =p\left( s\right) $ $%
\nu $-a.e. on $S$ and hence, $p=q$. This shows that the projection $p$ is
finite in $L_{\infty }^{wo}\left( \nu ;B\left( H\right) \right) $.

Assume now that $p\in P\left( L_{\infty }^{wo}\left( \nu ;B\left( H\right)
\right) \right) $ is finite and let $H\left( s\right) =p\left( s\right)
\left( H\right) $, $s\in S$. Setting $S_{\infty }=\left\{ s\in S:\dim
H\left( s\right) =\infty \right\} $, suppose that $\nu \left( S_{\infty
}\right) >0$ and define $p_{1}\in P\left( L_{\infty }^{wo}\left( \nu
;B\left( H\right) \right) \right) $ by $p_{1}=\chi _{S_{\infty }}p$. It
follows from \cite{Dix}, Proposition II.1.9 and Proposition II.1.1 that
there exist Bochner $\nu $-measurable functions $\eta _{n}:S_{\infty
}\rightarrow H$, $n\in \mathbb{N}$, such that $\left\{ \eta _{n}\left(
s\right) \right\} _{n=1}^{\infty }$ is an orthonormal basis in $H\left(
s\right) $ for all $s\in S_{\infty }$. Defining $q\in P\left( L_{\infty
}^{wo}\left( \nu ;B\left( H\right) \right) \right) $ by setting 
\begin{equation*}
q\left( s\right) \xi =\sum_{n=1}^{\infty }\left\langle \xi ,\eta _{2n}\left(
s\right) \right\rangle \eta _{2n}\left( s\right) ,\ \ \ s\in S_{\infty }
\end{equation*}%
and $q\left( s\right) =0$ if $s\in S\diagdown S_{\infty }$, it is easily
verified that $q<p_{1}$ and $q\sim p_{1}$. Hence, $p_{1}$ is not finite.
Since $p_{1}\leq p$, this is a contradiction. Therefore, we may conclude
that $\nu \left( S_{\infty }\right) =0$, that is, $p\left( s\right) $ is
finite for $\nu $-almost all $s\in S$. The proof is complete.\medskip
\end{proof}

We will show next that the algebra $L_{0}^{wo}\left( \nu ;B\left( H\right)
\right) $ may be identified with the algebra $LS\left( \mathcal{M}\right) $.
Let $f\in L_{0}^{wo}\left( \nu ;B\left( H\right) \right) $ be given. It
should be recalled that, by Lemma \ref{Lem01}, $fF\in L_{0}^{wo}\left( \nu
;H\right) $ for all $F\in L_{2}\left( \nu ;H\right) $. Define the linear
subspace $\mathfrak{D}\left( x_{f}\right) $ of $L_{2}\left( \nu ;H\right) $
by setting 
\begin{equation*}
\mathfrak{D}\left( x_{f}\right) =\left\{ F\in L_{2}\left( \nu ;H\right)
:fF\in L_{2}\left( \nu ;H\right) \right\}
\end{equation*}%
and define the linear operator $x_{f}:\mathfrak{D}\left( x_{f}\right)
\rightarrow L_{2}\left( \nu ;H\right) $ by 
\begin{equation}
x_{f}F=fF,\ \ \ F\in \mathfrak{D}\left( x_{f}\right) .  \label{Peq07}
\end{equation}

\begin{lemma}
If $f\in L_{0}^{wo}\left( \nu ;B\left( H\right) \right) $, then the operator 
$x_{f}$, defined by (\ref{Peq07}), is locally measurable with respect to $%
\mathcal{M}$.
\end{lemma}

\begin{proof}
Define the subsets $A_{n}$ of $S$ by setting%
\begin{equation}
A_{n}=\left\{ s\in S:\left\Vert f\left( s\right) \right\Vert _{B\left(
H\right) }\leq n\right\} ,\ \ \ n\in \mathbb{N},  \label{Peq08}
\end{equation}%
which are $\nu $-measurable by Lemma \ref{NLem04}. It is clear that 
\begin{equation*}
\chi _{A_{n}}L_{2}\left( \nu ;H\right) =\left\{ \chi _{A_{n}}F:F\in
L_{2}\left( \nu ;H\right) \right\} \subseteq \mathfrak{D}\left( x_{f}\right)
\end{equation*}%
and it follows from the dominated convergence theorem that 
\begin{equation}
\mathcal{D}=\bigcup\nolimits_{n=1}^{\infty }\chi _{A_{n}}L_{2}\left( \nu
;H\right)  \label{Peq12}
\end{equation}%
is dense in $L_{2}\left( \nu ;H\right) $. This shows that $\mathfrak{D}%
\left( x_{f}\right) $ is dense in $L_{2}\left( \nu ;H\right) $. Suppose that 
$\left\{ F_{n}\right\} _{n=1}^{\infty }$ is a sequence in $\mathfrak{D}%
\left( x_{f}\right) $ and that $F,G\in L_{2}\left( \nu ;H\right) $ such that 
$F_{n}\rightarrow F$ and $x_{f}F_{n}\rightarrow G$ in $L_{2}\left( \nu
;H\right) $. Passing to a subsequence, it may be assumed that $F_{n}\left(
s\right) \rightarrow F\left( s\right) $ and $f\left( s\right) F_{n}\left(
s\right) \rightarrow G\left( s\right) $ $\nu $-a.e. on $S$. Therefore, $%
f\left( s\right) F\left( s\right) =G\left( s\right) $ $\nu $-a.e. and hence, 
$F\in \mathfrak{D}\left( x_{f}\right) $ and $x_{f}F=G$. This shows that $%
x_{f}$ is closed.

Using that $\mathcal{M}^{\prime }=L_{\infty }\left( \nu \right) \otimes 
\mathbb{C}\mathbf{1}$, it is easy to verify that $x_{f}$ is affiliated with $%
\mathcal{M}$. Finally, defining the projections $\left\{ p_{n}\right\}
_{n=1}^{\infty }$ in $P\left( \mathcal{Z}\left( \mathcal{M}\right) \right) $
by $p_{n}=\chi _{A_{n}}\otimes \mathbf{1}$, it is evident that $%
p_{n}\uparrow \mathbf{1}$ and $x_{f}p_{n}\in \mathcal{M}$ for all $n\in 
\mathbb{N}$. Consequently, $x_{f}\in LS\left( \mathcal{M}\right) $.\medskip
\end{proof}

If $f\in L_{0}^{wo}\left( \nu ;B\left( H\right) \right) $ and the sets $%
A_{n}\subseteq S$ are defined by (\ref{Peq08}), then it follows easily that
the subspace $\mathcal{D}$, defined by (\ref{Peq12}), is a core for the
operator $x_{f}$, that is, $x_{f}$ is the closure of its restriction to $%
\mathcal{D}$. Using this observation, the first part of the following
proposition follows via a standard argument.

\begin{proposition}
\label{PProp08}The map $f\longmapsto x_{f}$, $f\in L_{0}^{wo}\left( \nu
;B\left( H\right) \right) $, is a unital $\ast $-isomorphism from $%
L_{0}^{wo}\left( \nu ;B\left( H\right) \right) $ onto $LS\left( \mathcal{M}%
\right) $.
\end{proposition}

\begin{proof}
It remains to be shown that for every $x\in LS\left( \mathcal{M}\right) $
there exists $f\in L_{0}^{wo}\left( \nu ;B\left( H\right) \right) $ such
that $x=x_{f}$. To this end, first assume that $x\in S\left( \mathcal{M}%
\right) $. By definition, there exists a sequence $\left\{ p_{n}\right\}
_{n=1}^{\infty }$ such that $p_{n}\uparrow \mathbf{1}$, $xp_{n}\in \mathcal{M%
}$ and $p_{n}^{\bot }$ is finite for all $n$. As observed in Lemma \ref%
{PLem08}, this implies that $p_{n}^{\bot }\left( s\right) $ is a finite
projection in $B\left( H\right) $ $\nu $-a.e., that is, $\func{tr}\left(
p_{n}^{\bot }\left( s\right) \right) \in \mathbb{N}\cup \left\{ 0\right\} $,
for $\nu $-almost all $s\in S$. Since $p_{n}^{\bot }\downarrow 0$ in $%
\mathcal{M}$ implies that $p_{n}^{\bot }\left( s\right) \downarrow 0$ in $%
B\left( H\right) $ $\nu $-a.e. on $S$ (see the observations preceding (\ref%
{Peq10})) and $\func{tr}\left( p_{n}^{\bot }\left( s\right) \right) \in 
\mathbb{N}\cup \left\{ 0\right\} $, it follows that for ($\nu $-almost)
every $s\in S$ there exists $n_{s}\in \mathbb{N}$ such that $p_{n}\left(
s\right) =\mathbf{1}$ for all $n\geq n_{s}$. Consequently, if the $\nu $%
-measurable sets $A_{n}$ are defined by setting 
\begin{equation*}
A_{n}=\left\{ s\in S:p_{n}\left( s\right) =\mathbf{1}\right\} ,
\end{equation*}%
then $A_{n}\uparrow S$. Note that $p_{n}\left( \chi _{A_{n}}\otimes \mathbf{1%
}\right) =\chi _{A_{n}}\otimes \mathbf{1}$. For each $n\in \mathbb{N}$,
there exists a unique $f_{n}\in L_{\infty }^{wo}\left( \nu ;B\left( H\right)
\right) $ such that $xp_{n}=x_{f_{n}}$. If $m\leq n$, then 
\begin{eqnarray*}
x_{\chi _{A_{m}}f_{m}} &=&x_{f_{m}}\left( \chi _{A_{m}}\otimes \mathbf{1}%
\right) =xp_{m}\left( \chi _{A_{m}}\otimes \mathbf{1}\right) =x\left( \chi
_{A_{m}}\otimes \mathbf{1}\right) \\
&=&x\left( \chi _{A_{n}}\otimes \mathbf{1}\right) \left( \chi
_{A_{m}}\otimes \mathbf{1}\right) =xp_{n}\left( \chi _{A_{n}}\otimes \mathbf{%
1}\right) \left( \chi _{A_{m}}\otimes \mathbf{1}\right) \\
&=&x_{f_{n}}\left( \chi _{A_{m}}\otimes \mathbf{1}\right) =x_{\chi
_{A_{m}}f_{n}}
\end{eqnarray*}%
and so, $\chi _{A_{m}}f_{m}=\chi _{A_{m}}f_{n}$. Therefore, there exists a
unique $f\in L_{0}^{wo}\left( \nu ;B\left( H\right) \right) $ such that $%
\chi _{A_{n}}f=\chi _{A_{n}}f_{n}$ for all $n$. We claim that $x=x_{f}$.
Indeed, 
\begin{equation*}
x\left( \chi _{A_{n}}\otimes \mathbf{1}\right) =xp_{n}\left( \chi
_{A_{n}}\otimes \mathbf{1}\right) =x_{f_{n}}\left( \chi _{A_{n}}\otimes 
\mathbf{1}\right) =x_{f}\left( \chi _{A_{n}}\otimes \mathbf{1}\right)
\end{equation*}%
for all $n$, from which it easily follows that $x=x_{f}$.

Assume now that $x\in LS\left( \mathcal{M}\right) $. There exists a sequence 
$\left\{ q_{n}\right\} _{n=1}^{\infty }$ of mutually orthogonal central
projections in $\mathcal{M}$ such that $\sum_{n=1}^{\infty }q_{n}=\mathbf{1}$
and $xq_{n}\in S\left( \mathcal{M}\right) $ for all $n$. Each $q_{n}$ is of
the form $q_{n}=\chi _{B_{n}}\otimes \mathbf{1}$, where $B_{n}$ is a $\nu $%
-measurable set. Since $q_{n}q_{m}=0$ whenever $n\neq m$, it follows that
the sets $\left\{ B_{n}\right\} _{n=1}^{\infty }$ are mutually disjoint
(modulo $\nu $-null sets). Furthermore, it is clear that $%
\bigcup\nolimits_{n}B_{n}=S$. From the first part of the proof it follows
that for every $n$ there exists a unique $f_{n}\in L_{0}^{wo}\left( \nu
;B\left( H\right) \right) $ such that $xq_{n}=x_{f_{n}}$. A moment's
reflection shows that $f_{n}=\chi _{B_{n}}f_{n}$ for all $n$ and so, there
exists a unique $f\in L_{0}^{wo}\left( \nu ;B\left( H\right) \right) $ such
that $f_{n}=\chi _{B_{n}}f$ for all $n\in \mathbb{N}$. Since $x\left( \chi
_{B_{n}}\otimes \mathbf{1}\right) =x_{f}\left( \chi _{B_{n}}\otimes \mathbf{1%
}\right) $ for all $n$, it is now clear that $x=x_{f}$, which completes the
proof of the proposition.\medskip
\end{proof}

\begin{remark}
\label{PRem03}In \cite{Ye} Definition 3.1, Yeadon introduces the \emph{%
topology of convergence locally in measure} in the algebra $LS\left( 
\mathcal{M}\right) $ for a general von Neumann algebra $\mathcal{M}$. If $%
\mathcal{M}=L_{\infty }\left( \nu \right) \overline{\otimes }B\left(
H\right) $ and the algebra $LS\left( \mathcal{M}\right) $ is identified with
the algebra $L_{0}^{wo}\left( \nu ;B\left( H\right) \right) $ via the above
proposition, then this topology in $LS\left( \mathcal{M}\right) $
corresponds to the topology of convergence in measure on sets of finite
measure in $L_{0}^{wo}\left( \nu ;B\left( H\right) \right) $ (as defined in
Section \ref{SectVecVal} of the present paper). To verify this claim, it
should be observed that $\mathcal{Z}\left( \mathcal{M}\right) \cong
L_{\infty }\left( \nu \right) $ and that the \emph{dimension function} $%
D:P\left( \mathcal{M}\right) \rightarrow M\left( \nu \right) ^{+}$ (in the
sense of Segal, \cite{Se} Definition 1.4) is given by $\left( Dp\right)
\left( s\right) =\func{tr}\left( p\left( s\right) \right) $, $s\in S$. Here, 
$M\left( \nu \right) ^{+}$ denotes the collection of all (equivalence
classes of) $\nu $-measurable functions with values in $\left[ 0,\infty %
\right] $. The details are left to the reader.
\end{remark}

If $f\in L_{0}^{wo}\left( \nu ;B\left( H\right) \right) $ is self-adjoint,
that is, $f\left( s\right) ^{\ast }=f\left( s\right) $ $\nu $-a.e. on $S$,
and if $\phi $ is a bounded Borel function on $\mathbb{R}$, then $\phi
\left( f\right) :s\longmapsto \phi \left( f\left( s\right) \right) $, $s\in
S $, is a bounded function from $S$ into $B\left( H\right) $. On the other
hand, the operator $x_{f}\in LS\left( \mathcal{M}\right) $ is self-adjoint
and so, the operator $\phi \left( x_{f}\right) \in \mathcal{M}$ is defined
via the functional calculus of $x_{f}$. The next lemma describes the
relationship between $\phi \left( f\right) $ and $\phi \left( x_{f}\right) $%
. The $\ast $-algebra of all bounded Borel functions on $\mathbb{R}$ is
denoted by $B_{b}\left( \mathbb{R}\right) $.

\begin{lemma}
\label{PLem09}If $f\in L_{0}^{wo}\left( \nu ;B\left( H\right) \right) $ is
self-adjoint and $\phi \in B_{b}\left( \mathbb{R}\right) $, then $\phi
\left( f\right) \in L_{\infty }^{wo}\left( \nu ;B\left( H\right) \right) $
and $x_{\phi \left( f\right) }=\phi \left( x_{f}\right) $.
\end{lemma}

\begin{proof}
Let $f\in L_{0}^{wo}\left( \nu ;B\left( H\right) \right) $ be fixed and
self-adjoint. Define $\mathcal{U}\subseteq B_{b}\left( \mathbb{R}\right) $
by setting 
\begin{equation*}
\mathcal{U}=\left\{ \phi \in B_{b}\left( \mathbb{R}\right) :\phi \left(
f\right) \in L_{\infty }^{wo}\left( \nu ;B\left( H\right) \right) \text{ and 
}x_{\phi \left( f\right) }=\phi \left( x_{f}\right) \right\} .
\end{equation*}%
It is easily verified that $\mathcal{U}$ is a unital $\ast $-subalgebra of $%
B_{b}\left( \mathbb{R}\right) $. If $\left\{ \phi _{n}\right\}
_{n=1}^{\infty }$ is a sequence in $\mathcal{U}$ and if $\phi \in
B_{b}\left( \mathbb{R}\right) $ is such that $\phi _{n}\rightarrow \phi $
uniformly on $\mathbb{R}$, then it follows from the properties of the
functional calculus that $\phi _{n}\left( f\right) \rightarrow \phi \left(
f\right) $ uniformly on $S$. Hence, $\phi \left( f\right) \in L_{\infty
}^{wo}\left( \nu ;B\left( H\right) \right) $ and $\left\Vert x_{\phi
_{n}\left( f\right) }-x_{\phi \left( f\right) }\right\Vert _{\mathcal{M}%
}\rightarrow 0$ as $n\rightarrow \infty $. Similarly, $\left\Vert \phi
_{n}\left( x_{f}\right) -\phi \left( x_{f}\right) \right\Vert _{\mathcal{M}%
}\rightarrow 0$ and so, $x_{\phi \left( f\right) }=\phi \left( x_{f}\right) $%
. This shows that $\mathcal{U}$ is uniformly closed in $B_{b}\left( \mathbb{R%
}\right) $. If $\lambda \in \mathbb{C}\diagdown \mathbb{R}$ and the function 
$r_{\lambda }\in B_{b}\left( \mathbb{R}\right) $ is defined by setting $%
r_{\lambda }\left( t\right) =\left( \lambda -t\right) ^{-1}$, $t\in \mathbb{R%
}$, then a straightforward calculation shows that $x_{r_{\lambda }\left(
f\right) }=r_{\lambda }\left( x_{f}\right) =\left( \lambda \mathbf{1}%
-x_{f}\right) ^{-1}$. Consequently, $r_{\lambda }\in \mathcal{U}$ for all $%
\lambda \in \mathbb{C}\diagdown \mathbb{R}$. The Stone-Weierstrass theorem
now implies that $C_{0}\left( \mathbb{R}\right) \subseteq \mathcal{U}$.

Suppose that the sequence $\left\{ \phi _{n}\right\} _{n=1}^{\infty }$ in $%
\mathcal{U}^{+}$ and $\phi \in B_{b}\left( \mathbb{R}\right) ^{+}$ are such
that $0\leq \phi _{n}\left( t\right) \uparrow _{n}\phi \left( t\right) $, $%
t\in \mathbb{R}$. From the properties of the functional calculus it follows
that $\phi _{n}\left( x_{f}\right) \uparrow _{n}\phi \left( x_{f}\right) $
in $\mathcal{M}$. On the other hand, 
\begin{equation*}
\phi _{n}\left( f\left( s\right) \right) \uparrow _{n}\phi \left( f\left(
s\right) \right) ,\ \ \ s\in S\text{,}
\end{equation*}%
in $B\left( H\right) $. Hence, $\phi \left( f\right) \in L_{\infty
}^{wo}\left( S;B\left( H\right) \right) $ and $\phi _{n}\left( f\right)
\uparrow \phi \left( f\right) $ in $L_{\infty }^{wo}\left( S;B\left(
H\right) \right) $, which implies that $x_{\phi _{n}\left( f\right)
}\uparrow _{n}x_{\phi \left( f\right) }$ in $\mathcal{M}$. This shows that $%
\phi \in \mathcal{U}^{+}$. From this observation it follows that $\mathcal{U}%
^{+}$ is closed in $B_{b}\left( \mathbb{R}\right) ^{+}$ under monotone
pointwise convergence of sequences. Together with the inclusion $C_{0}\left( 
\mathbb{R}\right) \subseteq \mathcal{U}$, this suffices to conclude that $%
\mathcal{U}=B_{b}\left( \mathbb{R}\right) $, which completes the proof of
the lemma.\medskip
\end{proof}

In the sequel, we shall frequently identify a function $f\in
L_{0}^{wo}\left( S;B\left( H\right) \right) $ with the corresponding
operator $x_{f}\in LS\left( \mathcal{M}\right) $. An immediate consequence
of Lemma \ref{PLem09} is the following.

\begin{corollary}
\label{PCor01}If $f\in L_{0}^{wo}\left( S;B\left( H\right) \right) $ is
self-adjoint, then the spectral measure $e^{f}$ of $f$ (as a self-adjoint
operator in $LS\left( \mathcal{M}\right) $) is given by 
\begin{equation*}
e^{f}\left( B\right) \left( s\right) =e^{f\left( s\right) }\left( B\right)
,\ \ \ s\in S,
\end{equation*}%
for every Borel set $B\subseteq \mathbb{R}$. In particular, for every Borel
set $B\subseteq \mathbb{R}$, the function $s\longmapsto e^{f\left( s\right)
}\left( B\right) $, $s\in S$, is $\nu $-measurable.
\end{corollary}

To apply the results of Section \ref{SectAdmis} in the present setting, the
following observations are relevant.

\begin{proposition}
\begin{enumerate}
\item[(i).] \label{PProp07}The algebras $LS\left( \mathcal{M}\right) $, $%
S\left( \mathcal{M}\right) $ and $S\left( \tau \right) $ are $B\left(
H\right) $-admissible.

\item[(ii).] The algebra $S_{0}\left( \tau \right) $ is $K\left( H\right) $%
-admissible.
\end{enumerate}
\end{proposition}

\begin{proof}
(i). It is evident that $LS\left( \mathcal{M}\right) \cong L_{0}^{wo}\left(
\nu ;B\left( H\right) \right) $ is $B\left( H\right) $-admissible. It is
clear that $S\left( \mathcal{M}\right) $ and $S\left( \tau \right) $ satisfy
condition (i) of Definition \ref{DDef01}. Since 
\begin{equation*}
L_{\infty }^{wo}\left( \nu ;B\left( H\right) \right) \cong \mathcal{M}%
\subseteq S\left( \tau \right) \subseteq S\left( \mathcal{M}\right) ,
\end{equation*}%
it is also clear that both $S\left( \mathcal{M}\right) $ and $S\left( \tau
\right) $ satisfy condition (ii). If $f\in L_{0}^{wo}\left( \nu ;B\left(
H\right) \right) $, then there exists a $\nu $-measurable partition $\left\{
A_{n}\right\} _{n=1}^{\infty }$ of $S$ such that $\chi _{A_{n}}f\in
L_{\infty }^{wo}\left( \nu ;B\left( H\right) \right) $, that is, $\chi
_{A_{n}}f\in \mathcal{M}$ for all $n$. Hence, $S\left( \mathcal{M}\right) $
and $S\left( \tau \right) $ also satisfy condition (iii) of Definition \ref%
{DDef01}.

(ii). It should be observed first that $S_{0}\left( \tau \right) \subseteq
L_{0}^{wo}\left( \nu ;K\left( H\right) \right) $. Indeed, if $f\in
S_{0}\left( \tau \right) $, then there exists for every $n\in \mathbb{N}$ a
projection $p_{n}\in P\left( \mathcal{M}\right) $ satisfying $\left\Vert
fp_{n}\right\Vert _{\infty }\leq 1/n$ and $\tau \left( p_{n}^{\bot }\right)
<\infty $. Using (\ref{Peq10}), it follows that there exists a $\nu $%
-measurable set $A\subseteq S$ satisfying $\nu \left( S\diagdown A\right) =0$%
, such that $\left\Vert f\left( s\right) p_{n}\left( s\right) \right\Vert
_{B\left( H\right) }\leq 1/n$ and $\func{tr}\left( p_{n}\left( s\right)
^{\bot }\right) <\infty $ for all $s\in A$. Consequently, $f\left( s\right)
\in K\left( H\right) $ for all $s\in A$ and so, $f\in L_{0}^{wo}\left( \nu
;K\left( H\right) \right) $. Since $S_{0}\left( \tau \right) $ is a
two-sided ideal in $S\left( \tau \right) $, it is clear that $S_{0}\left(
\tau \right) $ is an $L_{\infty }\left( \nu \right) $-module. To show that $%
S_{0}\left( \tau \right) $ satisfies condition (ii) of Definition \ref%
{DDef01}, suppose that $f\in L_{\infty }^{wo}\left( \nu ;B\left( H\right)
\right) $ is such that for every $\varepsilon >0$ there exists $p\in 
\mathcal{I}\left( F\left( H\right) \right) $ satisfying $\left\Vert f\left( 
\mathbf{1}\otimes p^{\bot }\right) \right\Vert _{\infty }\leq \varepsilon $
and let the $\nu $-measurable set $A\subseteq S$ with $\nu \left( A\right)
<\infty $ be given. We have to show that $\chi _{A}f\in S_{0}\left( \tau
\right) $. If $\varepsilon >0$, then there exists $p\in \mathcal{I}\left(
F\left( H\right) \right) $ such that $\left\Vert f\left( \mathbf{1}\otimes
p^{\bot }\right) \right\Vert _{\infty }\leq \varepsilon $. Replacing $p$ by
the orthogonal projection in $H$ onto $\func{Ker}\left( p\right) ^{\bot }$,
it is easily seen that the projection $p$ may be assumed to be self-adjoint,
that is, $p\in P\left( F\left( H\right) \right) $. Defining the projection $%
q\in P\left( \mathcal{M}\right) $ by $q=\chi _{A}\otimes p$, it is clear
that $\tau \left( q\right) =\nu \left( A\right) \func{tr}\left( p\right)
<\infty $. Furthermore, $q^{\bot }=\chi _{A}\otimes p^{\bot }+\chi
_{A^{c}}\otimes \mathbf{1}$ and so, 
\begin{equation*}
\left\Vert \left( \chi _{A}f\right) q^{\bot }\right\Vert _{\infty
}=\left\Vert \left( \chi _{A}f\right) \left( \mathbf{1}\otimes p^{\bot
}\right) \right\Vert _{\infty }\leq \varepsilon .
\end{equation*}%
This shows that $\chi _{A}f\in S_{0}\left( \tau \right) $.

Finally, it will be shown that $S_{0}\left( \tau \right) $ satisfies
condition (iii) of Definition \ref{DDef01} with respect to $\mathfrak{U}%
=K\left( H\right) $. Since the measure $\nu $ is assumed to be $\sigma $%
-finite, it suffices to consider the case that $\nu \left( S\right) <\infty $%
. Given $f\in L_{0}^{wo}\left( \nu ;K\left( H\right) \right) $, define the
sets $B_{k,n}$ ($k,n\in \mathbb{N}$) by setting 
\begin{equation*}
B_{k,n}=\left\{ s\in S:\func{tr}\left( e^{\left\vert f\left( s\right)
\right\vert }\left( 1/k,\infty \right) \right) >n\right\} .
\end{equation*}%
It follows from Corollary \ref{PCor01} that these sets are $\nu $%
-measurable. Furthermore, since $f\left( s\right) \in K\left( H\right) $, it
is clear that $\func{tr}\left( e^{\left\vert f\left( s\right) \right\vert
}\left( 1/k,\infty \right) \right) \in \mathbb{N}$ for all $k\in \mathbb{N}$
and $s\in S$. Consequently, $B_{k,n}\downarrow _{n}\emptyset $ for all $k$
and so, for each $k$ there exists $n_{k}\in \mathbb{N}$ such that $\nu
\left( B_{k,n_{k}}\right) \leq 2^{-k}$. Defining 
\begin{equation*}
A_{l}=\bigcap_{k=l}^{\infty }B_{k,n_{k}}^{c},\ \ \ l\in \mathbb{N},
\end{equation*}%
it is clear that $\nu \left( A_{l}^{c}\right) \leq 2^{-l+1}$ and so, $%
A_{l}\uparrow S$. Furthermore, if $s\in A_{l}$, then $\func{tr}\left(
e^{\left\vert f\left( s\right) \right\vert }\left( 1/k,\infty \right)
\right) \leq n_{k}$ for all $k\geq l$. We claim that $\chi _{A_{l}}f\in
S_{0}\left( \tau \right) $ for all $l\in \mathbb{N}$. Indeed, given $l$ and $%
\varepsilon >0$ let $k\in \mathbb{N}$ be such that $k\geq l$ and $1/k\leq
\varepsilon $. The projection $p=\chi _{A_{l}}e^{\left\vert f\right\vert
}\left( \varepsilon ,\infty \right) $ satisfies 
\begin{eqnarray*}
\tau \left( p\right) &=&\int_{A_{l}}\func{tr}\left( e^{\left\vert f\left(
s\right) \right\vert }\left( \varepsilon ,\infty \right) \right) d\nu \left(
s\right) \\
&\leq &\int_{A_{l}}\func{tr}\left( e^{\left\vert f\left( s\right)
\right\vert }\left( 1/k,\infty \right) \right) d\nu \left( s\right) \leq
n_{k}\nu \left( A_{l}\right) <\infty
\end{eqnarray*}%
and 
\begin{equation*}
\left\Vert \left( \chi _{A_{l}}f\right) p^{\bot }\right\Vert _{\infty
}=\left\Vert \left( \chi _{A_{l}}f\right) e^{\left\vert f\right\vert }\left[
0,\varepsilon \right] \right\Vert _{\infty }\leq \varepsilon .
\end{equation*}%
Hence, $\chi _{A_{l}}f\in S_{0}\left( \tau \right) $ for all $l$. This
suffices to show that $S_{0}\left( \tau \right) $ satisfies condition (iii)
of Definition \ref{DDef01} with respect to $\mathfrak{U}=K\left( H\right) $%
.\medskip
\end{proof}

A combination of Theorem \ref{PThm01}, Corollary \ref{PThm011} and
Proposition \ref{PProp07} immediately yield the following results.

\begin{corollary}
\label{PCor02}Let $H$ be a separable Hilbert space, $\left( S,\Sigma ,\nu
\right) $ be a $\sigma $-finite measure space and $\mathcal{M}=L_{\infty
}\left( \nu \right) \overline{\otimes }B\left( H\right) $. If $\mathfrak{A}$
is one of the algebras $LS\left( \mathcal{M}\right) $, $S\left( \mathcal{M}%
\right) $, $S\left( \tau \right) $ or $S_{0}\left( \tau \right) $ and if $%
\delta $ is a $\mathcal{Z}\left( \mathcal{M}\right) $-linear derivation in $%
\mathfrak{A}$, then there exists $w\in LS\left( \mathcal{M}\right) $ such
that $\delta \left( x\right) =wx-xw$, $x\in \mathfrak{A}$.
\end{corollary}

\begin{corollary}
\label{PCor021}Let $H$ be a separable infinite dimensional Hilbert space,
\linebreak $\left( S,\Sigma ,\nu \right) $ be a $\sigma $-finite measure
space and $\mathcal{M}=L_{\infty }\left( \nu \right) \overline{\otimes }%
B\left( H\right) $. If $\mathfrak{A}$ is one of the algebras $LS\left( 
\mathcal{M}\right) $, $S\left( \mathcal{M}\right) $, $S\left( \tau \right) $
or $S_{0}\left( \tau \right) $ and if $\delta $ is a derivation in $%
\mathfrak{A}$, then there exists $w\in LS\left( \mathcal{M}\right) $ such
that $\delta \left( x\right) =wx-xw$, $x\in \mathfrak{A}$.
\end{corollary}

If $\mathfrak{A}$ is equal to either $S\left( \mathcal{M}\right) $ or $%
S\left( \tau \right) $, then it actually follows that the element $w$ in
Corollary \ref{PCor02} may be chosen from the algebra $\mathfrak{A}$ itself
(that is, $\delta $ is inner). This will follow from the general result of
Proposition \ref{AProp02}, which is of interest in its own right. The proof
of this proposition is divided into a number of lemmas.

Let $\mathcal{M}$ be a von Neumann algebra on a Hilbert space $H$ and
suppose that $\mathcal{A}$ is a $\ast $-subalgebra of $LS\left( \mathcal{M}%
\right) $. If $w\in LS\left( \mathcal{M}\right) $ is such that $\left[ w,x%
\right] \in \mathcal{A}$ for all $x\in \mathcal{A}$, then the map $\delta
_{w}:\mathcal{A}\rightarrow \mathcal{A}$, defined by setting $\delta
_{w}\left( x\right) =\left[ w,x\right] $, $x\in \mathcal{A}$, is a
derivation in $\mathcal{A}$. In general, $\delta _{w}$ need not be inner in $%
\mathcal{A}$. It will be shown that if $\mathcal{A}$ is \textit{absolutely
solid} in $LS\left( \mathcal{M}\right) $ (that is, if $x\in LS\left( 
\mathcal{M}\right) $ and $y\in \mathcal{A}$ are such that $\left\vert
x\right\vert \leq \left\vert y\right\vert $, then $x\in \mathcal{A}$) and if 
$\mathcal{M}\subseteq \mathcal{A}$, then such a derivation $\delta _{w}$ is
always inner in $\mathcal{A}$, that is, there exists $w_{1}\in \mathcal{A}$
such that $\delta _{w}=\delta _{w_{1}}$.

The following theorem will be used (see \cite{KR}, Theorem 6.2.7).

\begin{theorem}
If $e,f\in P\left( \mathcal{M}\right) $, then there exist unique mutually
orthogonal central projection $p_{0}$ and $q_{0}$ which are maximal with
respect to the properties:

\begin{enumerate}
\item $q_{0}e\sim q_{0}f$;

\item $pe\prec pf$ whenever $p\in P\left( \mathcal{Z}\left( \mathcal{M}%
\right) \right) $ satisfies $p\leq p_{0}$;

\item $pf\prec pe$ whenever $p\in P\left( \mathcal{Z}\left( \mathcal{M}%
\right) \right) $ satisfies $p\leq \mathbf{1}-p_{0}-q_{0}$.
\end{enumerate}
\end{theorem}

For $e,f\in P\left( \mathcal{M}\right) $ define $z\left( e,f\right) \in
P\left( \mathcal{Z}\left( \mathcal{M}\right) \right) $ by setting $z\left(
e,f\right) =p_{0}+q_{0}$. It should be observed that $z\left( e,f\right) $
is also given by 
\begin{equation}
z\left( e,f\right) =\bigvee \left\{ p\in P\left( \mathcal{Z}\left( \mathcal{M%
}\right) \right) :pe\precsim pf\right\} .  \label{Aeq03}
\end{equation}%
Moreover, $z\left( f,e\right) =\mathbf{1}-p_{0}$. This implies, in
particular, that 
\begin{equation*}
z\left( f,e\right) \geq \mathbf{1}-z\left( e,f\right)
\end{equation*}%
and so, 
\begin{equation}
z\left( e,f\right) \vee z\left( f,e\right) =\mathbf{1}.  \label{Aeq02}
\end{equation}

Let $a\in LS_{h}\left( \mathcal{M}\right) $ with spectral measure $e^{a}$.
For $n\in \mathbb{Z}$, define the central projections $p_{n}^{a}$ and $%
q_{n}^{a}$ by setting 
\begin{equation*}
p_{n}^{a}=z\left( e^{a}\left( -\infty ,n\right] ,e^{a}\left( n+1,\infty
\right) \right) ,\ \ \ q_{n}^{a}=z\left( e^{a}\left( n+1,\infty \right)
,e^{a}\left( -\infty ,n\right] \right) .
\end{equation*}%
In the proof of Lemma \ref{ALem07} the following observation will be used.

\begin{lemma}
\label{ALem02}If $\left\{ e_{n}\right\} _{n=1}^{\infty }$ is a sequence of
finite projections in $\mathcal{M}$ such that $e_{n}\downarrow 0$ and if $%
p\in P\left( \mathcal{M}\right) $ is such that $p\precsim e_{n}$ for all $n$%
, then $p=0$.
\end{lemma}

\begin{proof}
Since $e_{1}$ is finite and $p\precsim e_{1}$, it follows that $p$ is finite
and so, the projection $e=e_{1}\vee p$ is also finite. The reduced von
Neumann algebra $e\mathcal{M}e$ is finite, $e_{n}\downarrow 0$ in $e\mathcal{%
M}e$ and $p\precsim e_{n}$ in $e\mathcal{M}e$ for all $n$. Let $T:e\mathcal{M%
}e\rightarrow Z\left( e\mathcal{M}e\right) $ be the center-valued trace (see 
\cite{KR}, Theorem 8.2.8). It follows from $e_{n}\downarrow 0$ that $T\left(
e_{n}\right) \downarrow 0$, and $p\precsim e_{n}$ implies that $T\left(
p\right) \leq T\left( e_{n}\right) $ for all $n$. Hence, $T\left( p\right)
=0 $ and so $p=0$, which proves the lemma.\medskip
\end{proof}

\begin{lemma}
\label{ALem07}The central projections $\left\{ p_{n}^{a}\right\} _{n\in 
\mathbb{Z}}$ and $\left\{ q_{n}^{a}\right\} _{n\in \mathbb{Z}}$ satisfy:

\begin{enumerate}
\item[(i).] $p_{n}^{a}\vee q_{n}^{a}=\mathbf{1}$ for all $n\in \mathbb{Z}$;

\item[(ii).] $p_{n}^{a}\downarrow _{n\in \mathbb{Z}}$ and $q_{n}^{a}\uparrow
_{n\in \mathbb{Z}}$;

\item[(iii).] $\bigwedge\nolimits_{n\in \mathbb{Z}}p_{n}^{a}=0$ and $%
\bigwedge\nolimits_{n\in \mathbb{Z}}q_{n}^{a}=0$.
\end{enumerate}
\end{lemma}

\begin{proof}
(i). This follows immediately from (\ref{Aeq02}).

(ii). If $n\in \mathbb{Z}$, then 
\begin{eqnarray*}
p_{n+1}^{a}e^{a}\left( -\infty ,n\right] &\leq &p_{n+1}^{a}e^{a}\left(
-\infty ,n+1\right] \\
&\precsim &p_{n+1}^{a}e^{a}\left( n+2,\infty \right) \leq
p_{n+1}^{a}e^{a}\left( n+1,\infty \right)
\end{eqnarray*}%
and so, it follows from (\ref{Aeq03}) that $p_{n+1}^{a}\leq p_{n}^{a}$.
Similarly, 
\begin{eqnarray*}
q_{n}^{a}e^{a}\left( n+2,\infty \right) &\leq &q_{n}^{a}e^{a}\left(
n+1,\infty \right) \\
&\precsim &q_{n}^{a}e^{a}\left( -\infty ,n\right] \leq q_{n}^{a}e^{a}\left(
-\infty ,n+1\right]
\end{eqnarray*}%
and hence, $q_{n}^{a}\leq q_{n+1}^{a}$.

(iii).Since $a\in LS\left( \mathcal{M}\right) $, there exists a sequence $%
\left\{ z_{k}\right\} _{k=1}^{\infty }$ in $P\left( \mathcal{Z}\left( 
\mathcal{M}\right) \right) $ such that $z_{k}\uparrow \mathbf{1}$ and $%
az_{k}\in S\left( \mathcal{M}\right) $ for all $k\in \mathbb{N}$. Defining $%
p_{\infty }\in P\left( \mathcal{Z}\left( \mathcal{M}\right) \right) $ by $%
p_{\infty }=\bigwedge\nolimits_{n\in \mathbb{Z}}p_{n}^{a}$, it follows that 
\begin{equation}
p_{\infty }e^{a}\left( -\infty ,n\right] \precsim p_{\infty }e^{a}\left(
n+1,\infty \right) ,\ \ \ n\in \mathbb{Z}.  \label{Aeq04}
\end{equation}%
(\cite{KR}, Proposition 6.2.3). It follows from (\ref{Aeq04}) that for each $%
k\in \mathbb{N}$,%
\begin{equation*}
p_{\infty }z_{k}e^{az_{k}}\left( -\infty ,n\right] \precsim p_{\infty
}z_{k}e^{az_{k}}\left( n+1,\infty \right) ,\ \ \ n\in \mathbb{Z}.
\end{equation*}%
Let $k\in \mathbb{N}$ be fixed for the moment. Since $az_{k}\in S\left( 
\mathcal{M}\right) $, there exists $n_{0}\in \mathbb{N}$ such that $%
e^{az_{k}}\left( n_{0}+1,\infty \right) $ is a finite projection.
Consequently, if $m\in \mathbb{Z}$ and $n\in \mathbb{N}$ satisfies $n\geq
\max \left( m,n_{0}\right) $, then 
\begin{equation*}
p_{\infty }z_{k}e^{az_{k}}\left( -\infty ,m\right] \leq p_{\infty
}z_{k}e^{az_{k}}\left( -\infty ,n\right] \precsim p_{\infty
}z_{k}e^{az_{k}}\left( n+1,\infty \right)
\end{equation*}%
and so, it follows from Lemma \ref{ALem02} that $p_{\infty
}z_{k}e^{az_{k}}\left( -\infty ,m\right] =0$. Letting $m\uparrow \infty $,
it follows that $p_{\infty }z_{k}=0$ for all $k\in \mathbb{N}$. Since $%
z_{k}\uparrow \mathbf{1}$, this implies that $p_{\infty }=0$. This shows
that $p_{n}^{a}\downarrow _{n\in \mathbb{Z}}0$.

Similarly, if $q_{-\infty }=\bigwedge\nolimits_{n\in \mathbb{Z}}q_{n}^{a}$,
then 
\begin{equation*}
q_{-\infty }e^{a}\left( n+1,\infty \right) \precsim q_{-\infty }e^{a}\left(
-\infty ,n\right] ,\ \ \ n\in \mathbb{Z},
\end{equation*}%
which implies that 
\begin{equation*}
q_{-\infty }z_{k}e^{az_{k}}\left( n+1,\infty \right) \precsim q_{-\infty
}z_{k}e^{az_{k}}\left( -\infty ,n\right] ,\ \ \ n\in \mathbb{Z},
\end{equation*}%
for all $k\in \mathbb{N}$. If $n<0$, then $e^{az_{k}}\left( -\infty ,n\right]
\leq e^{\left\vert az_{k}\right\vert }\left[ -n,\infty \right) $. Since $%
az_{k}\in S\left( \mathcal{M}\right) $, there exists $n_{0}\in \mathbb{Z}$
such that $e^{\left\vert az_{k}\right\vert }\left[ -n,\infty \right) $ is a
finite projection for all $n\leq n_{0}$. As above, it now follows that $%
q_{-\infty }z_{k}=0$ for all $k\in \mathbb{N}$. Since $z_{k}\uparrow \mathbf{%
1}$, this implies that $q_{-\infty }=0$. \medskip
\end{proof}

The following general observation concerning Boolean algebras will be needed.

\begin{lemma}
\label{ALem08}Let $\mathcal{B}$ be a complete Boolean algebra (with smallest
element $\mathbf{0}$ and largest element $\mathbf{1}$).

\begin{enumerate}
\item[(i).] If $\left\{ A_{n}\right\} _{n\in \mathbb{Z}}$ and $\left\{
B_{n}\right\} _{n\in \mathbb{Z}}$ are two sequences in $\mathcal{B}$ such
that $A_{n}\downarrow _{n\in \mathbb{Z}}$, $B_{n}\uparrow _{n\in \mathbb{Z}}$%
, $\bigwedge\nolimits_{n\in \mathbb{Z}}A_{n}=0$, $\bigwedge\nolimits_{n\in 
\mathbb{Z}}B_{n}=0$ and $A_{n}\vee B_{n}=\mathbf{1}$ for all $n$, then $%
\bigvee\nolimits_{n\in \mathbb{Z}}\left( A_{n}\wedge B_{n+1}\right) =\mathbf{%
1}$.

\item[(ii).] If $\left\{ A_{j}\right\} _{j\in J}$ is a countable family in $%
\mathcal{B}$ such that $\bigvee\nolimits_{j\in J}A_{j}=\mathbf{1}$, then
there exists a disjoint family $\left\{ B_{j}\right\} _{j\in J}$ in $%
\mathcal{B}$ such that $B_{j}\leq A_{j}$ for all $j\in J$ and $%
\bigvee\nolimits_{j\in J}B_{j}=\mathbf{1}$.
\end{enumerate}
\end{lemma}

\begin{proof}
(i). It follows from $A_{n}\vee B_{n}=\mathbf{1}$ that $A_{n}\geq B_{n}^{c}$
and $B_{n}\geq A_{n}^{c}$ for all $n$. This implies, in particular, that $%
B_{n}\uparrow _{n\in \mathbb{Z}}\mathbf{1}$. Hence, for each $n\in \mathbb{Z}
$, 
\begin{equation*}
B_{n}\vee \bigvee\nolimits_{k=n}^{\infty }\left( B_{k+1}\diagdown
B_{k}\right) =\mathbf{1}
\end{equation*}%
and so, 
\begin{equation*}
\bigvee\nolimits_{k=n}^{\infty }\left( B_{k+1}\diagdown B_{k}\right)
=B_{n}^{c}.
\end{equation*}%
Consequently, 
\begin{equation*}
\bigvee\nolimits_{n\in \mathbb{Z}}\left( B_{n+1}\diagdown B_{n}\right)
=\bigvee\nolimits_{n\in \mathbb{Z}}B_{n}^{c}=\mathbf{1}.
\end{equation*}%
It follows from $A_{n}\geq B_{n}^{c}$ that 
\begin{equation*}
A_{n}\wedge B_{n+1}\geq B_{n}^{c}\wedge B_{n+1}=B_{n+1}\diagdown B_{n}
\end{equation*}%
and hence, 
\begin{equation*}
\bigvee\nolimits_{n\in \mathbb{Z}}\left( A_{n}\wedge B_{n+1}\right) =\mathbf{%
1}.
\end{equation*}

(ii). Without loss of generality, it may be assumed that $J=\mathbb{N}$.
Defining $B_{1}=A_{1}$ and 
\begin{equation*}
B_{j}=A_{j}\diagdown \left( A_{1}\vee \cdots \vee A_{j-1}\right) ,\ \ \
j\geq 2,
\end{equation*}%
it is clear that the sequence $\left\{ B_{j}\right\} _{j\in \mathbb{N}}$ has
the desired properties.\medskip
\end{proof}

A combination of Lemma \ref{ALem07} and Lemma \ref{ALem08} (applied to the
complete Boolean algebra of central projections in $\mathcal{M}$) yields the
following result.

\begin{corollary}
\label{ACor01}If $a\in LS_{h}\left( \mathcal{M}\right) $, then there exists
a disjoint sequence \linebreak $\left\{ z_{n}\right\} _{n\in \mathbb{Z}}$ in 
$P\left( \mathcal{Z}\left( \mathcal{M}\right) \right) $ such that $z_{n}\leq
p_{n}^{a}q_{n+1}^{a}$ for all $n\in \mathbb{Z}$ and $\bigvee\nolimits_{n\in 
\mathbb{Z}}z_{n}=\mathbf{1}$.
\end{corollary}

\begin{lemma}
\label{ALem09}If $a\in LS_{h}\left( \mathcal{M}\right) $ and $n\in \mathbb{Z}
$ satisfies $e^{a}\left( -\infty ,n\right] \precsim e^{a}\left( n+1,\infty
\right) $ and $e^{a}\left( n+2,\infty \right) \precsim e^{a}\left( -\infty
,n+1\right] $, then there exists a partial isometry $v\in \mathcal{M}$ such
that 
\begin{equation*}
\left\vert av-va\right\vert +\mathbf{1}\geq \left\vert a-\left( n+1\right) 
\mathbf{1}\right\vert .
\end{equation*}
\end{lemma}

\begin{proof}
Let $u_{1},u_{2}\in \mathcal{M}$ be partial isometries such that 
\begin{eqnarray*}
u_{1}^{\ast }u_{1} &=&e^{a}\left( -\infty ,n\right] \ \ \ \text{and}\ \ \
u_{1}u_{1}^{\ast }\leq e^{a}\left( n+1,\infty \right) , \\
u_{2}^{\ast }u_{2} &=&e^{a}\left( n+2,\infty \right) \ \ \ \text{and}\ \ \
u_{2}u_{2}^{\ast }\leq e^{a}\left( -\infty ,n+1\right] .
\end{eqnarray*}%
It should be observed that 
\begin{eqnarray*}
u_{1}^{\ast }u_{2} &=&\left( e^{a}\left( n+1,\infty \right) u_{1}\right)
^{\ast }\left( e^{a}\left( -\infty ,n+1\right] u_{2}\right) \\
&=&u_{1}^{\ast }e^{a}\left( n+1,\infty \right) e^{a}\left( -\infty ,n+1 
\right] u_{2}=0
\end{eqnarray*}%
and hence, $u_{2}^{\ast }u_{1}=\left( u_{1}^{\ast }u_{2}\right) ^{\ast }=0$.
Defining $v=u_{1}+u_{2}$, this implies that 
\begin{equation*}
v^{\ast }v=u_{1}^{\ast }u_{1}+u_{2}^{\ast }u_{2}=e^{a}\left( -\infty ,n 
\right] +e^{a}\left( n+2,\infty \right)
\end{equation*}%
and so, $v$ is a partial isometry.

Since 
\begin{eqnarray*}
\left\vert v^{\ast }av-v^{\ast }va\right\vert ^{2} &=&\left\vert v^{\ast
}\left( av-va\right) \right\vert ^{2}=\left( av-va\right) ^{\ast }vv^{\ast
}\left( av-va\right) \\
&\leq &\left( av-va\right) ^{\ast }\left( av-va\right) =\left\vert
av-va\right\vert ^{2},
\end{eqnarray*}%
it follows that $\left\vert v^{\ast }av-v^{\ast }va\right\vert \leq
\left\vert av-va\right\vert $. Therefore, it is sufficient to show that 
\begin{equation}
\left\vert v^{\ast }av-v^{\ast }va\right\vert +\mathbf{1}\geq \left\vert
a-\left( n+1\right) \mathbf{1}\right\vert .  \label{Aeq07}
\end{equation}%
It should be observed that $v^{\ast }v$ is a projection which commutes with $%
a$ and that the operator $v^{\ast }av-v^{\ast }va$ is self-adjoint.

Note that 
\begin{eqnarray*}
u_{1}^{\ast }au_{2} &=&\left( e^{a}\left( n+1,\infty \right) u_{1}\right)
^{\ast }a\left( e^{a}\left( -\infty ,n+1\right] u_{2}\right) \\
&=&u_{1}^{\ast }e^{a}\left( n+1,\infty \right) ae^{a}\left( -\infty ,n+1 
\right] u_{2}=0
\end{eqnarray*}%
and hence, $u_{2}^{\ast }au_{1}=\left( u_{1}^{\ast }au_{2}\right) ^{\ast }=0$%
. This implies that 
\begin{equation*}
v^{\ast }av=u_{1}^{\ast }au_{1}+u_{2}^{\ast }au_{2}
\end{equation*}%
and so, 
\begin{equation}
v^{\ast }av-v^{\ast }va=\left[ u_{1}^{\ast }au_{1}-a\right] e^{a}\left(
-\infty ,n\right] +\left[ u_{2}^{\ast }au_{2}-a\right] e^{a}\left(
n+2,\infty \right) .  \label{Aeq05}
\end{equation}%
Since 
\begin{eqnarray*}
u_{1}^{\ast }au_{1} &=&u_{1}^{\ast }ae^{a}\left( n+1,\infty \right)
u_{1}\geq \left( n+1\right) u_{1}^{\ast }e^{a}\left( n+1,\infty \right) u_{1}
\\
&=&\left( n+1\right) u_{1}^{\ast }u_{1}=\left( n+1\right) e^{a}\left(
-\infty ,n\right]
\end{eqnarray*}%
and 
\begin{equation*}
ae^{a}\left( -\infty ,n\right] \leq ne^{a}\left( -\infty ,n\right] \leq
\left( n+1\right) e^{a}\left( -\infty ,n\right] ,
\end{equation*}%
it follows that 
\begin{equation*}
\left[ u_{1}^{\ast }au_{1}-a\right] e^{a}\left( -\infty ,n\right] \geq -%
\left[ a-\left( n+1\right) \mathbf{1}\right] e^{a}\left( -\infty ,n\right]
\geq 0.
\end{equation*}%
Similarly, 
\begin{equation*}
u_{2}^{\ast }au_{2}=u_{2}^{\ast }ae^{a}\left( -\infty ,n+1\right] u_{2}\leq
\left( n+1\right) u_{2}^{\ast }u_{2}=\left( n+1\right) e^{a}\left(
n+2,\infty \right)
\end{equation*}%
and 
\begin{equation*}
ae^{a}\left( n+2,\infty \right) \geq \left( n+2\right) e^{a}\left(
n+2,\infty \right) \geq \left( n+1\right) e^{a}\left( n+2,\infty \right)
\end{equation*}%
and so, 
\begin{equation*}
-\left[ u_{2}^{\ast }au_{2}-a\right] e^{a}\left( n+2,\infty \right) \geq %
\left[ a-\left( n+1\right) \mathbf{1}\right] e^{a}\left( n+2,\infty \right)
\geq 0.
\end{equation*}%
The above estimates, in combination with (\ref{Aeq05}), show that 
\begin{equation}
\left\vert v^{\ast }av-v^{\ast }va\right\vert \geq \left\vert a-\left(
n+1\right) \mathbf{1}\right\vert \left( e^{a}\left( -\infty ,n\right]
+e^{a}\left( n+2,\infty \right) \right) .  \label{Aeq06}
\end{equation}

Observe next that 
\begin{equation*}
ne^{a}\left( n,n+2\right] \leq ae^{a}\left( n,n+2\right] \leq \left(
n+2\right) e^{a}\left( n,n+2\right]
\end{equation*}%
and so, 
\begin{equation*}
-e^{a}\left( n,n+2\right] \leq \left[ a-\left( n+1\right) \mathbf{1}\right]
e^{a}\left( n,n+2\right] \leq e^{a}\left( n,n+2\right] .
\end{equation*}%
This implies that 
\begin{equation}
e^{a}\left( n,n+2\right] \geq \left\vert a-\left( n+1\right) \mathbf{1}%
\right\vert e^{a}\left( n,n+2\right] .  \label{Peq09}
\end{equation}%
Adding (\ref{Peq09}) and (\ref{Aeq06}) yields 
\begin{equation*}
\left\vert v^{\ast }av-v^{\ast }va\right\vert +e^{a}\left( n,n+2\right] \geq
\left\vert a-\left( n+1\right) \mathbf{1}\right\vert ,
\end{equation*}%
from which (\ref{Aeq07}) follows. The proof is complete.\medskip
\end{proof}

\begin{lemma}
\label{ALem10}If $a\in LS_{h}\left( \mathcal{M}\right) $, then there exist a
partial isometry $v\in \mathcal{M}$ and an operator $a_{1}\in LS_{h}\left( 
\mathcal{M}\right) $ such that $a-a_{1}\in \mathcal{Z}\left( LS\left( 
\mathcal{M}\right) \right) _{h}$ and 
\begin{equation*}
\left\vert av-va\right\vert +\mathbf{1}\geq \left\vert a_{1}\right\vert .
\end{equation*}
\end{lemma}

\begin{proof}
Let the sequence $\left\{ z_{n}\right\} _{n\in \mathbb{Z}}$ in $P\left( 
\mathcal{Z}\left( \mathcal{M}\right) \right) $ be as in Corollary \ref%
{ACor01}. For each $n$, the operator $az_{n}$ satisfies the conditions of
Lemma \ref{ALem09} (with respect to the reduced von Neumann algebra $%
\mathcal{M}z_{n}$). Consequently, for each $n\in \mathbb{Z}$ there exists a
partial isometry $v_{n}\in \mathcal{M}z_{n}$ such that 
\begin{equation*}
\left\vert \left( az_{n}\right) v_{n}-v_{n}\left( az_{n}\right) \right\vert
+z_{n}\geq \left\vert az_{n}-\left( n+1\right) z_{n}\right\vert .
\end{equation*}%
Defining the partial isometry $v\in \mathcal{M}$ by setting $v=\sum_{n\in 
\mathbb{Z}}v_{n}z_{n}$, the operator $b\in \mathcal{Z}\left( LS\left( 
\mathcal{M}\right) \right) _{h}$ by $b=\sum_{n\in \mathbb{Z}}\left(
n+1\right) z_{n}$ and $a_{1}=a-b$, it follows that $a-a_{1}=b\in \mathcal{Z}%
\left( LS\left( \mathcal{M}\right) \right) _{h}$ and 
\begin{equation*}
\left\vert av-va\right\vert +\mathbf{1}\geq \left\vert a_{1}\right\vert .
\end{equation*}%
The proof of the lemma is complete.\medskip
\end{proof}

We are now in a position to prove the already advertised proposition.

\begin{proposition}
\label{AProp02}Let $\mathcal{A}$ be a $\ast $-subalgebra of $LS\left( 
\mathcal{M}\right) $ such that $\mathcal{M}\subseteq \mathcal{A}$ and $%
\mathcal{A}$ is absolutely solid (that is, if $x\in LS\left( \mathcal{M}%
\right) $ and $y\in \mathcal{A}$ satisfy $\left\vert x\right\vert \leq
\left\vert y\right\vert $, then $x\in \mathcal{A}$). If $w\in LS\left( 
\mathcal{M}\right) $ is such that $\left[ w,x\right] \in \mathcal{A}$ for
all $x\in \mathcal{A}$, then there exists $w_{1}\in \mathcal{A}$ such that $%
\left[ w,x\right] =\left[ w_{1},x\right] $ for all $x\in \mathcal{A}$.
\end{proposition}

\begin{proof}
First observe that $\left[ w^{\ast },x\right] =-\left[ w,x^{\ast }\right]
^{\ast }\in \mathcal{A}$ for all $x\in \mathcal{A}$. Consequently, $\func{Re}%
w$ and $\func{Im}w$ satisfy the hypothesis of the proposition. Therefore,
without loss of generality, it may be assumed that $w$ is self-adjoint. It
follows from Lemma \ref{ALem10} that there exists a partial isometry $v\in 
\mathcal{M}$ and an operator $w_{1}\in LS\left( \mathcal{M}\right) $ such
that $w-w_{1}\in \mathcal{Z}\left( LS_{h}\left( \mathcal{M}\right) \right) $
and $\left\vert \left[ w,v\right] \right\vert +\mathbf{1}\geq \left\vert
w_{1}\right\vert $. Since $\left\vert \left[ w,v\right] \right\vert +\mathbf{%
1}\in \mathcal{A}$, the assumptions on $\mathcal{A}$ imply that $w_{1}\in 
\mathcal{A}$. Finally, since $w-w_{1}\in \mathcal{Z}\left( LS_{h}\left( 
\mathcal{M}\right) \right) $, it is also clear that $\left[ w,x\right] =%
\left[ w_{1},x\right] $ for all $x\in \mathcal{A}$.\medskip
\end{proof}

Combining Proposition \ref{AProp02} with Corollary \ref{PCor02} and
Corollary \ref{PCor021} we obtain the following results.

\begin{corollary}
\label{PCor03}Let $H$ be a separable Hilbert space, $\left( S,\Sigma ,\nu
\right) $ be a $\sigma $-finite measure space and $\mathcal{M}=L_{\infty
}\left( \nu \right) \overline{\otimes }B\left( H\right) $. Any $\mathcal{Z}%
\left( \mathcal{M}\right)$-linear derivation on each of the algebras $%
LS\left( \mathcal{M}\right) $, $S\left( \mathcal{M}\right) $ or $S\left(
\tau \right) $ is inner.
\end{corollary}

\begin{proof}
For the algebra $LS\left( \mathcal{M}\right) $, this is already stated in
Corollary \ref{PCor02}. Since the algebras $S\left( \mathcal{M}\right) $ and 
$S\left( \tau \right) $ are solid in $LS\left( \mathcal{M}\right) $ and $%
\mathcal{M}$ is contained in both algebras, the result follows from
Corollary \ref{PCor02} in combination with Proposition \ref{AProp02}.\medskip
\end{proof}

\begin{corollary}
\label{PCor031}Let $H$ be a separable infinite dimensional Hilbert space, $%
\left( S,\Sigma ,\nu \right) $ be a $\sigma $-finite measure space and $%
\mathcal{M}=L_{\infty }\left( \nu \right) \overline{\otimes }B\left(
H\right) $. Any derivation on each of the algebras $LS\left( \mathcal{M}%
\right) $, $S\left( \mathcal{M}\right) $ or $S\left( \tau \right) $ is inner.
\end{corollary}

It should be noted that a derivation on $S_{0}\left( \tau \right) $ need not
be inner. Indeed, assuming that $H$ is infinite dimensional, let $p\in
P\left( B\left( H\right) \right) $ be such that $\func{tr}\left( p\right) =%
\func{tr}\left( \mathbf{1}-p\right) =\infty $ and define $w\in P\left( 
\mathcal{M}\right) $ by $w=\chi _{S}\otimes p$. It is easily verified that
the derivation $\delta $ in $S_{0}\left( \tau \right) $, given by $\delta
\left( f\right) =wf-fw$, $f\in S_{0}\left( \tau \right) $, is not inner.

\begin{remark}
\label{PRem02}It should be observed that if $\mathcal{M}=L_{\infty }\left(
\nu \right) \overline{\otimes }B\left( H\right) $, where $S=\left[ 0,1\right]
$ equipped with Lebesgue measure $\nu $ and $H$ is an infinite dimensional
separable Hilbert space, then 
\begin{equation*}
\mathcal{M}\varsubsetneq S\left( \tau \right) \varsubsetneq S\left( \mathcal{%
M}\right) \varsubsetneq LS\left( \mathcal{M}\right) .
\end{equation*}%
To verify this claim, we make the following observations.

\begin{enumerate}
\item[(a).] Let $p_{0}\in B\left( H\right) $ be a rank one projection and $%
a\in L_{0}\left( \nu \right) $. Define $f\in L_{0}^{wo}\left( \nu ;B\left(
H\right) \right) =LS\left( \mathcal{M}\right) $ by setting $f=a\otimes p_{0}$%
. Since $\tau \left( \mathbf{1}\otimes p_{0}\right) =1$ and $f\left( \mathbf{%
1}\otimes p_{0}\right) ^{\bot }=0$, it is clear that $f\in S_{0}\left( \tau
\right) \subseteq S\left( \tau \right) $. If $a\notin L_{\infty }\left( \nu
\right) $, then $f\notin \mathcal{M}$. Hence, $\mathcal{M}\varsubsetneq
S\left( \tau \right) $.

\item[(b).] Let $f\in L_{0}^{wo}\left( \nu ;B\left( H\right) \right) $ be
such that $f\left( s\right) \in K\left( H\right) $ $\nu $-a.e. on $S$ and
let $e^{\left\vert f\right\vert }$ be the spectral measure of $\left\vert
f\right\vert $ (see Corollary \ref{PCor01}). Defining $p_{n}=e^{\left\vert
f\right\vert }\left( n,\infty \right) $, $n\in \mathbb{N}$, it follows that $%
p_{n}\left( s\right) =e^{\left\vert f\left( s\right) \right\vert }\left(
n,\infty \right) $ is a finite rank projection for $\nu $-almost all $s\in S$
and hence, by Lemma \ref{PLem08}, $p_{n}$ is a finite projection in $%
\mathcal{M}$. Since $fp_{n}^{\bot }\in \mathcal{M}$ for all $n$ and $%
p_{n}^{\bot }\uparrow \mathbf{1}$ in $\mathcal{M}$, it follows that $f\in
S\left( \mathcal{M}\right) $.

For $n\in \mathbb{N}$, define the intervals $I_{n}$ by $I_{n}=\left(
2^{-n},2^{-n+1}\right] $ and let $q_{n}\in P\left( B\left( H\right) \right) $
be such that $\limfunc{tr}\left( q_{n}\right) =2^{n}$. Let $f\in
L_{0}^{wo}\left( \nu ;B\left( H\right) \right) ^{+}$ be defined by setting 
\begin{equation*}
f=\sum_{n=1}^{\infty }n\chi _{I_{n}}\otimes q_{n},
\end{equation*}%
as a pointwise convergent series. If $0\leq \lambda \in \mathbb{R}$ and $%
s\in I_{n}$, then $e^{f\left( s\right) }\left( \lambda ,\infty \right)
=e^{nq_{n}}\left( \lambda ,\infty \right) $, which is equal to $q_{n}$ if $%
0\leq \lambda <n$ and equal to $0$ if $\lambda \geq n$. Consequently, 
\begin{equation*}
\tau \left( e^{f}\left( \lambda ,\infty \right) \right) =\sum_{n>\lambda }%
\limfunc{tr}\left( q_{n}\right) \nu \left( I_{n}\right) =\infty
\end{equation*}%
for all $\lambda \geq 0$. This shows that $f\notin S\left( \tau \right) $
and hence, $S\left( \tau \right) \varsubsetneq S\left( \mathcal{M}\right) $.

It should be observed that the function $g\in L_{\infty }^{wo}\left( \nu
;B\left( H\right) \right) ^{+}$, given by $g=\sum_{n=1}^{\infty }\chi
_{I_{n}}\otimes q_{n}$, takes its values in $F\left( H\right) $ but, $\tau
\left( e^{g}\left( \lambda ,\infty \right) \right) =\infty $ whenever $0\leq
\lambda <1$ and so, $g\notin S_{0}\left( \tau \right) $.

\item[(c).] If $a\in L_{0}\left( \nu \right) $, then it is clear that $%
a\otimes \mathbf{1}\in LS\left( \mathcal{M}\right) $. We claim that $%
a\otimes \mathbf{1}\in S\left( \mathcal{M}\right) $ only if $a\in L_{\infty
}\left( \nu \right) $. Indeed, if $a\otimes \mathbf{1}\in S\left( \mathcal{M}%
\right) $, then there exists a projection $p\in P\left( \mathcal{M}\right) $
such that $\left( a\otimes \mathbf{1}\right) p\in L_{\infty }^{wo}\left( \nu
;B\left( H\right) \right) $ and $p^{\bot }$ is finite (that is, $p\left(
s\right) ^{\bot }$ is finite rank $\nu $-a.e.). Define the $\nu $-measurable
set $A\subseteq S$ by setting 
\begin{equation*}
A=\left\{ s\in S:\left\vert a\left( s\right) \right\vert >\left\Vert \left(
a\otimes \mathbf{1}\right) p\right\Vert _{\infty }\right\} .
\end{equation*}%
Since $\left\vert a\left( s\right) \right\vert \left\Vert p\left( s\right)
\right\Vert _{B\left( H\right) }\leq \left\Vert \left( a\otimes \mathbf{1}%
\right) p\right\Vert _{\infty }$ $\nu $-a.e. on $S$, it follows that $%
p\left( s\right) =0$ $\nu $-a.e. on $A$, that is, $p\left( s\right) ^{\bot }=%
\mathbf{1}$ $\nu $-a.e. on $A$. Consequently, $\nu \left( S\right) =0$ and
hence, $a\in L_{\infty }\left( \nu \right) $. It is now clear that $S\left( 
\mathcal{M}\right) \varsubsetneq LS\left( \mathcal{M}\right) $.
\end{enumerate}
\end{remark}

\begin{remark}
As follows from Corollary \ref{PCor03}, any derivation $\delta $ in $%
LS\left( \mathcal{M}\right) $, $S\left( \mathcal{M}\right) $ or $S\left(
\tau \right) $, where $\mathcal{M}=L_{\infty }\left( \nu \right) \overline{%
\otimes }B\left( H\right) $, is inner and hence, $\delta $ is $\mathcal{Z}$%
-linear. The main ingredient to obtain this $\mathcal{Z}$-linearity has been
Proposition \ref{PProp03} (and Remark \ref{PRem01}). In the setting of von
Neumann algebra, the $\mathcal{Z}$-linearity of derivations also follows
from the following proposition, which is of interest in its own right. The
proof of this proposition follows the same lines as the proof of Proposition %
\ref{PProp03}, although there are some non-trivial other ingredients.
However, in view of the length of the present paper, a detailed proof of
this proposition will appear elsewhere.

\begin{proposition}
Suppose that $\mathcal{M}$ is properly infinite, that is, every non-zero
projection in $P\left( \mathcal{Z}\left( \mathcal{M}\right) \right) $ is
infinite (with respect to $\mathcal{M}$). Let $\mathcal{A}\subseteq LS\left( 
\mathcal{M}\right) $ be a $\ast $-subalgebra such that $\mathcal{M}\subseteq 
\mathcal{A}$. If $\delta :\mathcal{A}\rightarrow \mathcal{A}$ is a
derivation, then $\delta $ is $\mathcal{Z}$-linear.
\end{proposition}

\noindent A typical example of a von Neumann algebra $\mathcal{M}$ which
satisfies the conditions of the previous proposition is $\mathcal{M}%
=L_{\infty }\left( \nu \right) \overline{\otimes }\mathcal{N}$, where $%
\mathcal{N}$ is a factor for which $\mathbf{1}_{\mathcal{N}}$ is infinite.
Indeed, if $0\neq p\in \mathcal{Z}\left( \mathcal{M}\right) $, then $p=\chi
_{A}\otimes \mathbf{1}_{\mathcal{N}}$ for some projection $\chi _{A}\in
L_{\infty }\left( \nu \right) $. Since $\mathbf{1}_{\mathcal{N}}$ is
infinite, there exists $q\in P\left( \mathcal{N}\right) $ such that $q\neq 
\mathbf{1}_{\mathcal{N}}$ and $q\sim \mathbf{1}_{\mathcal{N}}$. This implies
that $\chi _{A}\otimes q<p$ and $\chi _{A}\otimes q\sim p$. Hence, $p$ is an
infinite projection in $\mathcal{M}$.
\end{remark}

\noindent \texttt{A.F. Ber}

\noindent \texttt{Chief software developer}

\noindent \texttt{ISV "Solutions"}

\noindent \texttt{Tashkent, Uzbekistan}

\noindent \texttt{e-mail: ber@ucd.uz\bigskip }

\noindent \texttt{B. de Pagter}

\noindent \texttt{Delft Institute of Applied Mathematics}

\noindent \texttt{Faculty EEMCS, Delft University of Technology}

\noindent \texttt{P.O. Box 5031, 2600 GA Delft, The Netherlands}

\noindent \texttt{e-mail: b.depagter@tudelft.nl\bigskip }

\noindent \texttt{F.A. Sukochev}

\noindent \texttt{School of Mathematics and Statistics}

\noindent \texttt{University of New South Wales}

\noindent \texttt{Kensington, NSW 2052, Australia}

\noindent \texttt{e-mail: f.sukochev@unsw.edu.au}

\end{document}